\documentclass[reqno, 11pt]{amsart}
\usepackage{amssymb,amsmath,amsthm}
\usepackage{enumerate,float,hyperref,graphicx}
\hypersetup{%
  colorlinks=true,
  citecolor=blue
}
\usepackage{constants}

\allowdisplaybreaks
\numberwithin{equation}{section}
\newtheorem{theorem}{Theorem}[section]
\newtheorem{proposition}[theorem]{Proposition}
\newtheorem{lemma}[theorem]{Lemma}
\newtheorem{corollary}[theorem]{Corollary}
\newtheorem{claim}[theorem]{Claim}

\theoremstyle{definition}
\newtheorem{definition}{Definition}[section]
\newtheorem{assumption}{Assumption}[section]
\newtheorem{notation}{Notation}[section]
\theoremstyle{remark}
\newtheorem{remark}{Remark}[section]

\DeclareMathOperator{\Var}{Var}
\DeclareMathOperator{\Cov}{Cov}
\DeclareMathOperator{\Cor}{Corr}
\DeclareMathOperator{\Exp}{\mathbb{E}}
\DeclareMathOperator{\Prob}{\mathbb{P}}
\DeclareMathOperator{\Diam}{Diam}

\newcommand{\RR}{\mathbb{R}}
\newcommand{\ZZ}{\mathbb{Z}}
\newcommand{\QQ}{\mathbb{Q}}
\newcommand{\Lp}{\left(}
\newcommand{\Rp}{\right)}
\newcommand{\LP}{\left\{}
\newcommand{\RP}{\right\}}
\newcommand{\LT}{\left[}
\newcommand{\RT}{\right]}

\newcommand\numberthis{\addtocounter{equation}
    {1}\tag{\theequation}}
\newcommand{\ltwo}[1]{\lvert #1\rvert}
\newcommand{\wandering}[4]{
    \mathcal{W}\left(#1,#2,#3,#4\right)}

\newcommand{\event}{\mathcal{E}}
\newcommand{\segment}{\mathcal{I}}
\newcommand{\diff}{\mathcal{D}}
\newcommand{\edge}{\mathtt{E}}

\newcommand{\point}[1]{\boldsymbol #1}
\newcommand{\Unit}[1]{\point{e}_{#1}}
\newcommand{\bu}{\point{u}}
\newcommand{\bv}{\point{v}}
\newcommand{\origin}{\point{0}}
\newcommand{\bx}{\point{x}}
\newcommand{\by}{\point{y}}
\newcommand{\bz}{\point{z}}
\newcommand{\ba}{\point{a}}
\newcommand{\bb}{\point{b}}
\newcommand{\bw}{\point{w}}
\newcommand{\FF}{\mathcal{F}}
\newcommand{\xl}{\bx_{\sf L}}
\newcommand{\xs}{\bx_{\sf S}}

\newcommand{\Sector}{\mathcal{S}}
\newcommand{\tht}{{\theta^t}}
\newcommand{\uth}{{\Unit{\theta}}}
\newcommand{\utht}{{\Unit{\tht}}}
\newcommand{\thn}{{\theta_0}}
\newcommand{\thnt}{{\theta_0^t}}
\newcommand{\uthn}{{\Unit{\thn}}}
\newcommand{\uthnt}{{\Unit{\thnt}}}
\newcommand{\Bb}{\mathcal{B}}

\newcommand{\lowconst}{\mathfrak{p}}
\newcommand{\lowconstinv}{\lowconst^{-1}}
\newcommand{\upconst}{\mathfrak{q}}

\newcommand{\xprojthnot}{\pi^1_{\theta_0,\theta_0^t}}
\newcommand{\xprojth}{\pi^1_{\theta,\theta^t}} 
 
\newcommand{\yprojthnot}{\pi^2_{\theta_0,\theta_0^t}}
\newcommand{\yprojth}{\pi^2_{\theta,\theta^t}} 
 
\newcommand{\gwq}
    {\Delta(\lambda^{\zeta q} k)
        (\log(\lambda^{\zeta q} k))^{\eta/2}}
\newcommand{\gwp}
    {\Delta(\lambda^{\zeta p} k)
        (\log(\lambda^{\zeta p} k))^{\eta/2}}
\newcommand{\gwm}
    {\Delta(\lambda^{\zeta m} k)
        (\log(\lambda^{\zeta m} k))^{\eta/2}}
\newcommand{\gwpp}
    {\Delta(\lambda^{\zeta (p+1)} k)
        (\log(\lambda^{\zeta (p+1)} k))^{\eta/2}}

\newconstantfamily{epsilon}{
symbol=\epsilon,
}

\newconstantfamily{nu}{
symbol=\nu,
}

\begin{document}

\title[Fluctuations of Transverse Increments in
FPP]{Fluctuations of Transverse Increments in Two-dimensional
First Passage Percolation}
\date{}
\author{Ujan Gangopadhyay}
\address{Ujan Gangopadhyay, \ Department of Mathematics, \ University of Southern California, \ Los Angeles, CA, USA.}
\email{ujangangopadhyay@gmail.com}

\begin{abstract} 
We consider a model of first passage percolation (FPP) where
the nearest-neighbor edges of the standard two-dimensional
Euclidean lattice are equipped with random variables. These
variables are i.i.d.\, nonnegative, continuous, and have a
finite moment generating function in a neighborhood of $0$.
We derive consequences about transverse increments of passage
times, assuming the model satisfies certain properties. 
Approximately, the assumed properties are the following: We
assume that the standard deviation of the passage time on
scale $r$ is of some order $\sigma(r)$, and
$\left\{\sigma(r), r > 0\right\}$ grows approximately as a
power of $r$. Also, the tails of the passage time
distributions for distance $r$ satisfy an exponential bound
on a scale $\sigma(r)$ uniformly over $r$. In addition, the
boundary of the limit shape in a neighborhood of some fixed
direction $\theta$ has a uniform quadratic curvature. By
transverse increment we mean the difference of passage times
from the origin to a pair of points which are located as
follows: they are approximately in the same direction, say
$\theta$, from the origin; the direction of one of them from
the other is the direction of the tangent of the boundary of
the limit shape at the point on the limit shape in the
direction $\theta$. The main consequence derived is the
following. If $\sigma(r)$ varies as $r^\chi$ for some
$\chi>0$, and $\xi$ is such that $\chi=2\xi-1$, then the
fluctuation of the transverse increment of passage time
between a pair of points situated at distance $r$ from each
other is of the order of $r^{\chi/\xi}$.
\end{abstract}

\subjclass[2010]{Primary 60K35; Secondary 82B43}

\keywords{first passage percolation, increments}

\maketitle

\tableofcontents

\section{Introduction}

In this paper, we investigate the transverse increments of
passage times in the classical model of first passage
percolation (FPP) on $\ZZ^2$, which was introduced in \cite{HammersleyWelsh65}. 

\subsection{A brief description of the model}

Let $\edge(\ZZ^2)$ be the set of nearest-neighbor edges in
$\ZZ^2$. On $\edge(\ZZ^2)$ we consider a collection of random
variables $\mathbb{T}:=\LP\tau_e:e\in \edge(\ZZ^2)\RP$, which
are called \emph{edge-weights.} We assume certain properties
of the edge-weights. We categorize these assumptions as being
basic or technical. The basic we use throughout the paper,
and the technical we use more selectively.

\paragraph{\textbf{The basic assumptions:}} We assume that
the edge-weights are i.i.d., nonnegative, and continuous. In 
addition, there exists $C>0$ such that $\Exp\exp\Lp C 
\tau_e\Rp<\infty$. 

Using the edge-weights, we define the \emph{passage time} of 
a self-avoiding lattice path $\gamma$, denoted by
$T(\gamma)$, as the sum of the edge-weights of all the edges 
on the path $\gamma$, i.e.,  
\[
T(\gamma):=\sum_{\gamma\mbox{ contains } e} \tau_e.
\]
In the above definition, we adopt the convention that a path 
is a continuous curve in $\RR^2$. Next, we define the passage
time between two points $\bu$ and $\bv$ in $\ZZ^2$ as
\[
T(\bu,\bv):=\inf\LP T(\gamma)\mid\gamma\;\mbox{is a path joining $\bu$ and $\bv$}\RP.
\]
It follows from the above definition that $T$ is a random
pseudo-metric on $\ZZ^2$. It was shown in 
\cite{WiermanReh1978}, under only the i.i.d. and 
nonnegative assumptions on the edge-weights, the infimum in 
the definition of $T$ is
attained for some paths, i.e., the infimum is a minimum.
Since we have assumed that the edge-weights are also
continuous, it follows that there is only one such minimizing
path almost surely. We call this path \emph{the geodesic}
between $\bu$ and $\bv$, and denote it by $\Gamma(\bu,\bv)$.
Since the edge-weights have finite expectation, the passage
times $T(\bu,\bv)$, for all $\bu$, $\bv$, also have finite
expectation. Therefore,  
\[
h(\bu):=\Exp T(\origin,\bu)
\] 
is well-defined. From the triangle inequality of $T$ it
follows that $h$ is subadditive, i.e., for any 
$\bu,\bv\in\ZZ^2$ we have
\[
h(\bu+\bv)\leq h(\bu)+h(\bv).
\]
The subadditive ergodic theorem in \cite{Kingman1984} implies
that for any $\bu\in\ZZ^2$ the following limits exist almost
surely and in $L^1$:
\[
g(\bu):=\lim_{n\to\infty}\frac{T(\origin,n\bu)}{n}
       =\lim_{n\to\infty}\frac{h(n\bu)}{n}
       =\inf_{n>0}\frac{h(n\bu)}{n}.
\]
The domain of $g$ can be extended to $\QQ^2$ by taking limit 
along appropriate subsequences in the above definition. By 
extending the domain in this way, $g$ becomes a norm on 
$\QQ^2$. Therefore, the domain of $g$ can be further extended
to $\RR^2$. The unit ball in the norm $g$ is 
\[
\Bb:=\LP\bx\in\RR^2:g(\bx)\leq 1\RP,
\]
This is called the \emph{limit shape}. The \emph{wet region} 
at time $t$ is defined as 
\[
\Bb(t):=\LP\bx\in\ZZ^2:T(\origin,\bx)\leq t\RP.
\]
The shape theorem in \cite{CoxDur81} implies, under 
conditions milder than our basic assumptions, $\Bb(t)$ 
approaches $\Bb$ in an appropriate sense as $t\to\infty$. In 
addition, $\Bb$ is compact, convex, has a nonempty interior, 
and has all the symmetries of the lattice. 

\begin{notation}\label{notn:nonlattice}
    Here, we define passage times between points in $\RR^2$.
    For $\bx\in\RR^2$, let $\lfloor\bx\rfloor$ be the
    down-left corner of the unit square containing $\bx$ in
    $\ZZ^2$. For $\bx,\by\in\RR^2$, let 
    \[
    T(\bx,\by):=T(\lfloor\bx\rfloor,\lfloor\by\rfloor).
    \]
    Similarly, by $\Gamma(\bx,\by)$ we mean the geodesic
    $\Gamma(\lfloor\bx\rfloor,\lfloor\by\rfloor)$.
    Furthermore, for $\bx\in\RR^2$, let 
    \[
    h(\bx):=h(\lfloor\bx\rfloor).
    \]
\end{notation}

\begin{remark} 
    Throughout the paper, we denote by $C, C_0, C_1, C_2, \dots$ constants that depend only on the distribution of the edge-weights. We restart numbering of $C_i$s in each proof. Often we break the proof of a theorem in propositions and claims. In these situations, we do not restart numbering the constants in the proof of the propositions and claims. Also, when we use a result which has been proved before, we use ``tilde-version" of the variables and the parameters.
\end{remark} 

\begin{remark}\label{remark:nonlattice1}
    Extending definition of $T$ from $\ZZ^2$ to $\RR^2$ yields a minor technical issue. Although $g(\bx)\leq h(\bx)$ for all $\bx\in\ZZ^2$, this may not be true for $\bx\in\RR^2$. Instead, we have, for some constant $C_1>0$ and for all $\bx\in\RR^2$
    \[
    g(\bx)-C_1\leq h(\bx).
    \]
    Similarly, we have, for some constant $C_2>0$ and for all $\bx,\by\in\RR^2$ 
    \[
    h(\bx+\by)-C_2\leq h(\bx)+h(\by).
    \]
\end{remark}

\subsection{Heuristics of the main results}

It is common in the literature, for instance, in the works 
\cite{Newman95}, \cite{NewmanPiza1995}, 
\cite{LiceaNewman1996}, \cite{DamronHanson14}, 
\cite{DamronHanson2017}, to assume specific unproved 
properties of the limit shape. Often properties such as 
differentiability, and curvature, either locally or globally,
which eliminates the possibility of facets or corners. These
properties are believed to be valid under our assumptions,
but there is no proof yet. We also make similar assumptions.

Suppose the boundary of the limit shape is differentiable at 
a direction $\theta$, and $\theta^t$ is the corresponding 
tangential direction. By \emph{transverse increments} we mean
differences of the form $T(\origin,\bx)-T(\origin,\by)$ where
$\bx$ has direction $\theta$, $\bx-\by$ has direction $\tht$.
Heuristically we can say what the order of the fluctuations
of transverse increments should be. For this, we need the
scaling exponents $\chi$ and $\xi$.

It is believed that for FPP on Euclidean lattices of any
dimension, that there exists an exponent $\chi$, called the
\emph{`fluctuation exponent,'} such that
$T(\origin,\bv)-h(\bv)$ is of the order of $\ltwo{\bv}^\chi$
($\ltwo{\cdot}$ is the Euclidean norm.) Also, it is believed
that there exists an exponent $\xi$, called the
\emph{`wandering exponent,'} such that the geodesic
$\Gamma(\origin,\bv)$ wanders $\ltwo{\bv}^\xi$ distance on
average from the line joining $\origin$ and $\bv$. The two
exponents are related by the equation $\chi=2\xi-1$ which was
proved in \cite{Chatterjee2013} assuming $\chi$ and
$\xi$ exist in a certain sense. In dimension $d=2$ it is
believed that $\chi=1/3$ and $\xi=2/3$. In $d=3$ it is
believed that $\chi$ is approximately $1/4$, and in higher
dimensions there does not seem to be a consensus even among
physicists about values of $\chi$ and $\xi$, see for example
\cite{Marinari_2002}, \cite{Le_Doussal_2005},
\cite{Fogedby_2008}, \cite{Kim_2013}, \cite{Alves_2014}. In the exactly solvable models of two-dimensional last passage percolation, it has been proved that $\chi=1/3$ and $\xi=2/3$, see \cite{KurtShape}, \cite{KurtTransverse}, \cite{BalazsCatorTimo2006}.

\begin{figure}[H]
    \centering
    \includegraphics[width=0.65\linewidth]{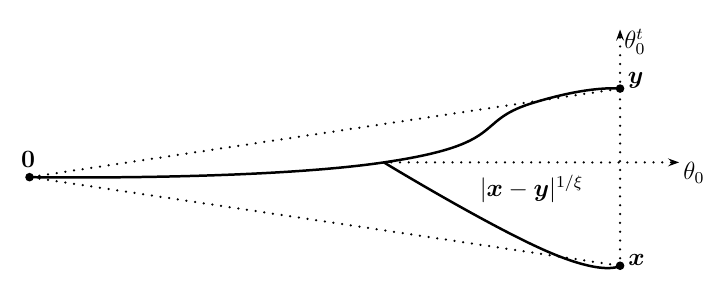}
    \caption{Illustration for the heuristic of the exponent
    $\chi/\xi$: directions of $\bx$ and $\by$ are
    approximately $\thn$; direction of $\by-\bx$ is $\thnt$;
    the two geodesics $\Gamma(\origin,\bx)$ and
    $\Gamma(\origin,\by)$ are expected to coalesce
    approximately at distance $\ltwo{\bx-\by}^{1/\xi}$ in
    $-\thn$ direction when traced starting from $\bx$ and
    $\by$ respectively.}
    \label{Fig1}
\end{figure}

If one assumes the existence of these exponents in some 
appropriate sense, then fluctuations of the transverse 
increment $T(\origin,\bx)-T(\origin,\by)$ should be of the 
order of $\ltwo{\bx-\by}^{\chi/\xi}$. The heuristic of this 
is the following, see Figure~\ref{Fig1}. We expect that the
geodesics $\Gamma(\origin,\bx)$ and $\Gamma(\origin,\by)$
stay disjoint after starting from $\bx$ and $\by$
respectively for a distance of the order of
$\ltwo{\bx-\by}^{1/\xi}$. Then these two branches should
contribute approximately independently 
$\ltwo{\bx-\by}^{\chi/\xi}$ to the fluctuation. The right
scale of the coalescence time as above has been proved in
\cite{BasuSarkarSly2019} for the exactly solvable model of
two-dimensional last passage percolation.

One reason for studying the fluctuations of transverse 
increments is the following. In $d=2$, it is believed that 
the transverse increments behave like increments of Brownian 
motion, that is, the increments are approximately 
uncorrelated. If this is true, then the exponent for 
fluctuation of transverse increments should be $1/2$ so that 
$\chi/\xi=1/2$. This with $\chi=2\xi-1$ would imply 
$\chi=1/3$ and $\xi=2/3$. 

\begin{figure}[H]
    \centering
    \includegraphics[width=0.65\linewidth]{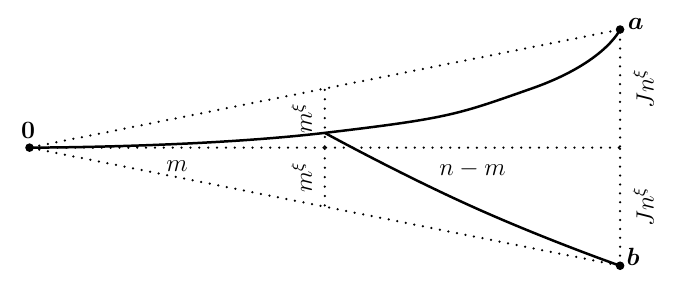}
    \caption{Illustration for the heuristic of the exponent 
    $2\chi/(1-\xi)$: directions of $\ba$ and $\bb$ are 
    approximately $\thn$; direction of $\ba-\bb$ is $\thnt$; 
    distances of $\ba$ and $\bb$ from origin are 
    approximately $n$; distance between $\ba$ and $\bb$ is 
    $2Jn^\xi$; the two geodesics $\Gamma(\origin,\ba)$ and 
    $\Gamma(\origin,\bb)$ are expected to branch apart at 
    approximately distance $m$ from the origin where 
    $m^\xi/m=Jn^\xi/n$.}
    \label{fig:longrangeheu}
\end{figure}%

As an application of the upper bound on fluctuations of 
transverse increments we get an upper bound on long-range 
correlations. By long-range correlation we mean the 
correlation between $T(\origin,\ba)$ and $T(\origin,\bb)$ 
where $\ba$ and $\bb$ are points approximately in the same 
direction from the origin and distance between $\ba$ and 
$\bb$ is large compared to typical wanderings of the 
geodesics $\Gamma(\origin,\ba)$ and $\Gamma(\origin,\bb)$. 
Heuristically we can say the following about this 
correlation, see Figure~\ref{fig:longrangeheu}. Assuming 
$\chi$ and $\xi$ exist, typical wandering of the geodesics 
$\Gamma(\origin,\ba)$ and $\Gamma(\origin,\bb)$ is of the 
order of $n^\xi$ where $n$ is the distance of the points 
$\ba$, $\bb$ from the origin. Suppose the distance between 
$\ba$ and $\bb$ is $J n^\xi$ for some large $J$. Then 
$\Gamma(\origin,\ba)$ and $\Gamma(\origin,\bb)$ are expected 
to branch apart at a distance $m$ from the origin such that 
the distance between the rays joining $\origin$ to $\ba$ and 
$\origin$ to $\bb$ at distance $m$ from the origin is of the 
order of the typical wandering of the geodesics of at 
distance $m$ from origin. So we have approximately 
$m^\xi/m=Jn^\xi/n$. Hence $m=n J^{-(1-\xi)^{-1}}$. Then the 
covariance between $T(\origin,\ba)$ and $T(\origin,\bb)$ is 
expected to be of the order of $m^{2\chi}$. So the 
correlation should be of the order of $J^{-2\chi/(1-\xi)}$. 

\subsection{Advanced assumptions}

Along with our basic assumptions we make the following 
assumptions. Similar assumptions have been used in 
\cite{KenGeoBiGeo}, \cite{KenUniform}, \cite{ganguly2020optimal}. We assume 
that there exists $\sigma:(0,\infty)\to(0,\infty)$ such that 
the following hold.

\resetconstant

\begin{assumption}\label{asm:sigmaexp}
    There exist positive constants $C_1$, $C_2$, $C_3$, such that for all $\bx,\by\in\RR^2$ with $\ltwo{\bx-\by}\geq C_1$, and all $t>0$, we have
    \begin{equation}\tag{A1}\label{A1}
    \Prob\Lp\,|T(\bx,\by)-h(\bx-\by)|\geq t\sigma(\ltwo{\bx-\by})\,\Rp
    \leq C_2 e^{-C_3 t}.
    \end{equation}
\end{assumption}

\begin{assumption}\label{asm:sigmareg}
    There exist constants $\lowconst>0$, $\upconst>0$, 
    $\alpha\in(0,1)$, $\beta\in(0,1)$, $\alpha\leq\beta$, 
    such that for all $x>y>0$ we have
    \begin{equation}\tag{A2}\label{A2}
        \lowconst\Lp\frac{x}{y}\Rp^\alpha
    \leq\frac{\sigma(x)}{\sigma(y)}
    \leq\upconst\Lp\frac{x}{y}\Rp^\beta.
    \end{equation}
\end{assumption}

\resetconstant

\begin{assumption}\label{asm:sigmavar}
    There exist positive constants $\epsilon$, $C$,
    such that for $\bx,\by\in\RR^2$ with
    $\ltwo{\bx-\by}\geq C$, we have
    \begin{equation}\tag{A3}\label{A3}
        \begin{aligned}
            \Prob\Lp\, 
            T(\bx,\by) \leq 
            h(\bx-\by)-\epsilon\sigma(\ltwo{\bx-\by})
            \,\Rp
            \geq \epsilon,\\
            \Prob\Lp\, 
            T(\bx,\by) \geq 
            h(\bx-\by)+\epsilon\sigma(\ltwo{\bx-\by})
            \,\Rp
            \geq \epsilon.
        \end{aligned}
    \end{equation}
\end{assumption}

\begin{remark} 
    By Remark~1.1 of \cite{KenGeoBiGeo}, we assume without 
    loss of generality that $\sigma$ is monotonically 
    increasing and continuous.
\end{remark} 

\resetconstant

\begin{remark}
    Under Assumption~\ref{asm:sigmaexp}, 
    Assumption~\ref{asm:sigmavar} is equivalent to saying 
    that there exist positive constants $C_1$, 
    $C_2$, such that for all $\bx,\by\in\RR^2$ 
    with $\ltwo{\bx-\by}\geq C_1$, we have 
    \[
        \Var(T(\bx,\by))\geq
        C_2\sigma^2(\ltwo{\bx-\by}).
    \]
\end{remark} 

\begin{remark} 
    The assumption $\beta<1$ is natural because, from results
    in \cite{Kesten1993}, the passage times are known to
    satisfy exponential concentration with scaling exponent
    $1/2$, which shows that $\chi$ must be $\leq 1/2$ under
    any reasonable definition.
\end{remark}  

\begin{remark} 
    The assumption $\alpha>0$ is natural because, for certain
    definition of $\chi$ and $\xi$ it was shown in
    \cite{WehrAizenman1990} that $\chi\geq (1-(d-1)\xi)/2$ in
    $d$-dimensions, which in two dimensions coupled with
    $\chi\geq 2\xi-1$ and $\chi\leq 1/2$ yields $\chi\geq
    1/8$.
\end{remark} 

\begin{notation} 
    For any direction $\theta$ the unit vector in direction 
    $\theta$ is denoted by $\uth$. By abuse of notation, we 
    denote the standard unit vectors in $\RR^2$ by $\Unit{1}$
    and $\Unit{2}$.
\end{notation}

\begin{definition} 
    We say that a direction $\thn$ is of \emph{type I} if 
    there exist constants $C>0$, $\delta_1>0$, $\delta_2>0$ 
    such that the following holds: the limit shape boundary 
    $\partial\Bb$ is differentiable in the sector 
    $(\thn-\delta_1,\thn+\delta_1)$; for 
    $|\delta|\leq\delta_2$ and 
    $\theta\in(\thn-\delta_1,\thn+\delta_1)$, we have
    \[
        g(\uth+\delta\utht)-g(\uth)\geq C\delta^2,
    \]
    where $\tht$ is the direction of the tangent to
    $\partial\Bb$ at the point in direction $\theta$.
\end{definition} 

\begin{remark} 
    We take the direction of the tangents in 
    counter-clockwise direction around the limit shape 
    boundary.
\end{remark}

\begin{remark} 
    An alternative formulation of type I direction is the 
    following: there exist constants $C>0$, $\delta_1>0$ such
    that $\partial\Bb$ is differentiable in the sector
    $(\thn-\delta_1,\thn+\delta_1)$; for
    $\theta\in(\thn-\delta_1,\thn+\delta_1)$ and all
    $\delta\in\RR$ we have
    \[
        g(\uth+\delta\utht)-g(\uth)
        \geq C\min\LP|\delta|,\delta^2\RP.
    \]
\end{remark}

\begin{definition} 
    We say that a direction $\thn$ is of \emph{type II} if 
    there exist constants $C>0$, $\delta_1>0$, $\delta_2>0$ 
    such that the following holds: the limit shape boundary 
    $\partial\Bb$ is differentiable in the  sector 
    $(\thn-\delta_1,\thn+\delta_1)$; for 
    $|\delta|\leq\delta_2$ and 
    $\theta\in(\thn-\delta_1,\thn+\delta_1)$ we have
    \begin{equation}\label{eq:defcurvup}
        g(\uth+\delta\utht)-g(\uth)\leq C\delta^2,
    \end{equation}
    where, as before, $\tht$ is the direction of the tangent 
    to $\partial\Bb$ at the point in direction $\theta$.
\end{definition}

\begin{remark} 
    In a neighborhood of a type I direction the limit shape boundary cannot have a facet. Similarly in a neighborhood of a type II direction the limit shape boundary cannot have a corner.
\end{remark}

\begin{remark}\label{remark:linear}
    As observed in Remark~1.2 of \cite{KenGeoBiGeo}, the 
    condition in \eqref{eq:defcurvup} can be alternatively 
    stated as follows. If $\bu_\theta$ is the point on 
    $\partial\Bb$ in direction $\theta$, then in a 
    neighborhood of $\bu_\theta$, the boundary is squeezed 
    between the tangent at $\bu_\theta$ and a parabola 
    tangent to $\partial\Bb$ at $\bu_\theta$. This implies 
    that the direction of the tangent grows at most linearly 
    in a neighborhood of $\theta$. So, if $\thn$ is a 
    direction of type II, then there exist constants $C>0$, 
    $\delta_1>0$, $\delta_2>0$ such that the following holds:
    $\partial\Bb$ is differentiable in the sector 
    $(\thn-\delta_1,\thn+\delta_1)$; for 
    $\theta_1,\theta_2\in(\thn-\delta_1,\thn+\delta_1)$ with 
    $|\theta_1-\theta_2|\leq\delta_2$, we have 
    $|\theta_1^t-\theta_2^t|\leq C|\theta_1-\theta_2|$. 
\end{remark} 

\subsection{Main results}

\begin{notation}
    Given two distinct directions $\theta_1$, $\theta_2$ 
    define the projections $\pi^1_{\theta_1,\theta_2}$ and 
    $\pi^2_{\theta_1,\theta_2}$ so that for any $\bv$ we have
    \[
        \bv
        =\pi^1_{\theta_1,\theta_2}(\bv)\Unit{\theta_1}
            +\pi^2_{\theta_1,\theta_2}(\bv)\Unit{\theta_2}.
    \]
\end{notation}

\begin{notation} 
    For $n>0$ let 
    \[
        \Delta(n):=\Lp n\sigma(n)\Rp^{1/2}.
    \]
\end{notation}

Our first main result is the following.

\begin{theorem}\label{thm:loglogupinc}
    Let $\thn$ be a direction of both type I and II. For
    $n>0$, $L>0$, define
    \[
        \segment(n,L):=
            \LP\,\bx\in\RR^2\,:\,
                \xprojthnot(\bx)=n,\,
                0\leq\yprojthnot(\bx)\leq L
            \,\RP,
    \]
    and 
    \[
        \diff(n,L):=
            \max\LP\,|T(\origin,\bx)-T(\origin,\by)|\,:\,
                \bx,\by\in\segment(n,L)
            \,\RP.
    \]
    Fix $\eta\in(0,1]$. Then, under the 
    Assumptions~\ref{asm:sigmaexp} and \ref{asm:sigmareg},
    there exist positive constants $C_1$, $C_2$, $L_0$, 
    $n_0$, $t_0$, such that for $L\geq L_0$, $n\geq n_0$, 
    $t\geq t_0$, $L\leq\Delta(n)$, we have 
    \[
        \Prob\Lp\,\diff(n,L)\geq 
            t(\log L)^\eta\sigma\Lp\Delta^{-1}(L)\Rp
        \,\Rp
    \leq C_1 \exp\Lp-C_2 t(\log L)^\eta\Rp.
    \]
\end{theorem} 

The following theorem is our lower bound on the fluctuations 
of transverse increments. In this theorem, we show that the 
standard deviation of the transverse increment between a pair
of points at a distance $L$ is at least of the order of
$\sigma\Lp\Delta^{-1}(L)\Rp$ with a correction factor smaller
than any power of $L$.

\begin{theorem}\label{thm:lowmain}
    Let $\thn$ be a direction of both type I and II. Fix 
    $\nu\in(1/2,1)$. Then, under the
    Assumptions~\ref{asm:sigmaexp}, \ref{asm:sigmareg}, and 
    \ref{asm:sigmavar}, there exist positive constants $L_0$,
    $n_0$, such that for $L\geq L_0$, $n\geq n_0$,
    $L\leq\Delta(n)$, we have
    \[
        \Var\Lp 
            T(\origin,n\uthn)-T(\origin,n\uthn+L\uthnt)
        \Rp
        \geq 
            e^{-(\log L)^{\nu}}
                \sigma^2\Lp\Delta^{-1}(L)\Rp.
    \]
\end{theorem} 

As a corollary of Theorems~\ref{thm:loglogupinc} and
\ref{thm:lowmain} we get the following result. It shows that
if we assume $\chi$ and $\xi$ exist in a certain sense, then
$\chi/\xi$ is the correct scaling exponent for the
fluctuations of the transverse increments. 

\begin{corollary} 
    Suppose there exists $\chi>0$ such that 
    \[
        \lim_{x\to\infty}
        \frac{\log\sigma(x)}{\log x}
        =\chi,
    \]
    and let 
    \[
        \xi:=
        \frac{1+\chi}{2}=
            \lim_{x\to\infty}
                \frac{\log\Delta(x)}{\log x}.
    \]
    Let $\thn$ be a direction of both type I and II. Then,
    under the Assumptions~\ref{asm:sigmaexp},
    \ref{asm:sigmareg}, and \ref{asm:sigmavar}, there exist
    functions $f_1$, $f_2$, $f_3$, which converge to $0$ at
    $\infty$, and positive constants $C_1$, $C_2$, $C_3$,
    $n_0$, $L_0$, $t_0$, such that for $n\geq n_0$, $L\geq
    L_0$, $t\geq t_0$, $L\leq n^{\xi+f_1(n)}$, we have
    \[
        \Prob\Lp\,
            |T(\origin,n\uthn)-T(\origin,n\uthn+L\uthnt)|
            \geq t L^{\chi/\xi + f_2 (L)}
            \,\Rp
        \leq C_1\exp\Lp - C_2 t\Rp,
    \]
    and
    \[
        \Var\Lp
            T(\origin,n\uthn)-T(\origin,n\uthn+L\uthnt)
            \Rp
        \geq C_3 L^{2\chi/\xi+f_3(L)}. 
    \]
\end{corollary}

\begin{proof}
    Fix $\eta\in(0,1)$ and $\nu\in(1/2,1)$. Define $f_1$,
    $f_2$, $f_3$ such that for all $x>1$,
    \begin{align*}
        & x^{\xi+f_1(x)} = \Delta(x),\\
        & x^{\chi/\xi+f_2(x)} = 
            \sigma(\Delta^{-1}(x))(\log x)^{\eta},\\
        & x^{2\chi/\xi+f_3(x)} = 
            e^{-(\log x)^{\nu_0}}\sigma^2(x).
    \end{align*}
    Then $f_1$, $f_2$, $f_3$ converge to $0$ at $\infty$ and
    the result readily follows from
    Theorems~\ref{thm:loglogupinc} and \ref{thm:lowmain}.
\end{proof} 

\begin{notation} 
    Let $f(n):=\Delta(n)(\log n)^{1/2}/n$ and 
        $f^{-1}(y):=\sup\LP x:f(x)\geq y\RP$.
\end{notation}

\begin{remark}
    Since $\beta<1$ by \eqref{A2}, and because we have 
    assumed $\sigma$ is monotonically increasing and 
    continuous, we get that $f$ is continuous and goes to $0$
    at $\infty$. Therefore, $f^{-1}$ is continuous,
    monotonically decreasing, and converges to $0$ at
    $\infty$.
\end{remark}

Now we state the result on upper bound of long-range
correlations. 

\begin{theorem}\label{thm:longrange}
    Let $\thn$ be a direction of both type I and II. 
    Recall $\beta$ and $\upconst$ from \eqref{A2}. Fix
    $\delta\in(0,(1-\beta)/2)$. Then, under the
    Assumptions~\ref{asm:sigmaexp} and \ref{asm:sigmareg},
    there exist positive constants $C$, $J_0$, $n_0$, such
    that for $n\geq n_0$, $J\in[\upconst^{1/2}J_0,n^\delta]$,
    we have
    \begin{multline*}
    \Cov\Lp\, 
        T(\origin,n\uthn-J\Delta(n)(\log n)^{1/2}\uthnt),
        T(\origin,n\uthn+J\Delta(n)(\log n)^{1/2}\uthnt)
    \,\Rp\\ 
    \leq 
    C\sigma^2\Lp f^{-1}\Lp\frac{J}{J_0}f(n)\Rp\Rp\log n.
    \end{multline*}
\end{theorem} 

The following corollary shows how we get the exponent 
$2\chi/(1-\xi)$ under further regularity assumptions on 
$\sigma$. 

\begin{corollary}\label{cor:longcor}
    Suppose $\sigma(n)=n^{\chi} L(n)$, where $L$ is a slowly 
    varying function. Let $\xi:=(1+\chi)/2$. Fix 
    $\delta_1\in(0,(1-\beta)/2)$. Let $\thn$ be a direction 
    of both type I and II. Then, under the 
    Assumptions~\ref{asm:sigmaexp} and \ref{asm:sigmareg}, 
    there exist positive constants $C$, $J_0$, such that the 
    following holds: for any $\delta_2>0$ there exists
    $n_0>0$ such that for $n\geq n_0$ and
    $J\in[\upconst^{1/2}J_0,n^{\delta_1}]$, we have
    \begin{multline*}
    \Cor\Lp\, 
        T(\origin,n\uthn-J\Delta(n)(\log n)^{1/2}\uthnt), 
        T(\origin,n\uthn+J\Delta(n)(\log n)^{1/2}\uthnt)
    \,\Rp\\ 
    \leq C J^{-2\chi/(1-\xi)+\delta_2}\log n.
    \end{multline*}
\end{corollary} 

\resetconstant 

\begin{proof}
    From Theorem~\ref{thm:longrange} we get positive
    constants $\Cl{longcor1}$, $J_0$, $n_0$, such that for
    $n\geq n_0$, $J\in[\upconst^{1/2}J_0,n^{\delta_1}]$, we
    have 
    \begin{equation}\label{eq:longrangecor0}
        \Cov\Lp\,
            T(\origin,\ba),T(\origin,\bb)
            \,\Rp
        \leq\Cr{longcor1}\sigma^2(m)\log n,
    \end{equation}
    where $\ba:=n\uthn+J\Delta(n)(\log n)^{1/2}\uthnt$, 
    $\bb:=n\uthn-J\Delta(n)(\log n)^{1/2}\uthnt$, 
    $m:=f^{-1}(Jf(n)/J_0)$. Using $J\leq n^{\delta_1}$, 
    $\delta_1<(1-\beta)/2$, and \eqref{A2}, we get 
    $\ltwo{\ba}\leq\Cl{longcor3}n$, 
    $\ltwo{\bb}\leq\Cr{longcor3}$. Hence, using 
    \eqref{eq:longrangecor0} and \eqref{A2}, we get
    \begin{equation}\label{eq:longrangecor1}
        \Cor\Lp\, 
            T(\origin,\ba),T(\origin,\bb)
            \,\Rp
    \leq\C\frac{\sigma^2(m)}{\sigma^2(n)}\log n.
    \end{equation}
    From $m=f^{-1}(Jf(n)/J_0)$ we get  
    \[
        J_0\frac{\Delta(m)(\log m)^{1/2}}{m} =J\frac{\Delta(n)(\log n)^{1/2}}{n}.
    \]
    Therefore, using $J\in[\upconst^{1/2}J_0,n^\delta_1]$, 
    $\delta_1<(1-\beta)/2$, and \eqref{A2}, we get $m\leq n$ 
    and $\log m\geq\Cl{longcor2}\log n$. Fix an $\delta_2>0$ 
    and let $\delta_3>0$ be such that 
    \begin{equation}\label{eq:longrangecor3}
            \frac{2\chi}{1-\xi}-\delta_2
        \leq\frac{2\chi-2\delta_3}{1-\xi+\delta_3/2}.
    \end{equation}
    Since $L$ is slowly-varying, by possibly increasing $n_0$
    based on $\delta_2$, we get
    \begin{equation}\label{eq:longrangecor2}
        \frac{L(n)}{L(m)}\geq\Lp\frac{n}{m}\Rp^{-\delta_3}.
    \end{equation}
    Therefore, using $\Delta(n)(\log n)^{1/2}=n^\xi L(n)^{1/2}(\log n)^{1/2}$ we get  
    \[
            \Lp\frac{n}{m}\Rp^{1-\xi}
        =   \frac{J}{J_0}\Lp\frac{L(n)\log n}
                {L(m)\log m}\Rp^{1/2}
        \geq\frac{J}{J_0}\Lp\frac{n}{m}\Rp^{-\delta_3/2}.
    \]
    Combining this with \eqref{eq:longrangecor2} and \eqref{eq:longrangecor3}, we get 
    \[
        \frac{\sigma^2(m)}{\sigma^2(n)}
        =\Lp\frac{m}{n}\Rp^{2\chi}\Lp\frac{L(m)}{L(n)}\Rp^2
        \leq\Lp\frac{m}{n}\Rp^{2\chi-2\delta_3}
        \leq\Lp\frac{J_0}{J}\Rp^{
            \frac{2\chi-2\delta_3}{1-\xi+\delta_3/2}}
        \leq J_0^{\frac{2\chi}{1-\xi}}J^{
            -\frac{2\chi}{1-\xi}+\delta_2}. 
        \]
    Combining this with \eqref{eq:longrangecor1} completes
    the proof of Corollary~\ref{cor:longcor}.
\end{proof}

\section{Wandering of geodesics}

\resetconstant

In this section we establish some upper bounds on the wandering of geodesics. 

\begin{notation} 
    For any set $A\subset\RR^2$ let 
    \[
        \Diam(A):=\sup\LP\,\ltwo{\bx-\by}:\bx,\by\in A\,\RP. 
    \]
\end{notation}

The following result provides a preliminary bound on the 
wandering of geodesics. The proof of this result follows from
Proposition~5.8 of \cite{Kesten86} under our basic
assumptions. The result has been shown to hold under milder
assumptions in Theorem~6.2 of 
\cite{AuffingerDamronHanson2015}.

\begin{lemma}\label{lem:boxwand}
    There exist positive constants $\Cl{bw1}$, $\Cl{bw2}$,
    $\Cl{bw3}$, such that the following holds. If
    $\ltwo{\bu-\bv}\geq\Cr{bw1}$ for some $\bu,\bv\in\RR^2$,
    then
    \[ 
    \Prob\Lp\,
        \Diam(\Gamma(\bu,\bv))\geq\Cr{bw2}\ltwo{\bu-\bv}
    \,\Rp\leq e^{-\Cr{bw3}\ltwo{\bu-\bv}}.
    \]
\end{lemma}

Utilizing the curvature of the limit shape we get a more
refined bound on the wandering of the geodesics. The
curvature of the limit shape is used in the following manner.
Consider two points in $\RR^2$. The shortest path between
them in the $g$-norm is, of course, the line joining them.
When the geodesic between them wander transversely too far
from the line joining them, the extra distance covered by the
geodesic in the $g$-norm can be thought of as a cost for
excessive wandering. A lower bound of this cost yields an
upper bound on the wandering of the geodesics. A lower bound
of the curvature of the limit shape provides a lower bound on
the cost associated with the $g$-norm. This is stated in the
following result, which is essentially same as Lemma~2.3 of
\cite{KenGeoBiGeo}. We state it without proof.  

\begin{lemma}\label{lem:auxgeom1}
    Let $\thn$ be a direction of type I. Then there exist positive constants $C$ and $\delta$ such that for $n>0$, $k$, $l$, $d$, satisfying $|l|/n\leq\delta$, we have     
    \[
    g(\bu) + g(\ba-\bu)-g(\ba)\geq C\min\LP |d|,\frac{d^2}{n}\RP,
    \]
    where $\ba:=n\uthn+l\uthnt$ and $\bu:=k\uthn+(l\frac{k}{n}+d)\uthnt$.
\end{lemma}

Geodesics cannot wander too much because the cost associated
with the $g$-norm becomes difficult to be compensated by the
fluctuations of passage times. Thus, bounds on the
fluctuations $T(\origin,\bx)-g(\bx)$, combined with
Lemma~\ref{lem:auxgeom1} yields upper bounds on geodesic
wanderings. By Assumption~\ref{asm:sigmaexp} we know
$T(\origin,\bx)-h(\bx)$ satisfies exponential concentration
in the scale $\sigma(\ltwo{\bx})$ uniformly over $\bx$. So we
need a bound on the differences $h(\bx)-g(\bx)$. These
differences are known as nonrandom fluctuations in the
literature. A general method of bounding the nonrandom
fluctuations was developed in \cite{Ken93,Ken97}.
There it was shown, using exponential concentration of
$T(\origin,\bx)-h(\bx)$ on the scale of $\ltwo{\bx}^{1/2}$
from \cite{Kesten1993}, that $h(\bx)-g(\bx)$ is at most of
the order of $\ltwo{\bx}^{1/2}\log\ltwo{\bx}$. In our case,
Alexander's method can be used with no significant alteration
to yield a bound of the order of 
$\sigma(\ltwo{\bx})\log\ltwo{\bx}$, which we state below.
Improvements to the logarithmic bound has been made in
\cite{Tessera2018}, \cite{DamronWang} and \cite{KenUniform} 
for some related 
models, which we briefly discuss in 
Section~\ref{sec:ghloglog}. We
also improve the logarithmic bound to 
$\sigma(\ltwo{\bx})(\log\ltwo{\bx})^\eta$ for arbitrary small
$\eta>0$ in Section~\ref{sec:ghloglog}. This improvement is
necessary to prove the lower bound result
Theorem~\ref{thm:lowmain}. To state the bound on the
nonrandom fluctuations we use the notion of `general
approximation property' from \cite{Ken97}.

\begin{notation}\label{notn:functionclass}
    Let $\Phi$ be the set of functions from $(0,\infty)$ to 
    $[0,\infty)$ such that for every $\phi$ there exists some
    $C\geq 0$ such that $\phi(x)>0$ for $x>C$ and
    $\inf_{x>y>C}\phi(x)/\phi(y)>0$. 
    For $\eta\in(0,1]$, define $\phi_\eta(k)=0$ for $k\leq 1$, and for $k>1$ 
    \[
    \phi_\eta(k):=k^{-\alpha}\sigma(k)(\log k)^\eta.
    \]
    Also define $\widehat{\phi}(k)=0$ for $k\leq 2$, and for 
    $k>2$
    \[
    \widehat{\phi}(k):=k^{-\alpha}\sigma(k)\log\log k.
    \]
    So $\widehat{\phi}$ and $\phi_\eta$ for all
    $\eta\in(0,1]$ belong to $\Phi$.
\end{notation}

\begin{definition}
    For $\nu\geq 0$ and $\phi\in\Phi$ we say that $h$
    satisfies the \emph{general approximation property} 
    with exponent $\nu$ and correction factor $\phi$ in a
    sector of directions $\Sector$ if there exist
    positive constants $C$ and $M$ such that for all
    $\bx\in\RR^2$ with $\ltwo{\bx}\geq M$ and having
    direction in $\Sector$ we have 
    \[
    h(\bx)\leq g(\bx)+C\ltwo{\bx}^\nu\phi(\ltwo{\bx}).
    \]  
    When we want to specify the relevant constants, we
    say $h$ satisfies $\text{GAP}(\nu,\phi,M,C)$ in
    sector $\Sector$. 
\end{definition}

In \cite{Ken97}, the class of correction factors consisted of
non-decreasing functions and the general approximation
property was not restricted to any particular set of
directions. These are some minor modifications we need in our
setup. As we mentioned before, we get the following result in
our context by following Alexander's method. 

\begin{proposition}\label{prop:nrflog}
    Under the Assumptions~\ref{asm:sigmaexp} and 
    \ref{asm:sigmareg}, there exist positive constants $C$, 
    $M$, such that $h$ satisfies 
    $\text{GAP}(\alpha,\phi_1,M,C)$ in all directions, i.e., 
    for all $\ltwo{\bx}\geq M$, we have 
    \[
    h(\bx) \leq g(\bx) + C\sigma(\ltwo{\bx})\log\ltwo{\bx}.
    \]
\end{proposition}

Let us now introduce a notation to measure wandering of 
geodesics. 

\begin{notation}\label{notn:wandering}
    Suppose $\thn$ is a direction where the boundary of the 
    limit shape is differentiable. Let $\thnt$ be the 
    direction of the tangent. Let 
    $\wandering{\bu}{\bv}{k}{\thn}$ denote the maximum 
    wandering of the geodesic $\Gamma(\bu,\bv)$, in 
    $\pm\thnt$ directions, from the line joining $\bu$ and 
    $\bv$, when the geodesic is at a distance $k$ from $\bu$ 
    in $\thn$ direction. We allow $k$ to be negative also. 
    More precisely  
    \begin{multline*}
        \wandering{\bu}{\bv}{k}{\thn}:=
        \max\left\{\,
            \left|
            \yprojthnot(\bw-\bu)
            -k\frac{\yprojthnot(\bv-\bu)}
            {\xprojthnot(\bv-\bu)}
        \right|\,:\,
            \bw\in\Gamma(\bu,\bv),\right.\\ \left.\xprojthnot(\bw-\bu)=k\,
            \vphantom{\left|\yprojthnot(\bw-\bu)
            -k\frac{\yprojthnot(\bv-\bu)}
            {\xprojthnot(\bv-\bu)}
            \right|}
        \right\}. 
    \end{multline*}
\end{notation}

In the next result, we show if distance between the endpoints
is approximately $n$, then wandering at a fixed but arbitrary
distance $k$ from one of the end points is at most of the
order of $\Delta(n)$ with some logarithmic correction
factors. Thus, roughly speaking, this theorem deals with
global wandering of a geodesic. The tail bound in this result
is sub-optimal for $k$ bigger than a large multiple of $n$,
because in that case we get a better bound using
Lemma~\ref{lem:boxwand}. In Theorem~\ref{thm:endwandlog} we
consider local wandering of geodesics. We show that at a
distance $k$ from an endpoint the wandering is of the order
of $\Delta(k)$ with some logarithmic factors. In Theorem~3 of
\cite{BasuSarkarSly2019} and Theorem~4.4 of \cite{Basu_2019}
similar bound has been proved for geodesics in the integrable
model of last passage percolation. Although, there the
results are sharper, i.e., there are no logarithmic
correction factors involved, because in the integrable model
of last passage percolation exact asymptotics of the
distribution of the passage times are known.

\begin{theorem}\label{thm:midptwand}
    Let $\thn$ be a direction of type I. Suppose $h$
    satisfies GAP with exponent $\alpha$ and correction
    factor $\phi_\eta$ for some $\eta\in(0,1]$ in a sector
    $(\thn-\delta,\thn+\delta)$. Then, under 
    Assumptions~\ref{asm:sigmaexp} and \ref{asm:sigmareg},
    there exist positive constants $C_1$, $C_2$, $\delta_1$,
    $\delta_2$, $n_0$, $t_0$, such that for $n\geq n_0$,
    $t\geq t_0$, $t\Delta(n)(\log n)^{\eta/2}\leq n\delta_1$,
    $|l|\leq n\delta_2$, we have
    \[
        \max_{k}\Prob\Lp\,
            \wandering{\origin}{n\uthn+l\uthnt}{k}{\thn}
            \geq t\Delta(n)(\log n)^{\eta/2}
        \,\Rp
        \leq C_1 \exp\Lp- C_2 t^2(\log n)^\eta\Rp.
    \]
\end{theorem}

\resetconstant

\begin{proof}
    Fix $\delta_1>0$, $\delta_2>0$, to be assumed
    appropriately small whenever required. Fix $n_0>0$,
    $t_0>0$, to be assumed appropriately large whenever
    required. Consider $n$, $t$, $l$ satisfying $n\geq n_0$,
    $t\geq t_0$, $t\Delta(n)(\log n)^{\eta/2}\leq n\delta_1$,
    $|l|\leq n\delta_2$. Let $\ba:=n\uthn+l\uthnt$. Assuming
    $\delta_2<1$ we get $\ltwo{\ba}\leq 2n$, so that, by
    Lemma~\ref{lem:boxwand}, the geodesic
    $\Gamma(\origin,\ba)$ stays inside a square of side
    length $\C n$ around $\origin$ with probability at least
    $1-e^{-\C n}$. Using $t\Delta(n)(\log n)^{\eta/2}\leq
    n\delta_1$, assuming $\delta_1$ is small enough, and
    using \eqref{A2}, we get $t^2(\log n)^\eta\leq \C
    n/\sigma(n)\leq\C n^{1-\alpha}$. Hence, the probability
    bound in the statement is trivial for $|k|\geq\Cl{244}
    n$. So let us consider $k$ satisfying $|k|\leq\Cr{244}
    n$. 
    
    We split the probability under consideration as
    \begin{multline}\label{eq:midptwand0}
        \Prob\Lp\,
            \wandering{\origin}{\ba}{k}{\thn}
            \geq 
            t\Delta(n)(\log n)^{\eta/2}
        \,\Rp
        \leq 
        \Prob\Lp\,
            \wandering{\origin}{\ba}{k}{\thn}\geq n
        \,\Rp
        \\+ 
        \Prob\Lp\,
            \wandering{\origin}{\ba}{k}{\thn}
            \in[t \Delta(n)(\log n)^{\eta/2},n]
        \,\Rp.
    \end{multline} 
    For any point $\bu$ on $\Gamma(\origin,\ba)$, we have  
    \begin{align*}
        0 
        = & T(\origin,\bu)+T(\bu,\ba)-T(\origin,\ba)\\
        = & \Lp T(\origin,\bu)-h(\bu)\Rp 
            + \Lp T(\bu,\ba)-h(\ba-\bu)\Rp 
            - \Lp T(\origin,\ba)-h(\ba)\Rp\\
          & + \Lp h(\bu)-g(\bu)\Rp 
            + \Lp h(\ba-\bu)-g(\ba-\bu)\Rp 
            - \Lp h(\ba)-g(\ba)\Rp\\
          & + \Lp g(\bu)+g(\ba-\bu)-g(\ba)\Rp,    
    \end{align*}
    so that
    \begin{align*}
             & |T(\origin,\bu)-h(\bu)|
              +|T(\bu,\ba)-h(\ba-\bu)|
              +|T(\origin,\ba)-h(\ba)|\\
        \geq &\Lp h(\bu)-g(\bu)\Rp 
            + \Lp h(\ba-\bu)-g(\ba-\bu)\Rp 
            - \Lp h(\ba)-g(\ba)\Rp
            + \Lp g(\bu) + g(\ba-\bu)-g(\ba)\Rp\\
        \geq & \Lp g(\bu)+g(\ba-\bu)-g(\ba)\Rp
            - (h(\ba)-g(\ba)).
        \numberthis\label{eq:midptwand1}
    \end{align*}
    Define a set of points $V$ in the $\RR^2$ as follows. If
    $\thnt$ is one of the axial directions, then let $V$ be
    the set of lattice points $\bu$ with 
    $\xprojthnot(\bu)=k$. If $\thnt$ is not an axial 
    direction, then let $V$ be the set of intersection points
    of the integer lattice-grid with the line
    $\xprojthnot(\bu)=k$. For $\bu\in V$ let 
    \[
    d(\bu):=\yprojthnot(\bu)-\xprojthnot(\bu)\frac{l}{n}.
    \]
    Let $V_1$ be the set of points $\bu\in V$ with 
    $|d(\bu)|\geq n$. Let $V_2$ be the set of points 
    $\bu\in V$ with 
    $t\Delta(n)(\log n)^{\eta/2}\leq |d(\bu)|\leq n$. Thus,
    if $\wandering{\origin}{\ba}{k}{\thn}\geq n$, then
    $\Gamma(\origin,\ba)$ goes through a point $\bu\in V_1$.
    Assuming $\delta_2$ is small enough and using
    Lemma~\ref{lem:auxgeom1} we get 
    \begin{equation}\label{eq:midptwand2}
        g(\bu)+g(\ba-\bu)-g(\ba)\geq\Cl{245}|d(\bu)|.
    \end{equation}
    Using $h$ satisfies GAP with exponent $\alpha$ and 
    correction factor $\phi_\eta$ in a neighborhood of $\thn$
    and assuming $\delta_2$ is small enough we get  
    \begin{equation}\label{eq:midptwand3}
        h(\ba)-g(\ba)\leq\Cl{246}\sigma(n)(\log n)^{\eta}.
    \end{equation}
    Using $|d(\bu)|\geq n$, $|l|/n\leq\delta_2<1$, and
    $|k|\leq\Cr{244}n$, we get 
    \[
        \max\LP\,
            \ltwo{\bu-\ba},\ltwo{\bu},\ltwo{\ba}
        \,\RP
        \leq\C|d(\bu)|.
    \]
    Therefore, using \eqref{A1} and \eqref{A2}, we get for all $t^\prime>0$
    \begin{multline*}
        \Prob\Lp\,\max\LP\,
        |T(\origin,\bu) - h(\bu)|,
        |T(\bu,\ba) - h(\ba-\bu)|,
        |T(\origin,\ba) - h(\ba)|
        \,\RP\geq t^\prime\sigma(|d(\bu)|)\,\Rp\\
        \leq\C e^{-\C t^\prime}.
    \end{multline*}
    Using this with \eqref{eq:midptwand1}-\eqref{eq:midptwand3}, and \eqref{A2}, we get 
    \begin{align*}
        &\Prob\Lp\,\wandering{\origin}{\ba}{k}{\thn}\geq n\,\Rp \\
        \leq & 
        \sum_{\bu\in V_1}\C\exp\Lp-\C\Lp\Cr{245}|d(\bu)| - \Cr{246} \sigma(n)(\log n)^{\eta}\Rp/\sigma(|d(\bu)|)\Rp\\
        \leq & \sum_{u\in V_1}\C\exp\Lp-\C |d(\bu)|/\sigma(|d(\bu)|)\Rp\\
        \leq & \C\exp\Lp -\C n / \sigma(n)\Rp.
    \numberthis\label{eq:midptwand4}
    \end{align*}
    If 
    $\wandering{\origin}{\ba}{k}{\thn}
    \in[t\Delta(n)(\log n)^{\eta/2},n]$, then $\Gamma(\origin,\ba)$ goes through a point $\bu\in V_2$.
    Assuming $\delta_2$ is small enough and using Lemma~\ref{lem:auxgeom1} we get
    \begin{equation}\label{eq:midptwand5}
    g(\bu) + g(\ba-\bu)-g(\ba)
    \geq \Cl{2416} \frac{|d(\bu)|^2}{n}     
    \geq \Cr{2416} t^2 \sigma(n)(\log n)^{\eta}.
    \end{equation}
    Since $|d(\bu)|\leq n$, we have 
    \[
        \max\LP\,\ltwo{\bu-\ba},\ltwo{\bu},\ltwo{\ba}\,\RP
        \leq\C n.
    \]
    Hence, using Assumptions~\ref{asm:sigmaexp} and \ref{asm:sigmareg} we get, for all $t^\prime>0$
    \[
        \Prob\Lp\,\max\LP\, 
        |T(\origin,\bu) - h(\bu)|, 
        |T(\bu,\ba) - h(\ba-\bu)|, 
        |T(\origin,\ba) - h(\ba)|\,\RP
        \geq t^\prime\sigma(n)\,\Rp\leq\C e^{-\C t^\prime}.
    \]
    Using this with \eqref{eq:midptwand1}, \eqref{eq:midptwand3}, \eqref{eq:midptwand5} we get 
    \begin{align*}
     & \Prob\Lp\,
         \wandering{\origin}{\ba}{k}{\thn}
         \in[t\Delta(n)(\log n)^{1/2},n]
        \,\Rp\\
    \leq & 
    \sum_{\bu\in V_2} 
    \C\exp\Lp-\C\Lp\Cr{2416}\frac{|d(\bu)|^2}{n}
        -\Cr{246}\sigma(n)(\log n)^\eta\Rp/\sigma(n)\Rp\\ 
    \leq & \C\exp\Lp-\C t^2(\log n)^{\eta}\Rp.
    \numberthis\label{eq:midptwand6}
    \end{align*}
    Assuming $\delta_1$ is small enough and combining 
    \eqref{eq:midptwand4}, \eqref{eq:midptwand6}, 
    \eqref{eq:midptwand0} we get
    \[
    \Prob\Lp\,
        \wandering{\origin}{\ba}{k}{\thn}
        \geq t\Delta(n)(\log n)^{\eta/2}\,\Rp
    \leq \C e^{-\C t^2(\log n)^{\eta}}.
    \]
    This completes the proof of Theorem~\ref{thm:midptwand}. 
\end{proof}

\begin{corollary}\label{cor:overallwand}
    Let $\thn$ be a direction of type I. Then, under Assumptions~\ref{asm:sigmaexp} and \ref{asm:sigmareg}, there exist positive constants $C_1$, $C_2$, $\delta_1$, $\delta_2$, $n_0$, $t_0$, such that for $n\geq n_0$, $t\geq t_0$, $t\Delta(n)(\log n)^{1/2}\leq n\delta_1$, $|l|\leq n\delta_2$, we have 
    \[
        \Prob\Lp\,
            \max_k
            \wandering{\origin}{n\uthn+l\uthnt}{k}{\thn}
            \geq
            t\Delta(n)(\log n)^{1/2}
        \,\Rp
        \leq C_1 e^{-C_2 t^2\log n}.
    \]
\end{corollary}

\resetconstant

\begin{proof} 
    Fix $\delta_1>0$, $\delta_2>0$, to be assumed
    appropriately small whenever required. Fix $n_0>0$,
    $t_0>0$, to be assumed appropriately large whenever
    required. Consider $n$, $t$, and $l$, satisfying 
    $n\geq n_0$, $t\geq t_0$, 
    $t\Delta(n)(\log n)^{1/2}\leq n\delta_1$, 
    $|l|\leq n\delta_2$.
    Let $\ba:=n\uthn+l\uthnt$. Assuming $\delta_2<1$ we get 
    $\ltwo{\ba}\leq 2n$, so that by Lemma~\ref{lem:boxwand}, 
    the geodesic $\Gamma(\origin,\ba)$ stays inside a square 
    of side length $\C n$ around $\origin$ with probability 
    at least $1-e^{-\Cl{252} n}$. On this event 
    $|\xprojthnot(\bu)|\leq\Cl{251}n$ for all $\bu$ in 
    $\Gamma(\origin,\ba)$.
    Assuming $\delta_1<1$ and using \eqref{A2} we get 
    $t^2\log n\leq n/\sigma(n)\leq\C n^{1-\alpha}$.
    Thus
    \begin{equation}\label{eq:251}
        \Prob\Lp\,
            \max_{|k|\geq\Cr{251}n}
            \wandering{\origin}{\ba}{k}{\thn}
            \geq t\Delta(n)(\log n)^{1/2}
        \,\Rp\leq e^{-\Cr{252}n}
        \leq \C e^{-\C t^2\log n}.
    \end{equation}
    Using Proposition~\ref{prop:nrflog},  
    Theorem~\ref{thm:midptwand}, and a union bound, we get
    \begin{equation}\label{eq:252}
        \Prob\Lp\,
            \max_{|k|\leq\Cr{251}n} 
            \wandering{\origin}{\ba}{k}{\thn}
            \geq t\Delta(n)(\log n)^{1/2}
            \,\Rp
        \leq \C e^{-\C t^2\log n}.
    \end{equation}
    Combining \eqref{eq:251} and \eqref{eq:252} completes the
    proof. 
\end{proof} 

Now we show that wandering of a geodesic at a distance $k$ 
from an endpoint is at most of the order of $\Delta(k)$ with 
some logarithmic correction factor. 

\begin{theorem}\label{thm:endwandlog}
    Let $\thn$ be a direction of type I. Suppose $h$
    satisfies GAP with exponent $\alpha$ and correction
    factor $\phi_\eta$ for some $\eta\in(0,1]$ in a
    neighborhood of $\thn$. Then, under the 
    Assumptions~\ref{asm:sigmaexp} and \ref{asm:sigmareg}, 
    there exist positive constants $C_1$, $C_2$, $\delta_1$, 
    $\delta_2$, $k_0$, $n_0$, $t_0$, such that for 
    $n\geq n_0$, $t\geq t_0$, $k\geq k_0$, $k<n$,
    $t\Delta(k)(\log k)^{\eta/2}\leq k\delta_1$, 
    $|l|\leq n\delta_2$, we have
    \[
    \Prob\Lp\,
        \wandering{\origin}{n\uthn+l\uthnt}{k}{\thn}
        \geq t\Delta(k)(\log k)^{\eta/2}
        \,\Rp
    \leq C_1 e^{-C_2 t^2(\log k)^{\eta}}.
    \]
\end{theorem} 

\resetconstant

\begin{proof}
    Fix $\delta_1>0$, $\delta_2>0$, to be assumed
    appropriately small whenever required. Fix $k_0>0$,
    $n_0>0$, $t_0>0$, to be assumed appropriately large
    whenever required. Consider $k$, $n$, $l$, $t$ satisfying
    $n\geq n_0$, $t\geq t_0$, $k\geq k_0$, 
    $t\Delta(k)(\log k)^{\eta/2}\leq k\delta_1$, 
    $|l|\leq n\delta_2$. Let $\ba:=n\uthn+l\uthnt$. 
    
    The strategy of the proof is as follows. We define two
    sequences of points $(\ba_p)_{p=0}^{m+1}$ and
    $(\bb_p)_{p=0}^{m+1}$ in such a manner that, if for all
    $p$, the geodesics $\Gamma(\origin,\ba_p)$ and
    $\Gamma(\origin,\bb_p)$ do not wander ``excessively" then
    $\wandering{\origin}{\ba}{k}{\thn}
    \leq t\Delta(k)(\log k)^{\eta/2}$ as we require. We
    define events $\LP\event_p^i:1\leq i\leq 4,0\leq p\leq
    m\RP$ in which these geodesics do not wander excessively.
    Then we show that these events have appropriately high
    probability.
    
    Fix $1<\zeta<2/(1+\beta)$ and let
    $\epsilon:=1-\zeta(1+\beta)/2$. Fix $\lambda>1$,
    $\lambda^\prime>0$. Later we choose $\lambda$ so that
    \eqref{eq:endwandlog3} holds, and we choose
    $\lambda^\prime$ based on $\lambda$ so that
    \eqref{eq:lambdaprime} holds. Let $m\geq 0$ be such that
    $\lambda^m k < n \leq \lambda^{m+1} k$. Let
    $\ba_{m+1}:=\ba$, $\bb_{m+1}:=\ba$, and for 
    $0\leq p\leq m$, let $\xprojthnot(\ba_p)=\lambda^{p} k$,
    $\xprojthnot(\bb_p)=\lambda^{p} k$,
    \begin{gather*}
        \yprojthnot(\ba_{p})=
        \frac{\xprojthnot(\ba_{p})}{\xprojthnot(\ba_{p+1})}
        \yprojthnot(\ba_{p+1})
        +\lambda^\prime 
        t\Delta(\lambda^{\zeta p} k)
        (\log\lambda^{\zeta p} k)^{\eta/2},\\
        \yprojthnot(\bb_{p})=
        \frac{\xprojthnot(\bb_{p})}{\xprojthnot(\bb_{p+1})}
        \yprojthnot(\bb_{p+1})
        -\lambda^\prime 
        t\Delta(\lambda^{\zeta p} k)
        (\log\lambda^{\zeta p} k)^{\eta/2}.
    \end{gather*}
    For any $0\leq p\leq q\leq m$, using \eqref{A2} we get
    \begin{multline*}
        \frac{\lambda^{-q}\gwq}{\lambda^{-p}\gwp} 
        \leq\C\lambda^{-\epsilon(q-p)}\Lp\frac{q\zeta\log\lambda+\log k}{p\zeta\log\lambda+\log k}\Rp^{\eta/2}\\
        \leq\C\lambda^{-\epsilon(q-p)}\Lp 1+\frac{(q-p)\zeta\log \lambda}{p\zeta\log\lambda+\log k}\Rp^{\eta/2}
        \leq\C\lambda^{-\epsilon(q-p)}\Lp 1+(q-p)\log\lambda\Rp^{\eta/2}.
    \end{multline*}
    Therefore, we can choose $\lambda^\prime$, depending on
    $\lambda$, such that 
    \begin{equation}\label{eq:lambdaprime}
        \sum_{q=p}^{m}\lambda^{-(q-p)}\gwq
        \leq\frac{1}{\lambda^\prime}\gwp.
    \end{equation}
    Therefore, for $0\leq p \leq m$
    \begin{gather*}
        \yprojthnot(\ba_{p}) = \lambda^{p} k \frac{l}{n} + 
        \lambda^\prime t \sum_{q=p}^{m} 
        \lambda^{-(q-p)}\gwq\\ \leq \lambda^{p} k\frac{l}{n} 
        + t \gwp,\\
        \yprojthnot(\bb_{p}) = \lambda^{p} k \frac{l}{n} -
        \lambda^\prime t \sum_{q=p}^{m}
        \lambda^{-(q-p)}\gwq\\
        \geq \lambda^{p} k\frac{l}{n} - t \gwp.
    \end{gather*}
    For $0\leq p\leq m$, let 
    \[
    \overline{\segment}_p
    :=\LP\,\bx\in\RR^2\,:\,\xprojthnot(\bx)
        =\lambda^p k\,\RP. 
    \]
    So $\ba_p$, $\bb_p$ are on this line. Let $\segment_p$ be
    the segment joining $\ba_p$ and $\bb_p$. Therefore,
    $\wandering{\origin}{\ba}{k}{\thn}
    \leq t\Delta(k)(\log k)^{\eta/2}$ if and only if, 
    whenever $\Gamma(\origin,\ba)$ intersects 
    $\overline{\segment}_0$, it intersects $\segment_0$. We 
    define the events 
    $\LP\event^i_p:1\leq i\leq 4, 0\leq p\leq m\RP$ 
    in a way that, if 
    $\mathbb{T}\in\cap_{i=1}^4\cap_{p=0}^m\event^i_p$ then 
    for all $0\leq p\leq m$, whenever $\Gamma(\origin,\ba)$ 
    intersects $\overline{\segment}_p$, it intersects 
    $\segment_p$, see Figure~\ref{fig:inductive}.  
    \begin{figure}[H]
        \centering
        \includegraphics[width=0.45\linewidth]{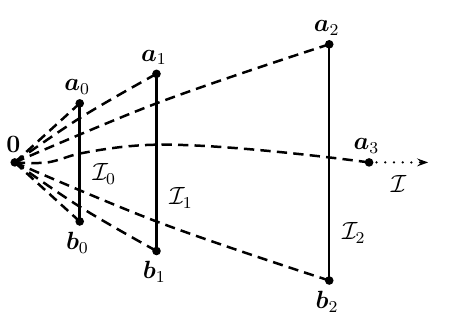}
        \caption{Rough sketch for $m=2$: 
        $\ba_{m+1}=\bb_{m+1}=\ba$; for $0\leq p\leq m$ the
        segment joining $\ba_p$ and $\bb_p$ is $\segment_p$,
        and extending $\segment_p$ we get the line 
        $\overline{\segment}_p$; for all $0\leq p\leq m$ if
        the geodesics $\Gamma(\origin,\ba_p)$ and 
        $\Gamma(\origin,\bb_p)$ do not wander excessively 
        then for all $0\leq p\leq m$ the geodesic 
        $\Gamma(\origin,\ba)$ intersects $\segment_p$ 
        whenever it intersects $\overline{\segment}_p$.}
        \label{fig:inductive}
    \end{figure}
    For $0\leq p \leq m$, define the events  
    \begin{gather*}
        \event^1_p:=\LP\,
        \wandering{\origin}{\ba_{p+1}}{\lambda^p k}{\thn}
        \leq\lambda^\prime t\gwp
        \,\RP,\\
        \event^2_p:=\LP\,
        \wandering{\origin}{\bb_{p+1}}{\lambda^p k}{\thn}
        \leq\lambda^\prime t\gwp
        \,\RP.
    \end{gather*}
    By construction of the sequences $(\ba_p)_{p=0}^{m+1}$ 
    and $(\bb_p)_{p=0}^{m+1}$, we get that for 
    $0\leq p\leq m$, if 
    $\mathbb{T}\in\event_p^1\cap\event_p^2$, then the 
    geodesics $\Gamma(\origin,\ba_{p+1})$ and 
    $\Gamma(\origin,\bb_{p+1})$ intersect $\segment_p$ 
    whenever they intersect $\overline{\segment}_p$. We want 
    the geodesic $\Gamma(\origin,\ba)$ to intersect 
    $\segment_p$ whenever it intersects 
    $\overline{\segment_p}$, for all $0\leq p\leq m$;
    Figure~\ref{fig:inductive} shows the ideal situation. But
    $\mathbb{T}\in\cap_{p=0}^m (\event_p^1\cap\event_p^2)$ is
    not sufficient for this because the geodesics can
    backtrack. So we define the events $\event_p^3$ and
    $\event_p^4$ for $0\leq p\leq m$ where the geodesics do
    not backtrack. Define the half-line 
    \[
        \segment:=\LP\,\bx\in\RR^2\,:\,
        \xprojthnot(\bx)\geq n,\, \frac{\yprojthnot(\bx)}{\xprojthnot(\bx)}
        =\frac{l}{n}\,\RP. 
    \]
    If in addition to 
    $\mathbb{T}\in\cap_{p=0}^m(\event_p^1\cap\event_p^2)$, we
    also have that for all $0\leq p\leq m$, the geodesics
    $\Gamma(\origin,\ba_p)$ and $\Gamma(\origin,\bb_p)$ do
    not intersect $\segment$, then we get that for all $0\leq
    p\leq m$ the geodesic $\Gamma(\origin,\ba)$ intersects
    $\segment_p$ whenever it intersects
    $\overline{\segment}_p$. Therefore we define the events
    $\event^3_p$, $\event^4_p$ for $0\leq p\leq m$ as 
    follows. Let
    \begin{gather*}
        \event^3_m:=
            \LP\,\max_{n^\prime\geq n}
            \wandering{\origin}{\ba_{m}}{n^\prime}{\thn}
            \leq\lambda^\prime t\gwm\,\RP,\\
        \event^4_m:=
            \LP\,\max_{n^\prime\geq n}
            \wandering{\origin}{\bb_{m}}{n^\prime}{\thn}
            \leq\lambda^\prime t\gwm\,\RP,
    \end{gather*}
    and for $0\leq p <m$ let
    \[
        \event^3_p:=
        \LP\,
        \text{Diam}\Lp\Gamma(\origin,\ba_{p})\Rp\leq n
        \,\RP,
        \quad
        \event^4_p:=
        \LP\,
        \text{Diam}\Lp\Gamma(\origin,\bb_{p})\Rp\leq n
        \,\RP.
    \]
    So, for $0\leq p\leq m$, if $\mathbb{T}\in\event^3_p$,
    then $\Gamma(\origin,\ba_p)$ does not intersect 
    $\segment$, and if $\mathbb{T}\in\event^4_p$, then 
    $\Gamma(\origin,\bb_p)$ does not intersect $\segment$. 
    Thus, if $\mathbb{T}\in\cap_{i=1}^4\cap_{p=0}^m\event_p^i
    $, then for all $0\leq p\leq m$, whenever 
    $\Gamma(\origin,\ba)$ intersects $\overline{\segment}_p$ 
    it intersects $\segment_p$. In this event 
    $\wandering{\origin}{\ba}{k}{\thn}\leq t\Delta(k)(\log 
    k)^{\eta/2}$ as required. So now we show these events 
    have appropriately high probability.
    
    We use Theorem~\ref{thm:midptwand} to bound 
    $\Prob((\event^i_p)^c)$ for $0\leq p\leq m$ and $i=1,2$. 
    For a fixed $p$ and $i=1$ we use the following 
    parameters. 
    \begin{align*}
        & \tilde{\eta}:=\eta,
        \quad
        \tilde{n}:=\xprojthnot(\ba_{p+1}),
        \quad
        \tilde{l}:=\yprojthnot(\ba_{p+1}),\\ 
        & \tilde{t}:=
        \lambda^\prime t 
            \frac{\gwp}
                {\Delta(\tilde{n})(\log\tilde{n})^{\eta/2}}.
    \end{align*}
    We need to verify $\tilde{n}\geq\tilde{n}_0$, 
    $\tilde{t}\geq\tilde{t}_0$, 
    $\tilde{t}\Delta(\tilde{n})
    (\log\tilde{n})^{\tilde{\eta}/2}
    \leq\tilde{n}\delta_1$, 
    $|\tilde{l}|\leq\tilde{n}\tilde{\delta}_2$. The condition
    $\tilde{n}\geq\tilde{n}_0$ holds by taking $k_0$ large
    enough because $\tilde{n}\geq k\geq k_0$. Using
    $\tilde{n}\leq\lambda^{p+1}k$ and \eqref{A2} we get
    $\tilde{t}\geq\C t$. So $\tilde{t}\geq\tilde{t}_0$ holds 
    by choosing $t_0$ large enough. Using 
    $\tilde{n}\geq\lambda^p k$, \eqref{A2}, 
    $t\Delta(k)(\log k)^{\eta/2}\leq k\delta_1$, 
    and assuming $\delta_1$ is small enough, we get
    \[
        \frac{1}{\tilde{n}}
        \tilde{t}\Delta(\tilde{n})(\log\tilde{n})^{\eta/2}
        \leq 
        \C\frac{1}{\lambda^p k} 
        t \Delta(\lambda^{p} k)(\log \lambda^{p} k)^{\eta/2}
        \leq
        \C \frac{1}{k} t \Delta(k)(\log k)^{\eta/2}
        \leq
        \C \delta_1\leq\tilde{\delta_1}.
    \]
    If $p=m$, then 
    $|\tilde{l}|/\tilde{n}
    =|l|/n\leq\delta_2\leq\tilde{\delta}_2$. 
    For $p<m$, using 
    $t\Delta(k)(\log k)^{\eta/2}\leq k\delta_1$, 
    $|l|\leq n\delta_2$, \eqref{A2}, and assuming $\delta_1$,
    $\delta_2$ are small enough, we get
    \begin{align*}
        \frac{|\tilde{l}|}{\tilde{n}}
        & = \frac{|\yprojthnot(\ba_{p+1})|}
                {\xprojthnot(\ba_{p+1})}\\ 
        & \leq \frac{1}{\lambda^{p+1} k}
            \Lp\lambda^{p+1} k\frac{|l|}{n}+t\gwpp\Rp\\
        & \leq \delta_2 +                    
            \delta_1 \frac{\gwpp}
                {\lambda^p\Delta(k)(\log k)^{\eta/2}}\\
        & \leq \delta_2 +                          
            \C \delta_1\lambda^{-\epsilon(p+1)}
            \Lp 1 + \frac{(p+1)\zeta\log\lambda}
                {\log k}\Rp^{\eta/2}\\
        & \leq \delta_2 + \C\delta_1\\
        & \leq \tilde{\delta_2}.
    \end{align*}
    So the conditions for applying 
    Theorem~\ref{thm:midptwand} hold. Using 
    $\tilde{n}\leq\lambda^{p+1}k$ and \eqref{A2} we get 
    \[
        \tilde{t}^2(\log\tilde{n})^\eta
    \geq\C t^2\lambda^{p\epsilon}(\log k)^{\eta}. 
    \]
    Therefore, applying Theorem~\ref{thm:midptwand} we get
    \begin{equation}\label{eq:endwandlog1}
        \Prob\Lp\,
        (\event^1_p)^c
        \,\Rp
        \leq\C\exp\Lp-\C t^2\lambda^{p\epsilon}
            (\log k)^\eta\Rp.
    \end{equation}
    Similar bounds hold for $\event^2_p$ for $0\leq p\leq m$,
    $\event^3_m$, $\event^4_m$. For $0\leq p<m$, 
    $\ltwo{\ba_p} \leq 2 \lambda^p k $. We have 
    $n>\lambda^m k$. Thus, assuming $\lambda$ large enough, 
    we get for $i=3,4$ and $0\leq p < m$, using 
    Lemma~\ref{lem:boxwand}, 
    \begin{equation}\label{eq:endwandlog3}
        \Prob\Lp\,
            (\event^i_p)^c
            \,\Rp
        \leq\C\exp\Lp-\C\lambda^p k\Rp.
    \end{equation}
    Combining \eqref{eq:endwandlog1} and
    \eqref{eq:endwandlog3} we get 
    \[
        \Prob\Lp\,
            \wandering{\origin}{\ba}{k}{\thn}
            \geq t\Delta(k)(\log k)^{\eta/2}
        \,\Rp
        \leq \Prob \Lp\,
            \cup_{i=1}^4\cup_{p=0}^m(\event^i_p)^c
        \,\Rp
        \leq \C e^{-\C t^2 (\log k)^\eta}.
    \]
    This completes the proof of Theorem~\ref{thm:endwandlog}.
\end{proof}

As a corollary we can deal with wandering of a geodesic
within a fixed distance from an endpoint. 

\begin{corollary}\label{cor:endwandlogspl}
    Let $\thn$ be a direction of type I. Then, under the
    Assumptions~\ref{asm:sigmaexp} and \ref{asm:sigmareg},
    there exist positive constants $C_1$, $C_2$, $\delta_1$,
    $\delta_2$, $k_0$, $n_0$, $t_0$, such that for 
    $k\geq k_0$, $n\geq n_0$, $t\geq t_0$, 
    $t\Delta(k)(\log k)^{1/2}\leq k\delta_1$, 
    $|l|\leq n\delta_2$, we have
    \[
    \Prob\Lp\,
    \max_{k^\prime\leq k}
    \wandering{\origin}{n\uthn+l\uthnt}{k^\prime}{\thn}
    \geq t\Delta(k)(\log k)^{1/2}
    \,\Rp
    \leq C_1 e^{-C_2 t^2\log k}.
    \]
\end{corollary} 

\resetconstant

\begin{proof} 
    Fix $\delta_1>0$, $\delta_2>0$, to be assumed 
    appropriately small whenever required. Fix $k_0>0$,
    $n_0>0$, $t_0>0$, to be assumed large enough whenever 
    required. Let $\ba:=n\uthn+l\uthnt$. Define 
    \[
        \segment:=\LP\,\bv\in\RR^2\,:\,
        \xprojthnot(\bv)=k,\,
        |\yprojthnot(\bv)-k\frac{l}{n}|
        \leq\frac{t}{2}\Delta(k)(\log k)^{1/2}\RP. 
    \]
    The event
    $\max_{k^\prime\leq k}
    \wandering{\origin}{\ba}{k^\prime}{\thn}
    \geq t\Delta(k)(\log k)^{1/2}$, 
    can happen in two ways. Either we have 
    $\wandering{\origin}{\ba}{k}{\thn}
    \geq(t/2)\Delta(k)(\log k)^{1/2}$, or
    $\Gamma(\origin,\ba)$ passes through some point
    $\bv$ in the segment $\segment$ and $\max_{k^\prime}
    \wandering{\origin}{\bv}{k^\prime}{\thnt}
    \geq (t/2)\Delta(k)(\log k)^{1/2}$.
    In the first case, by Theorem~\ref{thm:endwandlog} and 
    Proposition~\ref{prop:nrflog}, we get 
    \begin{equation}\label{eq:endlogspl1}
        \Prob\Lp\;
            \wandering{\origin}{\ba}{k}{\thn}
            \geq\frac{t}{2}\Delta(k)(\log k)^{1/2}
        \;\Rp
        \leq C_1 e^{-C_2 t^2\log k}.
    \end{equation}
    For the second case, consider $\bv\in\segment$. We apply
    Corollary~\ref{cor:overallwand} with 
    \[
        \tilde{\thn}:=\thn,\quad
        \tilde{n}:=k,\quad
        \tilde{l}:=\yprojthnot(\bv),\quad
        \tilde{t}:=t.
    \]
    Then 
    \[
        \frac{|\tilde{l}|}{\tilde{n}}
        \leq
        \frac{|l|}{n}
        +\frac{1}{k}\frac{t}{2}\Delta(k)(\log k)^{1/2}
        \leq 
        \delta_1+\frac{\delta_2}{2}.
    \]
    Also 
    \[
        \frac{1}{\tilde{n}}
        \tilde{t}\Delta(\tilde{n})(\log\tilde{n})^{1/2}
        \leq
        \frac{1}{k}t\Delta(k)(\log k)^{1/2}
        \leq
        \delta_2.
    \]
    Therefore, assuming $\delta_1$, $\delta_2$ are small 
    enough, we get
    \[
        \Prob\Lp\;
            \max_{k^\prime\in[0,k]}
            \wandering{\origin}{\bv}{k^\prime}{\thn}
            \geq\frac{t}{2}\Delta(k)(\log k)^{1/2}
            \;\Rp
        \leq C_3 e^{-C_4 t^2\log k}.
    \]
    Taking a union bound over $2t\Delta(k)(\log k)^{1/2}$ 
    many choices of $\bv$ we get 
    \begin{equation}\label{eq:endlogspl2}
        \Prob\Lp\;
        \max_{k^\prime}
        \wandering{\origin}{n\uthn+l\uthnt}{k^\prime}{\thn}
        \geq t\Delta(k)(\log k)^{1/2}
        \;\Rp
        \leq C_5 e^{-C_6 t^2\log k}.
    \end{equation}
    Combining \eqref{eq:endlogspl1}, \eqref{eq:endlogspl2} 
    completes the proof of Corollary~\ref{cor:endwandlogspl}.
\end{proof}

\begin{remark}\label{remark:nonlattice2}
    Recall from Notation~\ref{notn:nonlattice} that geodesic
    between points $\bu$, $\bv$ which are not necessarily 
    lattice points is defined as the geodesic between 
    $\lfloor\bu\rfloor$ and $\lfloor\bv\rfloor$. In 
    Theorem~\ref{thm:midptwand}, 
    Corollary~\ref{cor:overallwand}, 
    Theorem~\ref{thm:endwandlog}, and 
    Corollary~\ref{cor:endwandlogspl}, we are dealing with 
    geodesics having one endpoint $\origin$ and we are 
    measuring the wandering at a distance from $\origin$. 
    Later while applying these results we may have a point of
    $\RR^2$ in place of $\origin$. This does not cause any
    major complication, i.e., bounds that hold for wandering
    of $\Gamma(\origin,\bu)$ also hold for wandering of
    $\Gamma(\bv,\bu+\bv)$, where $\bu$ and $\bv$ are not
    necessarily lattice points. 
\end{remark} 

\section{Preliminary upper bound of the transverse 
increments}

In this section, our principal objective is to prove
Theorem~\ref{thm:logupinc}. This is a special case of
Theorem~\ref{thm:loglogupinc} which is our main upper bound
on the transverse increments. In Section~\ref{sec:ghloglog},
we use Theorem~\ref{thm:logupinc} to prove
Theorem~\ref{thm:nrf2} which is a refinement of the bound on
nonrandom fluctuations of Proposition~\ref{prop:nrflog}. We
use this refinement to prove Theorem~\ref{thm:loglogupinc} in
Section~\ref{sec:loglogupinc}.

We need the following result on curvature of the boundary of 
the limit shape. We skip the proof because the result is 
essentially same as Lemma~2.7 of \cite{KenGeoBiGeo}. 

\begin{lemma}\label{lem:auxgeom2}
    Let $\thn$ be a direction of type II. 
    Then there exist positive constants $C$, $\delta_1$,
    $\delta_2$, such that for $d$, $k>0$, $L$, satisfying
    $|d|\leq k\delta_1$, $|L|\leq k\delta_2$, we have
    \[
        |g(k\uthn+(d+L)\uthnt)-g(k\uthn+d\uthnt)|
        \leq C\Lp\frac{L^2}{k}+\frac{|d||L|}{k}\Rp.
    \]
\end{lemma}

The preliminary upper bound of the transverse increments is 
the following.

\begin{theorem}\label{thm:logupinc}
    Let $\thn$ be a direction of both type I and II. For
    $n>0$, $L>0$, let $\segment(n,L)$ and $\diff(n,L)$ be
    as defined in Theorem~\ref{thm:loglogupinc}. Then, under
    Assumptions~\ref{asm:sigmaexp} and \ref{asm:sigmareg},
    there exist positive constants $C_1$, $C_2$, $L_0$,
    $n_0$, $t_0$, such that for 
    $L\geq L_0$, $n\geq n_0$, $t\geq t_0$, $L\leq\Delta(n)$, 
    we have  
    \[
    \Prob\Lp\,
        \diff(n,L)\geq t\log L\sigma(\Delta^{-1}(L))
    \,\Rp
    \leq C_1 \exp\Lp - C_2 t\log L\Rp.
    \]
\end{theorem}  

\resetconstant

\begin{proof}
    Let us assume without loss of generality that
    $\thn\in[0,\pi/4]$. Fix $L_0>0$, $n_0>0$, $t_0>0$. 
    We assume $L_0$, $n_0$, $t_0$ are large whenever
    required. Consider $L\geq L_0$, $n\geq n_0$, $t\geq t_0$,
    $L\leq\Delta(n)$. Based on the values of $t$ we consider
    two cases. Because for suitably large values of $t$ we
    are in a large deviation regime, and the proof is 
    straightforward. 
    
    \vspace{\baselineskip}
    
    \textit{Case I:} Suppose 
    $t\geq 4\mu L(\sigma(\Delta^{-1}(L))\log L)^{-1}$, 
    where $\mu$ is the expected passage edge-weight.
    Since $\segment(n,L)$ has width $L$, the lattice points 
    corresponding to the points in $\segment(n,L)$ can be
    joined by a lattice path of at most 
    $\lceil 2 L \rceil$ edges. 
    Hence $\diff(n,L)\leq X_1+\cdots+X_{\lceil 2 L \rceil}$, 
    where $X_i$'s are i.i.d.\ random variables which have the
    same distribution as that of the edge-weights. Therefore,
    \begin{multline*}
        \Prob\Lp\,
            \diff(n,L)\geq t\log L\sigma(\Delta^{-1}(L))
        \,\Rp
        \leq 
        \Prob\Lp\, 
            X_1+\cdots+X_{\lceil 2 L\rceil}\geq
            t\log L\sigma(\Delta^{-1}(L))
        \,\Rp\\
        \leq \C\exp\Lp-\C t\log L\sigma(\Delta^{-1}(L))\Rp
        \leq \C\exp\Lp-\C t\log L\Rp.
    \end{multline*}
    This concludes the proof in this case.
    
    \vspace{\baselineskip}
    
    \textit{Case II:} Suppose 
    \begin{equation}\label{eq:logupinc0}
        t\leq 4\mu L(\sigma(\Delta^{-1}(L))\log L)^{-1}.
    \end{equation}
    Let 
    \[
        J:=\LT-t^{1/2}L(\log\Delta^{-1}(L))^{1/2},
        \Lp 1-\frac{\Delta^{-1}(L)}{n}\Rp L
        + t^{1/2}L\Lp\log\Delta^{-1}(L)\Rp^{1/2}\RT 
    \]
    and        
    \begin{equation}\label{eq:logupinc1}
        \segment^\ast:=
        \left\{\,
        \bx\in\RR^2
        \,:\,
        \xprojthnot(\bx)=n-\Delta^{-1}(L),\,
        \yprojthnot(\bx)\in J,
        \,\right\}.
    \end{equation}
    Define the event 
    \[
        \event: \mbox{$\Gamma(\origin,\bu)$ intersects
        $\segment^\ast$ for all $\bu\in\segment(n,L)$}.
    \]
    If $\mathbb{T}\not\in\event$, then 
    $\wandering{\bu}{\origin}{\Delta^{-1}(L)}{-\thn}
    \geq t^{1/2}L(\log\Delta^{-1}(L))^{1/2}$
    for some $\bu\in\segment(n,L)$. Consider a point 
    $\bu\in\segment(n,L)$. We want to apply 
    Theorem~\ref{thm:endwandlog} with the variables
    \begin{align*}
        &\tilde{k} :=\Delta^{-1}(L),\quad
        \tilde{l} :=\yprojthnot(\bu),\quad
        \tilde{n} :=\yprojthnot(\bu)=n,\\
        &\tilde{t} :=t^{1/2},\quad
        \tilde{\thn} :=-\thn,\quad
        \tilde{\eta} :=1.
    \end{align*}
    We need to verify that these variables satisfy the 
    conditions of Theorem~\ref{thm:endwandlog}.
    The point $\bu$ is not necessarily a lattice point.
    But this issue has been addressed in 
    Remark~\ref{remark:nonlattice2}.  
    The conditions on $\tilde{\thn}$ hold because by 
    assumption $\thn$ is of type I, and by 
    Proposition~\ref{prop:nrflog} $h$ satisfies GAP with 
    correction factor $\phi_1$ in all directions.
    Now, we need to verify the conditions 
    $\tilde{k}\geq\tilde{k_0}$,
    $\tilde{n}\geq\tilde{n_0}$,
    $\tilde{t}\geq\tilde{t_0}$, 
    $\tilde{t}\Delta(\tilde{k})(\log\tilde{k})^{1/2}
    \leq\tilde{k}\tilde{\delta_1}$, 
    $|\tilde{l}|\leq\tilde{n}\tilde{\delta_2}$.
    Assuming $n_0$, $L_0$, $t_0$ are large enough, we get  
    $\tilde{k}\geq\tilde{k_0}$, 
    $\tilde{n}\geq\tilde{n_0}$, 
    $\tilde{t}\geq\tilde{t_0}$.
    From $\bu\in\segment(n,L)$ 
    we get $|\yprojthnot(u)|\leq L$. 
    Therefore, using $L\leq\Delta(n)$, \eqref{A2}, and 
    assuming $n_0$ is large enough, we get 
    \[
        \frac{|\tilde{l}|}{\tilde{n}}
        \leq\frac{L}{n}
        \leq\frac{\Delta(n)}{n}
        \leq n_0^{-(1-\beta)/2}
        \leq\tilde{\delta_2}.
    \]
    Using \eqref{eq:logupinc0}, \eqref{A2}, and assuming 
    $L_0$ is large enough, we get
    \begin{multline*}
        \frac{1}{\tilde{k}}\tilde{t}
        \Delta(\tilde{k})(\log\tilde{k})^{1/2} =
        \frac{t^{1/2}L\Lp\log\Delta^{-1}(L)\Rp^{1/2}}
        {\Delta^{-1}(L)}
        \leq
            \C\frac{L^{3/2}(\log\Delta^{-1}(L))^{1/2}}
            {\Delta^{-1}(L)(\sigma(\Delta^{-1}(L))^{1/2}
                (\log L)^{1/2}}\\
        \leq
            \C\frac{L^{1/2}}{(\Delta^{-1}(L))^{1/2}} 
        \leq
            \C L^{-\frac{(1-\beta)}{2(1+\beta)}}
        \leq
            \C L_0^{-\frac{(1-\beta)}{2(1+\beta)}}
        \leq
            \tilde{\delta_1}.
    \end{multline*}
    Therefore, all the conditions for applying
    Theorem~\ref{thm:endwandlog} are satisfied, and we get
    \[
    \Prob\Lp\,
        \wandering{\bu}{\origin}{\Delta^{-1}(L)}{-\thn}
        \geq t^{1/2}L(\log\Delta^{-1}(L))^{1/2}
    \,\Rp\leq \C e^{-\C t\log L}.
    \]
    Therefore, taking a union bound over all 
    $\bu\in\segment(n,L)$ we get  
    \begin{equation}\label{eq:logupinc3}
        \Prob\Lp\,\event^c\,\Rp\leq \C e^{-\C t\log L}.
    \end{equation}
    So in order to complete the proof we consider 
    $\mathbb{T}\in\event$. 
    \begin{figure}[H]
        \centering
        \includegraphics[width=0.65\linewidth]{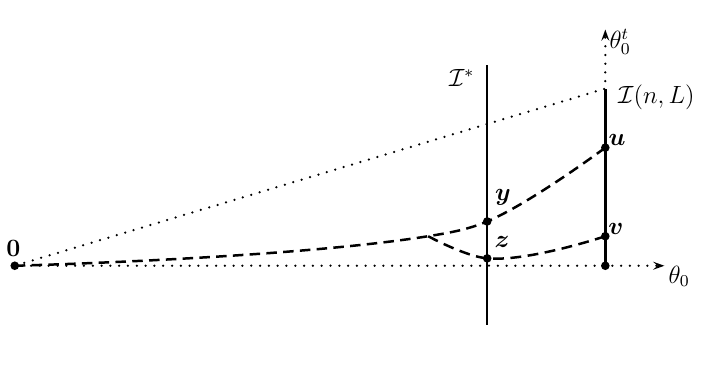}
        \caption{Illustration for Theorem~\ref{thm:logupinc}
        under case II: distance of $\segment(n,L)$ from
        $\origin$ in $\thn$ direction is $n$, width of
        $\segment(n,L)$ in $\thnt$ direction is $L$, distance
        of $\segment^\ast$ from $\origin$ is
        $n-\Delta^{-1}L$, if $\mathbb{T}\in\event$ then
        geodesics from $\origin$ to points in $\segment(n,L)$
        passes through $\segment^\ast$.}
    \end{figure}
    Consider two points $\bu$ and $\bv$ on $\segment(n,L)$. 
    Since $\mathbb{T}\in\event$, there exist points $\by$ and
    $\bz$ on $\segment^\ast$ such that the geodesic
    $\Gamma(\origin,\bu)$ passes through $\by$, and the
    geodesic $\Gamma(\origin,\bv)$ passes through $\bz$. Then \[
        T(\origin,\bu)-T(\origin,\bv)
        \leq 
        (T(\origin,\bz)+T(\bz,\bu))
        -(T(\origin,\bz)+T(\bz,\bv))
        =T(\bz,\bu)-T(\bz,\bv).
    \]
    Similarly, we get the opposite inequality with $\by$ in place of $\bz$. Therefore,
    \begin{equation}\label{eq:logupinc4}
         |T(\origin,\bu)-T(\origin,\bv)|
    \leq \max_{\bx\in\segment^\ast} |T(\bx,\bu)-T(\bx,\bv)|.
    \end{equation}
    Fix an $\bx\in \segment^\ast$. Then
    \begin{multline}\label{eq:logupinc5}
        |T(\bx,\bu)-T(\bx,\bv)| \leq 
          |T(\bx,\bu)-h(\bu-\bx)|
        + |T(\bx,\bv)-h(\bv-\bx)|\\
        + |h(\bu-\bx)-g(\bv-\bx)| 
        + |h(\bv-\bx)-g(\bv-\bx)|
        + |g(\bu-\bx)-g(\bv-\bx)|.
    \end{multline}
    Since $\bu$, $\bv$ are in $\segment(n,L)$, we have 
    \begin{equation}\label{eq:logupinc6}
        \yprojthnot(\bu-\bv)\leq L.
    \end{equation}
    From \eqref{eq:logupinc1} we get 
    \begin{equation}\label{eq:logupinc7}
        \xprojthnot(\bu-\bx) = \Delta^{-1}(L),
    \end{equation}
    and 
    \begin{equation}\label{eq:logupinc8}
         |\yprojthnot(\bu-\bx)|
    \leq \C t^{1/2}L(\log\Delta^{-1}(L))^{1/2}.
    \end{equation}
    Combining \eqref{eq:logupinc6}, \eqref{eq:logupinc7}, 
    \eqref{eq:logupinc0} we get 
    \begin{equation}\label{eq:logupinc9}
        \ltwo{\bu-\bx}\leq\C\Delta^{-1}(L).
    \end{equation}
    Hence, by Proposition~\ref{prop:nrflog}, and using 
    $\log L$ and $\log\Delta^{-1}(L)$ are of the same order, 
    we get 
    \begin{equation}\label{eq:logupinc10} 
        |h(\bu-\bx)-g(\bu-\bx)|
    \leq\C\sigma(\Delta^{-1}(L))\log L.
    \end{equation} 
    Similarly, \eqref{eq:logupinc7}-\eqref{eq:logupinc10} 
    hold for $\bu$ replaced with $\bv$. By 
    Lemma~\ref{lem:auxgeom2} and using \eqref{eq:logupinc6}, 
    \eqref{eq:logupinc7}, \eqref{eq:logupinc8}, we get
    \begin{equation}\label{eq:logupinc11}
        |g(\bu-\bx)-g(\bv-\bx)| 
    \leq \C t^{1/2} \frac{L^2}{\Delta^{-1}(L)}\log\Delta^{-1}(L) 
    \leq \C t^{1/2} \sigma(\Delta^{-1}(L))\log L.   
    \end{equation}
    Using \eqref{eq:logupinc10} and the same for $\bu$ 
    replaced with $\bv$, \eqref{eq:logupinc11}, and 
    \eqref{eq:logupinc5}, we get 
    \begin{align*}
        & \Prob\Lp\,|T(\bx,\bu)-T(\bx,\bv)|\geq 
        t\sigma(\Delta^{-1}(L))\log L\,\Rp\\
        \leq & \Prob\Lp\,|T(\bx,\bu)-h(\bu-\bx)|\geq\Cl{321} 
        t\sigma(\Delta^{-1}(L))\log L\,\Rp\\
        & + \Prob\Lp\,|T(\bx,\bv)-h(\bv-\bx)|\geq\Cr{321} 
        t\sigma(\Delta^{-1}(L))\log L\,\Rp.
    \end{align*}
    Therefore, using \eqref{eq:logupinc9} and the same for
    $\bu$ replaced by $\bv$, and using \eqref{A1} we get 
    \begin{equation}\label{eq:logupinc12}
        \Prob\Lp\,|T(\bx,\bu)-T(\bx,\bv)|\geq t\sigma(\Delta^{-1}(L))\log L\,\Rp
        \leq \C e^{-\C t\log L}.
    \end{equation}
    The number of choices of lattice points corresponding to 
    $\bx$, $\bu$, $\bv$ is $\C t L^3 (\log L)^{1/2}$. Using 
    \eqref{eq:logupinc4}, \eqref{eq:logupinc12}, and a union 
    bound, we get 
    \[
        \Prob\Lp\,\diff(n,L)\geq 
        t\sigma(\Delta^{-1}(L))\log L,\;\mathbb{T}\in\event\,\Rp
        \leq \C e^{-\C t\log L}.
    \]
    This concludes the proof because we already found in \eqref{eq:logupinc3} that $\event^c$ has appropriately small probability. 
\end{proof} 

We need the following variation of the last result. In the
last result we chose a direction of both type I and type II,
and considered the transverse increment over a segment with
one endpoint having that chosen direction. In the next
result, we consider transverse increment over a segment whose
endpoints have direction in a neighborhood of a fixed
direction which is of both type I and type II.

\begin{corollary}\label{cor:logupinc}
    Let $\thn$ be a direction of both type I and II. For $n>0$, $L>0$, $d$, let 
    \[
        \segment(n,L,d):=\LP\,\bx\in\RR^2\,:\,
        \xprojthnot(\bx)=n,\;d\leq\yprojthnot(\bx)\leq d+L\,\RP,
    \]
    and 
    \[
        \diff(n,L,d):=\max\LP\,
        |T(\origin,\bx)-T(\origin,\by)|
        \,:\,\bx,\by\in\segment(n,L,d)\,\RP.
    \]
    Then, under the Assumptions~\ref{asm:sigmaexp} and 
    \ref{asm:sigmareg}, there exist positive constants 
    $\delta_1$, $\delta_2$ $C_1$, $C_2$, $C_3$, $L_0$, $n_0$,
    $t_0$, such that for $L\geq L_0$, $n\geq n_0$, 
    $t\geq t_0$, $|d|\leq n\delta_1$, 
    $L\leq\delta_2\Delta(n)$, we have  
    \[
        \Prob\Lp\,
        \diff(n,L,d) \geq C_3 L\frac{|d|}{n} +
        t\sigma(\Delta^{-1}(L))\log L
        \,\Rp
        \leq C_1 e^{-C_2 t\log L}.
    \]
\end{corollary}

\resetconstant

\begin{proof}   
    Fix $\delta_1>0$, $\delta_2>0$, to be assumed
    appropriately small whenever required. Fix $n_0>0$,
    $L_0>0$, $t_0>0$, to be assumed appropriately large
    whenever required. Consider $n$, $L$, $t$, $d$ such that
    $L\geq L_0$, $n\geq n_0$, $|d|\leq n\delta_1$,
    $L\leq\delta_2\Delta(n)$. Let $\bu:=n\uthn+d\uthnt$,
    $\bv:=n\uthn+(d+L)\uthnt$. Let $\theta_1$, $\theta_2$ be
    the directions of $\bu$ and $\bv$ respectively. Let $\bw$
    be the projection of $\bu$ on the line joining $\origin$
    and $\bv$ in direction $\theta_1^t$ which exists assuming
    $\delta_1$ is small enough. Let $\segment^\ast$ be 
    segment joining $\bu$ and $\bw$.
    \begin{figure}[H]
    \centering
        \includegraphics[width=0.65\linewidth]{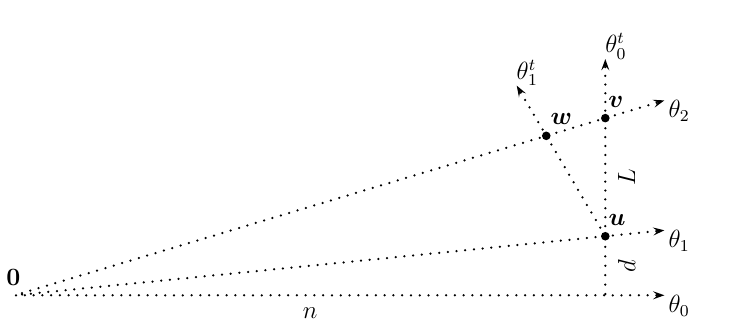}
        \caption{Illustration for 
        Corollary~\ref{cor:logupinc}. The segment joining
        $\bu$ and $\bv$ is $\segment(n,L,d)$. The segment 
        joining $\bu$ and $\bw$ is $\segment^\ast$.}
    \end{figure}%
    Let 
    \[
        \diff^\ast:=
        \max\LP\,|T(\origin,\bx)-T(\origin,\by)|\,:\,
        \bx,\by\in\segment^\ast\,\RP. 
    \]
    Let $|\segment^\ast|$ be the length of $\segment^\ast$. 
    To bound $\diff^\ast$ we use Theorem~\ref{thm:logupinc}
    with the following variables
    \[
        \tilde{\thn}=\theta_1,
        \quad
        \tilde{n}:=\ltwo{\bu},
        \quad
        \tilde{L}:=|\segment^\ast|,
        \quad
        \tilde{t}=\frac{t}{2},
        \quad
        \tilde{\eta}=1.
    \]
    Since $\thn$ is of both type I and type II, assuming
    $\delta_1$ is small enough, we get all possible values of
    $\theta_1$ are uniformly of both type I and type II i.e.,
    they satisfy the curvature conditions with same
    constants. Thus the condition on $\theta_1$ for applying
    Theorem~\ref{thm:logupinc} holds. Assuming $\delta_1$ and
    $\delta_2$ are small enough, we get 
    $\C L\leq|\segment^\ast|\leq \Cl{logupinc1} L$, 
    $\C n\leq\tilde{n}\leq\C n$. 
    Therefore, $\tilde{L}\geq\tilde{L_0}$, 
    $\tilde{n}\geq\tilde{n_0}$ hold assuming $L_0$ and $n_0$ 
    are large enough. Also from $L\leq\delta_2\Delta(n)$ we 
    get $\tilde{L}\leq\Delta(\tilde{n})$.
    Hence all the conditions for applying 
    Theorem~\ref{thm:logupinc} hold and we get 
    \begin{equation}\label{eq:corlogupinc1}
        \Prob\Lp\,
        \diff^\ast\geq\frac{t}{2}\sigma(\Delta^{-1}
        (L))\log L
        \,\Rp
        \leq\C e^{-\C t \log L}.
    \end{equation}
    Now let us consider the difference 
    $|\diff^\ast-\diff(n,L,d)|$.
    Considering the triangle with vertices $\bu$, $\bv$, 
    $\bw$, we get
    \begin{equation}\label{eq:corlogupinc2}
     \ltwo{\bv-\bw}
    =L\frac{|\sin(\theta_1^t-\thnt)|}
        {|\sin(\theta_1^t-\theta_2)|}.
    \end{equation}
    Assuming $\delta_1$ and $\delta_2$ are small enough, we 
    get  
    \begin{equation}\label{eq:corlogupinc3}
        |\sin(\theta_1^t-\theta_2)|\geq
        \C|\sin(\thnt-\thn)|,
    \end{equation}
    and
    \begin{multline}\label{eq:corlogupinc4}
        |\sin(\tht-\thnt)| 
        \leq\C|\theta_1^t-\thnt|
        \leq\C|\theta_1-\thn|\\
        \leq\C|\sin(\theta_1-\thn)|
        \leq\C\frac{|d|}{n}|\sin(\thnt-\theta)|
        \leq\C\frac{|d|}{n}|\sin(\thnt-\thn)|,
    \end{multline}
    where the second inequality holds by 
    Remark~\ref{remark:linear}. Combining 
    \eqref{eq:corlogupinc2}-\eqref{eq:corlogupinc4} we get
    \[
        \ltwo{\bv-\bw}
    \leq\Cl{curvup1}L\frac{|d|}{n}.
    \]
    Therefore, if $\bx$ is a point on $\segment(n,L,d)$ and 
    $\by$ is its projection on $\segment^\ast$ in direction 
    $\theta_2$, then $\ltwo{\bx-\by}\leq\C L|d|/n$. Since $h$
    is subadditive and therefore sublinear, we get 
    $h(\bx-\by)\leq\Cl{317} L|d|/n$. Assuming $\delta_1<1$,
    taking $L_0$ large enough, and using \eqref{A2}, we get
    $|d|/n\leq\delta_1<1\leq\Delta^{-1}(L)/L$, so that
    $L|d|/n\leq\Delta^{-1}(L)$. Therefore, using \eqref{A1} 
    we get for all $t^\prime>0$ 
    \[
    \Prob\Lp\,
        T(\bx,\by)\geq
        \Cr{317}L\frac{|d|}{n}+t^\prime
        \sigma(\Delta^{-1}(L))
    \,\Rp
    \leq\C e^{-\C t^\prime}.
    \]
    Using $|\segment^\ast|\leq\Cr{logupinc1} L$, 
    $|\segment(n,L,d)|=L$, and using a union bound, we get
    \[
        \Prob\Lp\,
            |\diff^\ast-\diff(n,L)|
            \geq\Cr{317}
            L\frac{|d|}{n}+t^\prime\sigma(\Delta^{-1}(L))
        \,\Rp
        \leq \C L^2 \exp\Lp -\C t^\prime\Rp.
    \]
    Therefore, taking $t^\prime=(t/2)\log L$, and assuming
    $t_0$ and $L_0$ are large enough, we get
    \begin{equation}\label{eq:corlogupinc5}
        \Prob\Lp\,
        |\diff^\ast-\diff(n,L)|
        \geq\Cr{317} L\frac{|d|}{n}+\frac{t}{2}\sigma(\Delta^{-1}(L))
        \log L\,\Rp
    \leq \C \exp\Lp -\C t\log L\Rp.
    \end{equation}
    Combining \eqref{eq:corlogupinc1} and 
    \eqref{eq:corlogupinc5} completes the proof.
\end{proof} 

\section{Refined upper bound on nonrandom fluctuations}\label{sec:ghloglog}

Theorem~\ref{thm:logupinc} of the previous section provides a
preliminary upper bound of the fluctuations of the transverse
increments. To prove the refined bound of
Theorem~\ref{thm:loglogupinc}, we need to reduce the
correction factor in Proposition~\ref{prop:nrflog} from
$\log$ to fixed but arbitrary small power of $\log$. Related
results are known in the literature. In \cite{Tessera2018}
it has been shown that $(\log x)^{1/2}$ is a valid correction
factor in FPP on Cayley Graphs on integer lattices. In 
\cite{DamronWang} it have shown that any iterate of $\log$ is
a valid correction factor in a spherically symmetric model of
FPP. In \cite{KenUniform} a bound with only a constant correction factor has been shown in a spherically symmetric model. In our case, to reduce the correction factor, we use a modified 
version of the procedure of \cite{Ken97}. We now introduce
the concept of convex hull approximation property from
\cite{Ken97}.

\begin{notation}\label{Notation:CHAP}
    Let $\Sector_0$ be the set of directions where the boundary of $\Bb$ is differentiable. Consider $\bx$ with direction in $\Sector_0$. Let $H_{\bx}$ be the tangent to $\partial g(\bx)\Bb$ at $\bx$. Let $H_{\bx}^0$ be the line through $\origin$ parallel to $H_{\bx}$. Let $g_{\bx}$ be the unique linear functional on $\RR^2$ satisfying $g_{\bx}(\by)=0$ for all $y\in H_{\bx}^0$, and $g_{\bx}(\bx)=g(\bx)$. Recall $\Phi$ from Notation~\ref{notn:functionclass}. Define for $\phi\in\Phi$, $\nu\geq 0$, $C>0$, $K>0$,
    \[
        Q_{\bx}(\nu,\phi,C,K):=
        \LP\by\in\ZZ^2:\ltwo{\by}\leq K\ltwo{\bx}, g_{\bx}(\by)\leq g(\bx), h(\by)\leq g_{\bx}(\by) + C\ltwo{\bx}^\nu\phi\Lp\ltwo{\bx}\Rp\RP.
    \]
\end{notation}

\begin{definition} 
    We say that $h$ satisfies the \textit{convex hull 
    approximation property} (CHAP) with exponent $\nu\geq 0$ 
    and correction factor $\phi\in\Phi$ in a set of 
    directions $\Sector\subset\Sector_0$, if there exist 
    constants $M>0$, $C>0$, $K>0$, $a>1$ such that 
    $x/\alpha\in\textsf{Co}(Q_{\bx}(\nu,\phi,C,K))$ for some 
    $\alpha\in[1,a]$, for all $\bx\in\QQ^2$ with 
    $\ltwo{\bx}\geq M$ and direction of $\bx$ in $\Sector$,
    where $\textsf{Co}$ denotes the convex hull. When we want
    to specify the specific constants, we say $h$ satisfies 
    $\text{CHAP}(\nu,\phi,M,C,K,a)$ in sector 
    $\Sector$.
\end{definition}

The procedure of \cite{Ken97} in our terminology as follows.
The objective of \cite{Ken97} was to prove GAP with exponent 
$\alpha$ and correction factor $\phi_1$. To achieve this, 
first, it is shown that CHAP with exponent $\alpha$ and 
correction factor $\phi_1$ holds. This is done in an 
iterative way. GAP with exponent $1$ and correction factor 
$\phi_1$ holds trivially because $h$ is sublinear. Then the 
exponent of GAP is reduced from $1$ to $\alpha$ iteratively 
using CHAP with exponent $\alpha$. Here we are not concerned 
about the exponent $\alpha$. 

In contrast to \cite{Ken97}, here we want to change the 
correction factor of GAP from $\phi_1$ to $\phi_\eta$ for 
some small $\eta>0$ while keeping the exponent $\alpha$ 
unchanged. We do this in two steps. In the first step, we 
prove that CHAP holds with exponent $\alpha$ and correction 
factor $\widehat{\phi}$ (recall $\widehat{\phi}$ from 
Notation~\ref{notn:functionclass}). From 
Proposition~\ref{prop:nrflog} we get $h$ satisfies GAP with 
exponent $\alpha$ and correction factor $\phi_1$.  Using 
this, we reduce the correction factor of GAP from $\phi_1$ to
$\phi_\eta$, which is the second step. The result on CHAP
with exponent $\alpha$ and correction factor $\widehat{\phi}$
is the following.

\begin{theorem}\label{thm:chaploglog}
    Let $\thn$ be a direction of both type I and type II.
    Then, under Assumptions~\ref{asm:sigmaexp} and 
    \ref{asm:sigmareg}, there exists $\delta>0$ such that 
    CHAP holds for $h$ with exponent $\alpha$ and correction 
    factor $\widehat{\phi}$ in the sector 
    $(\thn-\delta,\thn+\delta)$.
\end{theorem} 

After carrying out the second step we get the refined upper 
bound on the nonrandom fluctuations stated below. 

\begin{theorem}\label{thm:nrf2}
    Let $\thn$ be a direction of both type I and type II.
    Fix $\eta\in(0,1]$. Then, under 
    Assumptions~\ref{asm:sigmaexp} and \ref{asm:sigmareg}, 
    there exist constants $C>0$, $M>0$, $\delta>0$ such that 
    for all $\ltwo{\bx}\geq M$ with direction of $\bx$ in 
    $(\thn-\delta,\thn+\delta)$, we have
    \[
        h(\bx)
        \leq g(\bx)+C\sigma(\ltwo{\bx})(\log\ltwo{\bx})^\eta,
    \]
    i.e., $\text{GAP}(\alpha,\phi_\eta,C,M)$ holds in the sector $(\thn-\delta,\thn+\delta)$.
\end{theorem} 

In the first subsection below we prove 
Theorem~\ref{thm:chaploglog}, then in the second subsection 
we prove Theorem~\ref{thm:nrf2} using 
Theorem~\ref{thm:chaploglog}.

\subsection{Proof of Theorem~\ref{thm:chaploglog}}

\resetconstant

First, we choose $\delta>0$ which is fixed throughout the
proof of Theorem~\ref{thm:chaploglog}. Since $\thn$ is of 
both type I and II, we choose $\delta>0$ such that 
$\partial\Bb$ is differentiable in 
$(\thn-2\delta,\thn+2\delta)$ and there exists $\delta_1>0$ 
such that for all $\theta\in(\thn-2\delta,\thn+2\delta)$ and 
$|\delta_2|\leq\delta_1$
\begin{equation}\label{eq:deltachoice}
    \C\delta_2^2\leq g(\uth+\delta_2\utht)-g(\uth) \leq\C\delta_2^2.
\end{equation}
So all $\theta\in(\thn-\delta,\thn+\delta)$ are of both type 
I and type II with same constants.
This allows us to use results which hold in type I or type II
directions with same constants for all 
$\theta\in(\thn-\delta,\thn+\delta)$.

We extract a sufficient condition from \cite{Ken97} for the 
Theorem~\ref{thm:chaploglog} to hold. We state the condition
as Proposition~\ref{prop:sufficient}. Since it is essentially
proved in \cite{Ken97}, we do not prove it here. To state the
condition we need the concept of skeletons of paths defined 
below.

\vspace{\baselineskip}

\paragraph{\textbf{Construction of fine skeletons:}} For 
$\bx\in\RR^2$, $n>0$, $\lambda>0$, $K>0$, the \textit{fine 
$Q_{\bx}(\alpha,\widehat{\phi},\lambda,K)$-skeleton} of a 
self-avoiding path $\gamma$ from $\origin$ to $n\bx$ is the 
sequence of marked points $\bv_0,\dots,\bv_m$ on $\gamma$ 
constructed as follows. Let $\bv_0:=\origin$, and given 
$\bv_i$, let $\bv^\prime_{i+1}$ be the first point (if any) 
in $\gamma$ such that $\bv^\prime_{i+1}-\bv_i\not\in 
Q_{\bx}(\alpha,\widehat{\phi},\lambda,K)$; then let 
$\bv_{i+1}$ be the last lattice point in $\gamma$ before 
$\bv^\prime_{i+1}$ if $\bv^\prime_{i+1}$ exists; otherwise 
let $\bv_{i+1}=\lfloor n\bx\rfloor$ and end the construction.

\begin{proposition}\label{prop:sufficient}
    Consider an infinite sequence of i.i.d.\ copies of the 
    passage-time configuration $(\hat{T}^i)_{i=0}^\infty$ on
    the lattice. Suppose for some positive constants 
    $\lambda_1$, $\lambda_2$, $\lambda_3$, we have 
    \begin{align*}
        &\Prob\left(\,
        \sum_{i=0}^{m-1}
        \LT
        h(\bv_i-\bv_i+1)-\hat{T}^i(\bv_i,\bv_{i+1})
        \RT
        \geq
        \frac{\lambda_1}{16}m
            \sigma(\ltwo{\bx})\log\log\ltwo{\bx}\right.\\  
        &\left.\qquad\text{for some $m\geq 1$ and some 
        $Q_{\bx}(\alpha,\widehat{\phi},\lambda_1,5)$-skeleton
        $(\bv_j)_{j=0}^m$ of a 
        path}
        \vphantom{\sum_{i=0}^{m-1}}\,\right)\\
        \leq & \exp\Lp-\lambda_2\log\log\ltwo{\bx}\Rp
        \numberthis\label{eq:nrf23}
    \end{align*}
    for all $\bx$ with $\ltwo{\bx}\geq\lambda_3$ and
    direction of $\bx\in(\thn-\delta,\thn+\delta)$. Then
    there exists $M>0$ such that $h$ satisfies
    $\text{CHAP}(\alpha,\widehat{\phi},M,\lambda_1,4,2)$ in
    the sector $(\thn-\delta,\thn+\delta)$. 
\end{proposition}

In order to verify \eqref{eq:nrf23} we need the concept of
`coarse skeletons.'

\begin{figure}[H]
    \centering
    \includegraphics[width=0.55\linewidth]{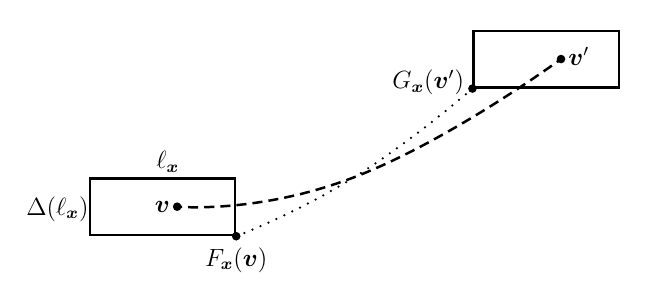}
    \caption{Construction of the coarse skeleton: for every 
    pair of consecutive points $\bv$, $\bv^\prime$ in a fine
    skeleton of a path, we have $F_{\bx}(\bv)$ and
    $G_{\bx}(\bv^\prime)$ as consecutive points in the coarse
    skeleton of the path.}
    \label{Fig:6}
\end{figure}

\paragraph{\textbf{Construction of coarse skeletons:}}
Consider
$\bx\in\RR^2$ with direction $\theta\in(\thn-\delta,\thn+\delta)$. Define
\begin{equation}\label{eq:coarsedef1}
    \ell_{\bx}:=\ltwo{\bx}(\log\ltwo{\bx})^{-2/\alpha}.    
\end{equation}
Also, for $i,j\in\mathbb Z$, define
\begin{equation}\label{eq:coarsedef2}
    B_{ij}:=\LP\,\by\in\RR^2\,:\,
    \xprojth(\by)\in
    [i\ell_{\bx},(i+1)\ell_{\bx}),\,
    \yprojth(\by)\in[j\Delta(\ell_{\bx}),
    (j+1)\Delta(\ell_{\bx}))
    \,\RP,
\end{equation}
where $\alpha$ is defined in Assumption~\ref{asm:sigmareg}. 
So $B_{ij}$ is a parallelogram with side lengths $\ell_{\bx}$
and $\Delta(\ell_{\bx})$, and these parallelograms cover the
whole plane. Given $\bv\in B_{ij}$, let 
\[
    F_{\bx}(\bv):=
    \lfloor i\ell_{\bx}\uth+j\Delta(\ell_{\bx})\utht
    \rfloor,\quad
    G_{\bx}(\bv):=
    \lfloor (i+1)\ell_{\bx}\uth+j\Delta(\ell_{\bx})\utht
    \rfloor.
\]
So, $F_{\bx}(\bv)$ is the lattice point corresponding to 
down-left corner of the parallelogram $B_{ij}$ containing 
$\bv$ and $G_{\bx}(\bv)$ is the down-right corner, see 
Figure~\ref{Fig:6}. Suppose $(\bv_i)_{i=0}^m$ is a fine 
$Q_{\bx}(\alpha,\widehat{\phi},\lambda,5)$-skeleton of some 
path for some $\lambda>0$. Then its \textit{coarse skeleton} 
$(\bw_j)_{j=0}^{2m}$ is defined as follows. For 
$0\leq i\leq m-1$, let $\bw_{2i}:=F_{\bx}(\bv_i)$ and 
$\bw_{2i+1}:=G_{\bx}(\bv_{i+1})$. 

\begin{remark}\label{remark:coarse} 
    If $(\bv_i)_{i=0}^m$ is a fine 
    $Q_{\bx}(\alpha,\widehat{\phi},\lambda,5)$ skeleton of 
    some path and $(\bw_j)_{j=0}^{2m}$ the corresponding 
    coarse skeleton, then  $\ltwo{\bv_i-\bv_{i+1}}\leq 
    5\ltwo{\bx}$ and $\ltwo{\bw_{2i-1}-\bw_{2i}}\leq 
    6\ltwo{\bx}$ for large enough $\ltwo{\bx}$.
\end{remark} 

We state two propositions which in combination establishes 
\eqref{eq:nrf23}. We define a few constants first. Let 
$\Cl{0.1}$, $\Cl{0.2}$, $\Cl{1}$, $\Cl{2}$, $\Cl{0}$, be 
positive constants such that for all $\bx\in\RR^2$ with 
$\ltwo{\bx}\geq\Cr{0.1}$, and for all $\bu,\bv\in\RR^2$ with 
$\Cr{0.2}\leq|\bu-\bv|\leq 6\ltwo{\bx}$, we have 
\begin{equation}\label{eq:defineczero}
    \Prob\Lp\,|T(\bu,\bv)-h(\bu-\bv)|\geq t\sigma(\ltwo{\bx})\,\Rp
    \leq\Cr{1} e^{-\Cr{2} t},
\end{equation}
and 
\begin{equation}\label{eq:defineczero2}
    \Cr{0}:=256\alpha^{-1}(1+\Cr{1})/\Cr{2}.
\end{equation}
For all $\bx\in\RR^2$, we use the shorthand notation 
\[
    Q_{\bx}:=Q_{\bx}\Lp\alpha,\widehat{\phi},\Cr{0},5\Rp.
\]

\begin{proposition}\label{prop:coarseentropy}
    Under the assumptions of Theorem~\ref{thm:chaploglog}, 
    there exist positive constants $\Cl{3}$, $\Cl{4}$ such
    that for $\bx$ with direction in
    $(\thn-\delta,\thn+\delta)$ and $\ltwo{\bx}\geq\Cr{3}$,
    and for all $m\geq 1$, we have
    \begin{align*}
        & \Prob\left(\,\sum_{i=0}^{m-1}\LT
        h(\bw_{2i}-\bw_{2i+1})-\hat{T}^i(\bw_{2i},\bw_{2i+1})
        \RT
        >
        \frac{\Cr{0}}{32}m
        \sigma\Lp\ltwo{\bx}\Rp\log\log\ltwo{\bx}\right.\\ 
        &\left.\qquad\text{for some coarse $Q_{\bx}$-skeleton
        $(\bw_j)_{j=0}^{2m}$}\vphantom{\sum_{i=0}^{m-1}}\Rp\\
        \leq 
        & \exp\Lp-m \Cr{4}\log\log\ltwo{\bx}\,\right).
    \end{align*}
\end{proposition}

\begin{proposition}\label{prop:finecoarsediffsum}
    Under the assumptions of Theorem~\ref{thm:chaploglog}, 
    there exist positive constants $\Cl{5}$, $\Cl{6}$ such
    that for $\bx$ with direction in 
    $(\thn-\delta,\thn+\delta)$ and 
    $\ltwo{\bx}\geq\Cr{5}$, and for all $m\geq 1$, we have 
    \begin{align*}
        & \Prob\left(\,
        \sum_{i=0}^{m-1}
        \LT 
        h(\bv_i-\bv_{i+1})-\hat{T}^i(\bv_i,\bv_{i+1})-
        h(\bw_{2i}-\bw_{2i+1})+\hat{T}^i(\bw_{2i},\bw_{2i+1})
        \RT
        \right.\\
        &\left.\qquad\geq\frac{\Cr{0}}{32} 
        m\sigma(\ltwo{\bx})\log\log\ltwo{\bx}
        \text{ for some fine $Q_{\bx}$-skeleton 
        $(\bv_j)_{j=0}^m$}\right.\\ 
        &\left.\qquad\text{and the corresponding coarse 
        skeleton 
        $(\bw_j)_{j=0}^{2m}$}
        \vphantom{\sum_{i=0}^{m-1}}\,\right)
        \\
        \leq & \exp\Lp-m \Cr{6}\log\log\ltwo{\bx}\Rp.
    \end{align*}
\end{proposition} 

So \eqref{eq:nrf23} holds by
Propositions~\ref{prop:finecoarsediffsum} and
\ref{prop:coarseentropy}. Hence, to complete the proof of 
Theorem~\ref{thm:chaploglog}, we only need to prove these two
propositions. We need the following lemma first.

\begin{lemma}\label{lem:width}
    Under the assumptions of Theorem~\ref{thm:chaploglog},
    there exist positive constants $\Cl{61}$, $\Cl{7}$ such
    that for $\bx$ with direction in
    $(\thn-\delta,\thn+\delta)$ and $\ltwo{\bx}\geq\Cr{61}$
    we have the following.
    \begin{enumerate}[(i)]
    \item For all $\by\in Q_{\bx}$, $|\yprojth(\by)|\leq 
    \Cr{7}\Delta(\ltwo{\bx})(\log\log\ltwo{\bx})^{1/2}$,
    where $\theta$ is the direction of $\bx$. 
    \item The number of coarse $Q_{\bx}$-skeletons of length
    $2m+1$ is at most $(\log\ltwo{\bx})^{4m/\alpha}$, where
    $\alpha$ is defined in \eqref{A2}.
    \end{enumerate}
\end{lemma}

\begin{proof} 
    Fix $\bx$ with direction $\theta$ in
    $(\thn-\delta,\thn+\delta)$. We assume that $\ltwo{\bx}$ 
    is large enough whenever required.
    
    \vspace{\baselineskip}
    
    \textit{Proof of (i):} Consider $\by\in Q_{\bx}$. Then 
    $g(\by) \leq h(\by) \leq g_{\bx}(\by)
    + \Cr{0}\sigma(\ltwo{\bx})\log\log\ltwo{\bx}$,
    so 
    \begin{equation}\label{eq:l551}
        g(\by)-g_{\bx}(\by)
    \leq\Cr{0}\sigma(\ltwo{\bx})\log\log\ltwo{\bx}.
    \end{equation}
    Consider three cases.
    
    \vspace{\baselineskip}
    
    \textit{Case I:} Suppose $\xprojth(\by)>0$ and
    $|\yprojth(\by)|\leq\delta_1\xprojth(\by)$, where 
    $\delta_1$ is defined in \eqref{eq:deltachoice}. 
    From \eqref{eq:deltachoice} we get 
    \begin{equation}\label{eq:l552}
        g(\by)-g_{\bx}(\by) = 
            g(\xprojth(\by)\uth+\yprojth(\by)\utht)
                - g(\xprojth(\by)\uth)
        \geq \C \frac{\yprojth(\by)^2}{\xprojth(\by)}. 
    \end{equation}
    Since $\by\in Q_{\bx}$, we have 
    \[
        \xprojth(\by)\leq\C\ltwo{\by}\leq\C\ltwo{\bx}.
    \]
    This with \eqref{eq:l552} and \eqref{eq:l551} implies 
    \[
        |\yprojth(\by)|
        \leq\C\Delta(\ltwo{\bx})(\log\log\ltwo{\bx})^{1/2}.
    \]
    
    \vspace{\baselineskip}
    
    \textit{Case II:} Now suppose $\xprojth(\by)>0$ and 
    $|\yprojth(\by)|\geq\delta_1\xprojth(\by)$. 
    Let us consider $\yprojth(\by)>0$, the case $\yprojth(\by)<0$ is similar.
    Using convexity of $g$ and \eqref{eq:deltachoice} we get
    \begin{align*}
    & g(\by)-g_{\bx}(\by)\\ 
    = & g(\xprojth(\by)\uth+\yprojth(\by)\utht)-g(\xprojth(\by)\uth)\\
    \geq & \frac{g(\xprojth(\by)\uth+\delta_1\xprojth(\by)\utht)-g(\xprojth(\by)\uth)}{\delta_1\xprojth(\by)/\yprojth(\by)}\\
    = & \yprojth(\by)\delta_1^{-1}(g(\uth+\delta_1\utht)-g(\uth))\\
    \geq & \C\yprojth(\by).
    \end{align*}
    Hence using \eqref{eq:l551} and \eqref{A2} we get
    \[
        \yprojth(\by)
    \leq \C\sigma(\ltwo{\bx})\log\log\ltwo{\bx}
    \leq \C\Delta(\ltwo{\bx})(\log\log\ltwo{\bx})^{1/2}.
    \]
    
    \vspace{\baselineskip}
    
    \textit{Case III:} Suppose $\xprojth(\by)<0$. 
    Then $g_{\bx}(\by)<0$. Using $\by\in Q_{\bx}$, 
    \eqref{eq:l551}, and \eqref{A2} we get  
    \[
        \yprojth(\by)
        \leq \C \ltwo{\by}
        \leq \C g(\by)
        \leq \C \sigma(\ltwo{\bx})\log\log\ltwo{\bx}
        \leq \C \Delta(\ltwo{\bx})(\log\log\ltwo{\bx})^{1/2}.
    \]
    
    \vspace{\baselineskip}
    
    \textit{Proof of (ii):} Given any $\bv$ let 
    $Q_{\bx}(\bv)$ denote the translate of $Q_{\bx}$ by
    $\bv$. Suppose part of a coarse skeleton is given as
    $(\bw_0,\dots,\bw_{2i})$. We find an upper bound on
    number of possibilities of $(\bw_{2i+1},\bw_{2i+2})$.
    Consider a fine skeleton $(\bv_0,\dots,\bv_i)$
    corresponding to $(\bw_0,\dots,\bw_{2i})$. Since
    $(\bw_{2i-1},\bw_{2i})$ is fixed, all choices of $\bv_i$
    lie in the same $B_{i_0j_0}$. Consider the union of
    $Q_{\bx}(\bv)$ over all $\bv\in B_{i_0 j_0}$. Each
    $Q_{\bx}(\bv)$ is contained in a parallelogram of length
    $\C\ltwo{\bx}$ in $\theta$ direction and length
    $\C\Delta(\ltwo{\bx})(\log\log\ltwo{\bx})^{1/2}$ in
    $\tht$ direction. Hence, the union of all such
    $Q_{\bx}(\bv)$ as $\bv$ varies in a parallelogram 
    $B_{i_0 j_0}$ is contained in a parallelogram of length 
    $\C\ltwo{\bx}$ in $\theta$ direction and length
    $\C\Delta(\ltwo{\bx})(\log\log\ltwo{\bx})^{1/2}$ in 
    $\tht$ direction (using (i)). From \eqref{eq:coarsedef1} 
    we have
    \[
        \frac{\Delta(\ltwo{\bx})(\log\log\ltwo{\bx})^{1/2}}
        {\Delta(\ell_{\bx})}
        \leq\C
        (\log\ltwo{\bx})^{(1+\beta)/\alpha}
        (\log\log\ltwo{\bx})^{1/2},
    \]
    and,
    \[
        \frac{\ltwo{\bx}}{\ell_{\bx}}
        =(\log\ltwo{\bx})^{2/\alpha}.
    \]
    Therefore, using $\beta<1$, the number of parallelograms
    $B_{ij}$ that cover the aforementioned union is at most
    $(\log\ltwo{\bx})^{4/\alpha}$. Hence, the number of
    choices of $(\bw_{2i+1},\bw_{2i+2})$ is at most
    $(\log\ltwo{\bx})^{4/\alpha}$. Iterating this from $i=0$
    to $m$, we get the result.
\end{proof} 

\subsubsection{Proof of Proposition~\ref{prop:coarseentropy}}

Fix a point $\bx\in\RR^2$ with direction 
$\theta\in(\thn-\delta,\thn+\delta)$. 
We assume $\ltwo{\bx}$ is large enough whenever required. 
Fix a coarse $Q_{\bx}$-skeleton $(\bw_j)_{j=0}^{2m}$ for some
$m\geq 1$. By Remark~\ref{remark:coarse} and
equation~\eqref{eq:defineczero} we get 
\[
    \Prob\Lp\,
    h(\bw_{2i}-\bw_{2i+1})-\hat{T}^i(\bw_{2i},\bw_{2i+1})
    \geq t\sigma\ltwo{\bx}
    \,\Rp
    \leq\Cr{1} e^{-\Cr{2}t}.
\]
For $\Cl{sp}:=\Cr{2}/(1+\Cr{1})$ we get 
\[
    \Exp\exp\Lp\,\Cr{sp}\Lp h(\bw_{2i}-\bw_{2i+1})-\hat{T}^i(\bw_{2i},\bw_{2i+1})
    \Rp^+/\sigma\Lp\ltwo{\bx}\Rp\,\Rp\leq 2.
\]
Using the independence of $\hat{T}^i$s we get
\[
    \Exp\exp\Lp\,\Cr{sp}
    \sum_{i=0}^{m-1}
    \Lp h(\bw_{2i}-\bw_{2i+1})-\hat{T}^i(\bw_{2i},\bw_{2i+1})\Rp^+/\sigma\Lp\ltwo{\bx}\Rp\,\Rp
    \leq 2^m.
\]
Hence, for all $t>0$, we have
\[
    \Prob\Lp\,\sum_{i=0}^{m-1} h(\bw_{2i}-\bw_{2i+1})-\hat{T}^i(\bw_{2i},\bw_{2i+1})
    >t m\sigma\Lp\ltwo{\bx}\Rp\log\log\ltwo{\bx}\,\Rp
    \leq 2^m e^{-\Cr{sp} m t\log\log\ltwo{\bx}}.
\]
From Lemma~\ref{lem:width} we get that the number of coarse
skeletons of length $2m+1$ is 
$\leq (\log\ltwo{\bx})^{m(4/\alpha)}$. Therefore, by
\eqref{eq:defineczero2} we get 
\begin{align*}
    & \Prob\left(\,\sum_{i=0}^{m-1}
    \LT 
    h(\bw_{2i}-\bw_{2i+1})-\hat{T}^i(\bw_{2i},\bw_{2i+1})
    \RT>
    \frac{\Cr{0}}{32}m
    \sigma\Lp\ltwo{\bx}\Rp\log\log\ltwo{\bx}\right.
    \\ 
    & \left.\quad\text{for some coarse $Q_{\bx}$-skeleton $(\bw_j)_{j=0}^{2m}$ of a path}
    \vphantom{\sum_{i=0}^{m-1}}\,\right)
    \\
    \leq & 2^m(\log\ltwo{\bx})^{m(4/\alpha)}
    \exp\Lp-\Cr{sp}\frac{\Cr{0}}{32}m\log\log\ltwo{\bx}\Rp
    \\
    \leq & \exp\Lp-m\C\log\log\ltwo{\bx}\Rp.
\end{align*}
This completes the proof of
Proposition~\ref{prop:coarseentropy}.

\subsubsection{Proof of Proposition~\ref{prop:finecoarsediffsum}}

Fix a point $\bx\in\RR^2$ with direction $\theta\in(\thn-\delta,\thn+\delta)$. We assume $\ltwo{\bx}$ is large enough whenever required. Fix a coarse $Q_{\bx}$-skeleton $(\bw_j)_{j=0}^{2m}$ for some $m\geq 1$. For $0\leq i \leq m-1$, define the set
\[
    V_i := \LP\,(\bv,\bv^\prime)\in\ZZ^2\times\ZZ^2\,:\, F_{\bx}(\bv)=\bw_{2i},\,
    G_{\bx}(\bv^\prime)=\bw_{2i+1},\,
    \bv^\prime-\bv\in Q_{\bx}\,\RP,
\]
and also define
\[
    X_i := 
    \max_{(\bv,\bv^\prime)\in V_i}
    \frac{h(\bv-\bv^\prime)-\hat{T}^i(\bv,\bv^\prime)-
    h(\bw_{2i}-\bw_{2i+1})+\hat{T}^i(\bw_{2i},\bw_{2i+1})}
    {\sigma(\ltwo{\bx})}.
\]
From the definition of $Q_{\bx}$ we have 
$\ltwo{\bv-\bv^\prime}\leq 5\ltwo{\bx}$ for all 
$0\leq i\leq m-1$ and $(\bv,\bv^\prime)\in V_i$. Therefore, 
the number of elements in $V_i$ is at most 
$\Cl{p440}\ltwo{\bx}^4$. Hence, for all $t>0$ we have
\begin{equation}\label{eq:p4.41}
    \Prob\Lp\,X_i\geq t\,\Rp
    \leq\Cl{p441}\ltwo{\bx}^4 e^{-\Cl{p442}t}.
\end{equation}
If we show that for some constant $\Cl{p443}>0$ 
\begin{equation}\label{eq:p4.42}
    \Prob\Lp\,\sum_{i=0}^{m-1} X_i\geq\frac{\Cr{0}}{32}m\log\log\ltwo{\bx}\,\Rp
    \leq e^{-\Cr{p443}m\log\ltwo{\bx}},
\end{equation}
then Proposition~\ref{prop:finecoarsediffsum} follows using Lemma~\ref{lem:width}. Therefore, we prove \eqref{eq:p4.42} now. Let 
\begin{equation}\label{eq:4.425}
    \Cr{p444}:=8\Cr{p442}^{-1},
\end{equation}
and let $N_0$, $N_1$ be positive integers such that 
\begin{align}
    2^{N_0}<\frac{\Cr{0}}{96}\log\log\ltwo{\bx}
    \leq 2^{N_0+1},
    \label{eq:p4.43.1}\\
    2^{N_1-1}<\Cl{p444}\log\ltwo{\bx}\leq 2^{N_1}.
    \label{eq:p4.43.2}
\end{align}
Then
\begin{multline}\label{eq:p4.44}
    \Prob\Lp\,
    \sum_{i=0}^{m-1} 
    X_i\geq\frac{\Cr{0}}{32}m\log\log\ltwo{\bx}
    \,\Rp
    \leq
    \Prob\Lp\,
    \sum_{i=0}^{m-1} 
    X_i\mathbb{1}(2^{N_0}\leq X_i<2^{N_1})
    \geq\frac{\Cr{0}}{96}m\log\log\ltwo{\bx}
    \,\Rp
    \\
    +\Prob\Lp\,
    \sum_{i=0}^{m-1} 
    X_i\mathbb{1}(X_i\geq 2^{N_1})
    \geq\frac{\Cr{0}}{96}m\log\log\ltwo{\bx}
    \,\Rp.
\end{multline}
For the second term in the right-hand side of 
\eqref{eq:p4.44} we have
\begin{align*}
    & \Prob\Lp\,
        \sum_{i=0}^{m-1} 
        X_i\mathbb{1}(X_i\geq 2^{N_1})
        \geq\frac{\Cr{0}}{96}m\log\log\ltwo{\bx}
    \,\Rp\\ 
    \leq 
    & \Prob\Lp\,
        \sum_{q=N_1}^{\infty}
        \sum_{i=0}^{m-1} 
        2^{q+1}
        \mathbb{1}(2^{q+1}>X_i\geq 2^q)
        \geq\frac{\Cr{0}}{96}m\log\log\ltwo{\bx}
        \sum_{q=N_1}^\infty 2^{-(q-N_1+1)}
    \,\Rp\\
    \leq 
    & \sum_{q=N_1}^\infty
    \Prob\Lp\,
        \sum_{i=0}^{m-1}
        \mathbb{1}(X_i\geq 2^q)
        \geq\frac{\Cr{0}}{96}m\log\log\ltwo{\bx}
        2^{N_1-2(q+1)}
    \,\Rp\\
    \leq &
    \sum_{q=N_1}^\infty e^{-m I(a_1(q)|b_1(q))},
    \numberthis
    \label{eq:p4.45}
\end{align*}
where
\[
    a_1(q)
    :=\frac{\Cr{0}}{96}m\log\log\ltwo{\bx}2^{N_1-2(q+1)},
    \quad 
    b_1(q)
    :=\max_{0\leq i\leq m-1}\Prob(X_i\geq 2^q), 
\]
and $I$ is the large deviation rate function for Bernoulli
random variables:
\begin{equation}\label{eq:p4.48}
    I(x|y):=x\log\frac{x}{y}+(1-x)\log\frac{1-x}{1-y}.
\end{equation}
Using \eqref{eq:p4.41}, \eqref{eq:4.425}, and 
\eqref{eq:p4.43.2}, we get for $q\geq N_1$ 
\[
    b_1(q)\leq\Cr{p441}\ltwo{\bx}^4 e^{-\Cr{p442}2^q}
    \leq e^{-\C 2^q}.
\]
So $b_1(q)$ is much smaller than $a_1(q)$. Therefore, using
\eqref{eq:p4.48} for $I(a_1(q)|b_1(q))$ and expanding the log
terms we see that the term $a_1(q)\log(b_1(q)^{-1})$
dominates the others. Hence
\[
    I(a_1(q)|b_1(q))
    \geq\C a_1(q)\log(b_1(q)^{-1})
    \geq\C\log\log\ltwo{\bx}\log\ltwo{\bx}2^{-q}. 
\]
Therefore, continuing from \eqref{eq:p4.45} and using 
\eqref{eq:p4.41} we get
\begin{equation}\label{eq:p4.46}
    \Prob\Lp\,\sum_{i=0}^{m-1} 
    X_i\mathbb{1}(X_i\geq 2^{N_1}) 
    \geq\frac{\Cr{0}}{96}m\log\log\ltwo{\bx}\,\Rp
    \leq e^{-\C m \log\ltwo{\bx}}.
\end{equation}
For the first term in the right-hand side of \eqref{eq:p4.44}
we have
\begin{align*}
    & \Prob\Lp\,\sum_{i=0}^{m-1} X_i\mathbb{1}(2^{N_0}\leq X_i<2^{N_1})\geq\frac{\Cr{0}}{96}m\log\log\ltwo{\bx}\,\Rp\\
    \leq & \Prob\Lp\,\sum_{i=0}^{m-1}\sum_{q=N_0}^{N_1-1}2^{q+1}\mathbb{1}(X_i\geq 2^q)\geq\frac{\Cr{0}}{96} m\log\log\ltwo{\bx}\,\Rp\\
    \leq & \sum_{q=N_0}^{N_1-1}\Prob\Lp\,\sum_{i=0}^{m-1}\mathbb{1}(X_i\geq 2^q)\geq\frac{\Cr{0}}{96}m(\log\log\ltwo{\bx}) 2^{-(q+1)}(N_1-N_0)^{-1}\,\Rp\\ 
    \leq & \sum_{q=N_0}^{N_1-1} e^{-m I(a_2(q)|b_2(q))},
    \numberthis\label{eq:p4.47}
\end{align*}
where 
\[
a_2(q):=\frac{\Cr{0}}{96}(\log\log\ltwo{\bx})2^{-(q+1)}(N_1-N_0)^{-1},\quad b_2(q):=\max_{0\leq i\leq m-1}\Prob(X_i\geq 2^q).
\]
We use the following claim to derive a bound of $b_2(q)$. The proof of the claim is presented later.

\begin{claim}\label{claim:finecoarsediffsum}
    For $0\leq i\leq m-1$, $(\bv,\bv^\prime)\in V_i$, and $t\in[2^{N_0},2^{N_1-1}]$, we have 
    \[
        \Prob\Lp\,
        h(\bv-\bv^\prime)-\hat{T}^i(\bv,\bv^\prime)
        -h(\bw_{2i}-\bw_{2i+1})+\hat{T}^i(\bw_{2i},\bw_{2i+1}))
        \geq t\sigma(\ltwo{\bx})\,\Rp
        \leq e^{-C_9 t\log\ltwo{\bx}}.
    \]
\end{claim}

By this claim and using that the number of elements in $V_i$ 
is at most $\Cr{p440}\ltwo{\bx}^4$ we get, for 
$N_0\leq q\leq N_1-1$
\[
    b_2(q)
    = \max_{0\leq i\leq m-1} 
    \Prob\Lp\,X_i\geq 2^q\,\Rp
    \leq \C \ltwo{\bx}^4 e^{-\C 2^q\log\ltwo{\bx}}
    \leq e^{-\C 2^q\log\ltwo{\bx}}.
\]
From \eqref{eq:p4.43.2} we get 
$\C 2^{-q}\leq a_2(q) \leq \C 2^{-q}$. Therefore, using 
\eqref{eq:p4.48} for $I(a_2(q)|b_2(q))$, and expanding the 
log terms, we see that the term $a_2(q)\log(b_2(q)^{-1})$ 
dominates the others. Hence 
\[
    I(a_2(q)|b_2(q))\geq\C a_2(q)\log(b_2(q)^{-1})\geq\C\log\ltwo{\bx}.
\]
Therefore, continuing from \eqref{eq:p4.47} we get
\[
    \Prob\Lp\,
    \sum_{i=0}^{m-1} 
    X_i\mathbb{1}(2^{N_0}\leq X_i<2^{N_1}) \geq\frac{\Cr{0}}{96}m\log\log\ltwo{\bx}
    \,\Rp
    \leq e^{-m\C\log\ltwo{\bx}}.
\]
Combining this with \eqref{eq:p4.46} we get \eqref{eq:p4.42}.
Therefore, to complete the proof of 
Proposition~\ref{prop:finecoarsediffsum} we only need to
prove Claim~\ref{claim:finecoarsediffsum}.

\begin{proof}[Proof of Claim~\ref{claim:finecoarsediffsum}] 
    Fix $0\leq i\leq m-1$ and $(\bv,\bv^\prime)\in V_i$. 
    Define
    \[
        \diff(\bv,\bv^\prime):=
            h(\bv-\bv^\prime)-
            \hat{T}^i(\bv,\bv^\prime)
            -h(\bw_{2i}-\bw_{2i+1})
            +\hat{T}^i(\bw_{2i},\bw_{2i+1})).
    \]
    Define $\ell_0:=\ltwo{\bx}(\log\ltwo{\bx})^{-1/\alpha}$.
    We consider two cases, either 
    $\xprojth(\bv^\prime-\bv)\leq\ell_0$, or 
    $\xprojth(\bv^\prime-\bv)\geq\ell_0$.
    
    \vspace{\baselineskip}
    
    \textit{Case I:} Suppose 
    $\xprojth(\bv^\prime-\bv)\leq\ell_0$. 
    Using Lemma~\ref{lem:width}, $\yprojth(\bv-\bw_{2i})$ and
    $\yprojth(\bv^\prime-\bw_{2i+1})$ are at most
    $\C\Delta(\bx)\log\log\ltwo{x}$ which is smaller than
    $\ell_0$. Hence $\ltwo{\bv-\bw_{2i}}$ and 
    $\ltwo{\bv^\prime-\bw_{2i+1}}$ are at most of the order 
    of $\ell_0$. Therefore, using \eqref{A1} and \eqref{A2}, we get 
    \[
        \Prob\Lp\,
            \diff(\bv,\bv^\prime)\geq 
            t\sigma(\ltwo{\bx})
            \,\Rp
        \leq
        \C e^{-\C t\sigma(\ltwo{\bx})/\sigma(\ltwo{\ell_0})}
        \leq e^{-\C t\log\ltwo{\bx}}.
    \]
    So the claim is proved in this case.  
    
    \vspace{\baselineskip}
    
    \textit{Case II:} Suppose 
    \begin{equation}\label{case2ass}
        \xprojth(\bv^\prime-\bv)\geq\ell_0.    
    \end{equation} 
    Let
    \begin{align*}
    & \diff_1(\bv,\bv^\prime):=         
        h(\bv-\bv^\prime)-T(\bv,\bv^\prime)-
            h(\bv-\bw_{2i+1})+T(\bv,\bw_{2i+1}),\\
    & \diff_2(\bv,\bv^\prime):=             
        h(\bv-\bw_{2i+1})-T(\bv,\bw_{2i+1})-
            h(\bw_{2i}-\bw_{2i+1})+T(\bw_{2i},\bw_{2i+1}).
    \end{align*}
    Therefore, 
    $\diff(\bv,\bv^\prime)=
        \diff_1(\bv,\bv^\prime)+\diff_2(\bv,\bv^\prime)$.
    Hence
    \begin{equation}\label{eq:p542}
        \Prob\Lp\,\diff(\bv,\bv^\prime)\geq 
            t\sigma(\ltwo{\bx})\,\Rp\leq 
        \Prob\Lp\,\diff_1(\bv,\bv^\prime)\geq             
            \frac{t}{2}\sigma(\ltwo{\bx})\,\Rp+ 
        \Prob\Lp\,\diff_2(\bv,\bv^\prime)\geq 
            \frac{t}{2}\sigma(\ltwo{\bx})\,\Rp.
    \end{equation}
    We only consider the first term on the right-hand side, 
    the second term can be dealt with similarly. 
    Suppose $i_0$ and $j_0$ are such that $\bv^\prime\in 
    B_{i_0 j_0}$, where recall that $B_{ij}$ is defined in \eqref{eq:coarsedef2}. 
    Let  
    \[
    \begin{aligned}
        J&:=\LT j_0\Delta(\ell_{\bx})
            -t^{1/2}\Delta(\ell_{\bx})(\log\ell_{\bx})^{1/2},
            (j_0+1)\Delta(\ell_{\bx})+t^{1/2}
            \Delta(\ell_{\bx})(\log\ell_{\bx})^{1/2}
        \RT,\\
        R(t)&:= \LP\, \by\in\RR^2\,:\,
        \xprojth(\by)=i_0\ell_{\bx},\,
        \yprojth(\by)\in J\,\RP.
    \end{aligned}
    \]   
    So $R(t)$ is an extension of a side of the parallelogram 
    $B_{i_0j_0}$. Define the event
    \[
        \mbox{$\event(t)$: $\Gamma(\bv,\bv^\prime)$ intersects $R(t)$}.
    \]
    Since $\bv^\prime\in B_{i_0j_0}$, distance of the segment
    $R(t)$ from $\bv^\prime$ in $-\theta$ direction is less
    than $\ell_{\bx}$. Therefore, if 
    $\mathbb{T}\not\in\event(t)$, then 
    $\wandering{\bv^\prime}{\bv}{k}{-\theta}\geq 
    t^{1/2}\Delta(\ell_{\bx})(\log\ell_{\bx})^{1/2}$
    for some $k\leq\ell_{\bx}$. Hence, to bound the 
    probability of $\mathbb{T}\not\in\event(t)$ we use 
    Corollary~\ref{cor:endwandlogspl} with
    \[
        \tilde{\thn}:=-\theta,\quad
        \tilde{n}:=\xprojth(\bv^\prime-\bv),\quad
        \tilde{l}:=\yprojth(\bv-\bv^\prime),\quad
        \tilde{k}:=\ell_{\bx},\quad
        \tilde{t}:=t^{1/2}.
    \]
    We now verify the conditions of 
    Corollary~\ref{cor:endwandlogspl}. Recall that by our
    choice of $\delta$ from \eqref{eq:deltachoice} $\theta$ 
    is a direction of type I, hence so is $-\theta$. Using
    \eqref{case2ass} we get 
    $\tilde{n}\geq\ell_0\geq\tilde{n_0}$.
    By Lemma~\ref{lem:width} we have 
    $|\tilde{l}|\leq\C\Delta(\ltwo{x})\log\log\ltwo{\bx}$.
    Hence, $|\tilde{l}|\leq\tilde{n}\tilde{\delta_2}$, as 
    required. Using $t\leq 2^{N_1-1}$ and \eqref{eq:p4.43.2},
    we get 
    $t^{1/2}\Delta(\tilde{k})
    (\log\tilde{k})\leq\tilde{k}\tilde{\delta_3}$, as 
    required. Thus, all the conditions for applying 
    Corollary~\ref{cor:endwandlogspl} hold, and we get
    \begin{equation}\label{eq:p543}
        \Prob\Lp\,\event(t)^c\,\Rp
        \leq e^{-\C t\log\ell_{\bx}}
        \leq e^{-\C t\log\ltwo{\bx}}.
    \end{equation}
    Let $R:=R(2^{N_1-1})$. 
    For any $\by$ in $R$ let 
    \begin{align*}
        \diff^\prime_1(\by) & := h(\bv-\by)-T(\bv,\by)-
        h(\bv-\bw_{2i+1})+T(\bv,\bw_{2i+1}),\\
        \diff^\prime_2(\by) & := h(\by-\bv^\prime)-T(\by,\bv^\prime).
    \end{align*}
    If $\Gamma(\bv,\bv^\prime)$ passes through $\by$ then 
    $\diff_1(\bv,\bv^\prime)
        \leq\diff^\prime_1(\by)+\diff^\prime_2(\by)$.
    Hence
    \begin{multline}\label{eq:p544}
        \Prob\Lp\,
            \diff_1(\bv,\bv^\prime)
                \geq\frac{t}{2}\sigma(\ltwo{\bx})
                \,\Rp
        \leq 
        \Prob\Lp\,
        \max_{\by\in R}             
            \diff^\prime_1(\by)
                \geq\frac{t}{4}\sigma(\ltwo{\bx})
                \,\Rp\\ + 
        \Prob\Lp\,
        \max_{\by\in R}             
            \diff^\prime_2(\by)
                \geq\frac{t}{4}\sigma(\ltwo{\bx})
                \,\Rp
        + \Prob\Lp\,\event(t)^c\,\Rp.
    \end{multline} 
    Let us consider the first term in the right-hand side
    first. We use Corollary~\ref{cor:logupinc} with 
    \begin{align*}
        &\tilde{\thn}:=\theta,\quad
        \tilde{n}:=i_0\ell_{\bx}-\xprojth(\bv),\\
        &\tilde{L}:=
            \Delta(\ell_{\bx})
            (1+2^{(N_1-1)/2+1}(\log\ell_{\bx})^{1/2}),\\
        &\tilde{d}:=\Delta(\ell_{\bx})
            (j_0-2^{(N_1-1)/2}(\log\ell_{\bx})^{1/2})
            -\yprojth(\bv).
    \end{align*}
    We now verify the conditions of 
    Corollary~\ref{cor:logupinc}.
    By our choice of $\delta$ in \eqref{eq:deltachoice}, 
    $\theta$ is of both type I and type II.
    Since $\tilde{n}$ is the distance of $R$ from $\bv$ in 
    direction $-\theta$, we have
    \[
        \tilde{n}\geq \ell_0-\ell_{\bx}
                    \geq \Cl{331}\ell_0
                        \geq \tilde{n_0}.
    \]
    Using \eqref{eq:p4.43.2}, \eqref{A2}, and 
    $\tilde{n}\geq\Cr{331}\ell_0$ from above, we get 
    \begin{equation}\label{eq:lengthofR1}
        \tilde{L}
        \leq \Cl{332}\Delta(\ell_{\bx})\log\ltwo{\bx}
        \leq \C\Delta(\ell_0)
            (\log\ltwo{\bx})^{1-(1+\beta)/(2\alpha)}
        \leq \tilde{\delta_2}\Delta(\tilde{n}).
    \end{equation}
    By Lemma~\ref{lem:width} we have 
    \[
        |\yprojth(\bv^\prime-\bv)|
            \leq \C \Delta(\bx)(\log\log\ltwo{\bx})^{1/2}.
    \]
    Since $\bv^\prime\in B_{i_0j_0}$, we have 
    \[
        |\yprojth(\bv^\prime)-j_0 \Delta(\ell_{\bx})| 
        \leq \Delta(\ell_{\bx}).
    \]
    Therefore, using \eqref{A2} and 
    $\tilde{n}\geq\Cr{331}\ell_0$, we get
    \begin{equation}\label{eq:d}
        |\tilde{d}|
            \leq
            \C\Delta(\ltwo{\bx})(\log\log\ltwo{\bx})^{1/2}
            + 2^{N_1/2}\Delta(\ell_{\bx})
                (\log\ell_{\bx})^{1/2} 
            \leq 
        \C\Delta(\ltwo{\bx})(\log\log\ltwo{\bx})^{1/2}
            \leq\tilde{\delta_1}\tilde{n}.
    \end{equation}
    By Corollary~\ref{cor:logupinc} we get for 
    $\tilde{t}\geq\tilde{t_0}$
    \begin{equation}\label{eq:101520@1}
        \Prob\Lp\,
            \max_{\by\in R} \diff^\prime_1(\by)
            \geq\C\tilde{L}\frac{|\tilde{d}|}{\tilde{n}}
            +\tilde{t}\sigma(\Delta^{-1}(\tilde{L}))
            \log\Delta^{-1}(\tilde{L})
        \,\Rp
    \leq\C e^{-\C\tilde{t}\log\tilde{L}}.
    \end{equation}
    Let $\tilde{t}$ be such that 
    \begin{equation}\label{eq:101520@2}
        \C\tilde{L}\frac{|\tilde{d}|}{\tilde{n}}
        +\tilde{t}\sigma(\Delta^{-1}(\tilde{L}))
        \log\Delta^{-1}(\tilde{L})
        =\frac{t}{4}\sigma(\ltwo{\bx}).
    \end{equation}
    We need to verify $\tilde{t}\geq\tilde{t_0}$. Using 
    \eqref{eq:lengthofR1} and \eqref{A2}, we have
    \begin{equation}
    \tilde{L}\leq\Cr{332}\Delta(\ell_{\bx})\log\ltwo{\bx}\leq\C\Delta(\ltwo{\bx})(\log\ltwo{\bx})^{-1/\alpha}.
    \label{eq:lengthofR}
    \end{equation}
    Using this with \eqref{eq:d} we get 
    \begin{equation}\label{eq:101520@3}
        \tilde{L}\frac{|\tilde{d}|}{\tilde{n}}
        \leq \C\sigma(\ltwo{\bx})
        (\log\log\ltwo{\bx})^{1/2}
        (\log\ltwo{\bx})^{-1/\alpha}.
    \end{equation}
    Using $t\geq 2^{N_0}$, 
    \eqref{eq:p4.43.1}, 
    \eqref{eq:101520@2}, 
    \eqref{eq:101520@3}, 
    and \eqref{eq:lengthofR}, 
    we get
    \begin{multline*}
    \tilde{t} 
        \geq \C t\frac{\sigma(\ltwo{\bx})}
            {\sigma(\Delta^{-1}(\tilde{L}))
                \log\Delta^{-1}(\tilde{L})}
        \geq \C t\Lp\frac{\ltwo{\bx}}
            {\Delta^{-1}(\tilde{L})}\Rp^{\alpha}
                (\log\ltwo{\bx})^{-1}\\
        \geq \C t\Lp\frac{\Delta(\ltwo{\bx})}
            {\tilde{L}}\Rp^{2\alpha/(1+\beta)}
                (\log\ltwo{\bx})^{-1}
        \geq \C t 
            (\log\ltwo{\bx})^{(1-\beta)/(1+\beta)}
        \geq \tilde{t_0}.
    \end{multline*}
    Therefore, from \eqref{eq:101520@1} and
    \eqref{eq:101520@2}, we get 
    \begin{equation}\label{eq:p546}
        \Prob\Lp\,
        \max_{\by\in R}
        \diff^\prime_1(\by)\geq\frac{t}{4}\sigma(\ltwo{\bx})
        \,\Rp
        \leq \C e^{-\C t(\log\ltwo{\bx})^{2/(1+\beta)}}
        \leq \C e^{-\C t(\log\ltwo{\bx})}.
    \end{equation}
    Now we consider the second term in the right-hand side of
    \eqref{eq:p544}. By \eqref{eq:lengthofR} width of $R$ is 
    less than $\Delta(\bx)$. Distance of $R$ from 
    $\bv^\prime$ in $\theta$ direction is less than 
    $\ell_{\bx}$. So $\ltwo{\by-\bv^\prime}\leq\C\ell_{\bx}$ 
    for all $\by\in R$. Thus, using \eqref{A1}, \eqref{A2}, 
    and a union bound, we get
    \[
        \Prob\Lp\,\max_{\by\in R}\diff^\prime_2(\by)
        \geq\frac{t}{4}\sigma(\ltwo{\bx})\,\Rp
        \leq 
        \C\ltwo{\bx}
        e^{-\C t\sigma(\ltwo{\bx})/\sigma(\ell_{\bx})}
        \leq 
        \C e^{-\C t(\log\ltwo{\bx})^2}.
    \]
    Using this in \eqref{eq:p544} together with 
    \eqref{eq:p543} and \eqref{eq:p546} we get appropriate 
    bound for the first term in the right-hand side of 
    \eqref{eq:p542}. The second term can be dealt with 
    similarly. This completes the proof of 
    Claim~\ref{claim:finecoarsediffsum}.
\end{proof}

This also completes the proof of 
Proposition~\ref{prop:finecoarsediffsum}.

\subsection{Proof of Theorem~\ref{thm:nrf2}}

\resetconstant

We assume that $\eta<1$ because for $\eta=1$ the result is
same as Proposition~\ref{prop:nrflog}. 
Let $n$ be a positive integer such that 
$(1-\alpha)^n\leq\eta/2$.
Define for $0\leq m\leq n$ and $k\geq 3$
\[
    \psi_m(k):=k^{-\alpha}\sigma(k)(\log k)^{(1-\alpha)^m}(\log\log k)^{1-(1-\alpha)^m}.
\]
Because $\psi_n(k)\leq\phi_\eta(k)$ for large enough $k$, to 
prove Theorem~\ref{thm:nrf2} it is enough to show that $h$ 
satisfies GAP with exponent $\alpha$ and correction factor 
$\psi_n$.

By assumptions of Theorem~\ref{thm:nrf2} $\thn$ is a
direction of both type I and II. Therefore, by
Theorem~\ref{thm:chaploglog}, there exist positive constants
$\delta$, $C_c$, $M_c$, $K$, $a$, such that $h$ satisfies
$\text{CHAP}(\alpha,\widehat{\phi},M_c,C_c,K,a)$ in the
sector of directions $(\thn-\delta,\thn+\delta)$. Define
$\Sector_0:=[0,2\pi]$, and for $1\leq m\leq n$ define 
\[
    \Sector_m:=
    \LT\thn-\delta\frac{n-m+1}{n},
    \thn+\delta\frac{n-m+1}{n}\RT.
\]
We show that $h$ satisfies GAP with exponent $\alpha$ and 
correction factor $\psi_n$ in $\Sector_n$. 

By Proposition~\ref{prop:nrflog}, $h$ satisfies GAP with 
exponent $\alpha$ and correction factor $\phi_1=\psi_0$ in 
all directions. Hence, there exist constants $C_g>0$ and 
$M_g>0$ such that for $\ltwo{\bx}\geq M_g$ we have 
\begin{equation}\label{eq:gapcgmg}
    h(\bx)\leq 
    g(\bx)+C_g\sigma\Lp\ltwo{\bx}\Rp\log\ltwo{\bx}.
\end{equation}
We use an inductive argument. Fix $0\leq m<n$. Suppose $h$
satisfies GAP with exponent $\alpha$, correction factor
$\psi_m$, in the sector $\Sector_{m}$, with constants $C$ and
$M$. We show that $h$ satisfies GAP with exponent $\alpha$, 
correction factor $\psi_{m+1}$, in the sector 
$\Sector_{m+1}$, with constant $C$ and $M$. This establishes
that $h$ satisfies GAP with exponent $\alpha$, correction
factor $\psi_n$ in the sector $\Sector_n$. The constants $C$
and $M$ need to remain unchanged. We will see that if $C$ and
$M$ are chosen large enough then the inductive step works. We
assume without loss of generality $K>1$, $M>3$, $M>M_c$, 
$M>M_g$. Also we assume $M$ is large enough, independent of 
$m$, so that $\psi_{m+1}(\ltwo{\bx})\geq 1$ for all 
$\ltwo{\bx}\geq M$, which is possible by \eqref{A2}. Since 
$h$ has sublinear growth, there exists constant $r>0$ such 
that for all $\bx$ we have $h(\bx)\leq r\ltwo{\bx}$. Let  
$\nu:=(1-\beta)/4$, 
$c_0:=3C_c$, 
$c_1:=C_c a$, 
$c_2:=3K$, 
$c_3:=(c_1+c_0+C_g c_2(\alpha\nu)^{-1})\lowconstinv$, 
$c_4:=c_2^\alpha\lowconstinv$, 
\[
    c_5:=
    c_3^\alpha c_4^{1-\alpha}
    \Lp\Lp\frac{\alpha}{1-\alpha}\Rp^{1-\alpha}
    +\Lp\frac{1-\alpha}{\alpha}\Rp^{\alpha}\Rp,
\]
and $c_6:=(1-\alpha)c_3(\alpha c_4)^{-1}$. We start the 
inductive step now. Consider $\bx$ with direction 
$\theta\in\Sector_{m+1}$ and $\ltwo{\bx}\geq M$. We need to 
show
\begin{equation}\label{eq:toverify}
    h(\bx)\leq 
    g(\bx)+C\ltwo{\bx}^\alpha\psi_{m+1}(\ltwo{\bx}).
\end{equation}
We are free to choose $C$ and $M$ large enough, independent 
of $m$. In various steps we assume $C$ is large enough 
depending on $M$, and $M$ is chosen to be large enough 
without depending on $C$.

\vspace{\baselineskip}

\paragraph{\textbf{Bounding $h(\bx)$ when $\ltwo{\bx}$ is 
small:}} 
Suppose $\ltwo{\bx}\leq c_2 M$. Assuming $C\geq r c_2 M$, we 
get  
\[
    h(\bx)\leq r\ltwo{\bx}\leq r c_2 M\leq C \leq 
    C\ltwo{\bx}^\alpha\psi_{m+1}(\bx) \leq g(\bx) + 
    C\ltwo{\bx}^\alpha\psi^m(\ltwo{\bx}).
\]
Thus \eqref{eq:toverify} is verified.

\vspace{\baselineskip}

\paragraph{\textbf{Defining $\bx^\ast$, $\xl$ and $\xs$ when 
$\ltwo{\bx}$ is large:}}
Suppose  $\ltwo{\bx}\geq c_2 M$. Take
$q\in[c_2,\ltwo{\bx}/M]\cap\QQ$. Then 
$\ltwo{\bx/q}\geq M\geq M_c$. Applying 
$\text{CHAP}(\alpha,\widehat{\phi},M_c,C_c,K,a)$ to $\bx/q$
we get 
\begin{equation}\label{eq:xqdecomp}
    \bx/q=\sum_{i=1}^{3}\alpha_{qi}\by_{qi}
    \mbox{ with }
    \alpha_{qi}\geq 0, \;\sum_{i=1}^{3}\alpha_{qi}\in[1,a]
    \mbox{ and }\by_{qi}\in Q_{\bx/q}(\alpha,\widehat{\phi},C_c,K).
\end{equation}
Let 
\[
L(q):=\LP 1\leq i\leq 3:\ltwo{\by_{qi}}\geq\ltwo{\bx/q}^{1-\nu}\RP
\]
and
\begin{equation}
\bx^\ast:=\sum_{i=1}^{3}\lfloor q\alpha_{qi}\rfloor\by_{qi},\;
\xl:=\sum_{i\in L(q)}\gamma_{qi}\by_{qi},\; 
\xs:=\sum_{i\not\in L(q) }\gamma_{qi}\by_{qi},
\label{eq:xstxlxs}\end{equation}
where
\[
\gamma_{qi}:=q\alpha_{qi}-\lfloor q\alpha_{qi}\rfloor\in[0,1).
\]
Therefore,
\[
\bx=\bx^\ast+\xl+\xs.
\]

\vspace{\baselineskip}

\paragraph{\textbf{Direction of $\xl$:}} 
Consider $i\in L(q)$. Then 
\begin{equation}\label{eq:dirxl1}
    \ltwo{\by_{qi}}\geq\ltwo{\bx/q}^{1-\nu}.
\end{equation}
Using Lemma~\ref{lem:width}, 
$\by_{qi}\in Q_{\bx/q}(\alpha,\widehat{\phi},C_c,K)$, 
$\ltwo{\bx/q}\geq M$, assuming $M$ is large enough, and 
\eqref{A2}, we get 
\[
    \yprojth(\by_{qi})
    \leq\C\Delta(\bx/q)(\log\log\ltwo{\bx/q})^{1/2}
    \leq\C\ltwo{\bx/q}^{(1+\beta)/2}
        (\log\log\ltwo{\bx/q})^{1/2}.
\]
Therefore, using \eqref{eq:dirxl1}, \eqref{A2}, and 
$1-\nu=1-(1-\beta)/4>(1+\beta)/2$, we get 
\[
    \xprojth(\by_{qi})
    \geq\ltwo{\by_{qi}}-\yprojth(\by_{qi})  \geq\C\ltwo{\bx/q}^{1-\nu},
\]
and, further, using $1-\nu-(1+\beta)/2=(1-\beta)/4$ we get 
\[
    \frac{\yprojth(\by_{qi})}{\xprojth(\by_{qi})}
    \leq\C\frac{\ltwo{\bx/q}^{(1+\beta)/2}
        (\log\log\ltwo{\bx/q})^{1/2}}{\ltwo{\bx/q}^{1-\nu}}
    \leq\C\frac{(\log\log M)^{1/2}}{M^{(1-\beta)/4}}.
\]
Since $\theta\in\Sector_{m+1}$, assuming $M$ is large enough,
we get direction of $\by_{qi}$ is in $\Sector_m$. This implies $\xl$ has direction in $\Sector_m$.

\vspace{\baselineskip}

\paragraph{\textbf{Bounding $h(\bx^\ast)$:}} Using
subadditivity of 
$h$, \eqref{eq:xstxlxs}, \eqref{eq:xqdecomp}, and \eqref{A2},
we get
\begin{multline}\label{eq:hxt2}
    h(\bx^\ast)
    \leq\sum_{i=1}^{3}\lfloor q\alpha_{qi}\rfloor h(\by_{qi})
    \leq\sum_{i=1}^{3}\lfloor q\alpha_{qi}\rfloor
    \LT g_{\bx}(\by_{qi})+
    C_c\ltwo{\bx/q}^\alpha\widehat{\phi}(\ltwo{\bx/q})\RT\\
    \leq g_{\bx}(\bx^\ast)+c_1\lowconstinv
    q^{1-\alpha}\ltwo{\bx}^\alpha\widehat{\phi}(\ltwo{\bx}).
\end{multline}

\vspace{\baselineskip}

\paragraph{\textbf{Bounding $h(\xs)$ when $\ltwo{\xs}$ is 
large:}}
Suppose $\ltwo{\xs}\geq M_g$. Using $q\geq c_2=3K>3$, 
$\ltwo{\bx/q}\geq M\geq 1$, and \eqref{eq:xstxlxs} we get
\[
         \ltwo{\xs}
    =    \ltwo{\sum_{i\not\in L(q)}\gamma_{qi}\by_{qi}}
    \leq \sum_{i\not\in L(q)}\gamma_{qi}\ltwo{\by_{qi}}
    \leq 3\ltwo{\bx/q}^{1-\nu}
    \leq 3\ltwo{\bx/q}
    \leq \ltwo{\bx}.
\]
Using this and 
$\log\ltwo{\bx}\leq\ltwo{\bx}^{\alpha\nu}(\alpha\nu)^{-1}$ 
we get 
\begin{multline*}
    \frac{\sigma(\ltwo{\xs})\log\ltwo{\xs}}
    {\sigma(\ltwo{\bx})\log\log\ltwo{\bx}}
    \leq
    \lowconstinv\Lp\frac{\ltwo{\xs}}{\ltwo{\bx}}\Rp^\alpha
    \frac{\log\ltwo{\bx}}{\log\log\ltwo{\bx}}
    \leq\lowconstinv
    \frac{c_2^\alpha q^{\alpha\nu}\log\ltwo{\bx}}
    {q^\alpha\ltwo{\bx}^{\alpha\nu}\log\log\ltwo{\bx}}
    \leq\lowconstinv 
    c_2^\alpha(\alpha\nu)^{-1}q^{-\alpha(1-\nu)}.
\end{multline*}
Since $\ltwo{\xs}\geq M_g$, using \eqref{eq:gapcgmg} we get
\begin{align*}
    h(\xs) 
    & \leq g(\xs)+C_g\ltwo{\xs}^\alpha\phi_1(\ltwo{\xs})\\
    ~
    & \leq \sum_{i\not\in L(q)}\gamma_{qi} g(\by_{qi}) + 
    C_g\sigma(\ltwo{\xs})\log\ltwo{\xs}\\
    ~ 
    & \leq g_{\bx}(\xs)+\sum_{i\not\in L(q)}\gamma_{qi}\LT 
    g(\by_{qi})-g_{\bx}(\by_{qi})\RT+ 
    C_g\sigma(\ltwo{\xs})\log\ltwo{\xs} \\
    ~
    & \leq g_{\bx}(\xs)+\sum_{i\not\in L(q)}\gamma_{qi}\LT 
    g(\by_{qi})-g_{\bx}(\by_{qi})\RT+C_g\lowconstinv 
    c_2^\alpha (\alpha\nu)^{-1} 
    q^{-\alpha(1-\nu)}\sigma(\ltwo{\bx})\log\log\ltwo{\bx}\\
    ~
    & \leq g_{\bx}(\xs)+\sum_{i\not\in L(q)}\gamma_{qi}\LT 
    g(\by_{qi})-g_{\bx}(\by_{qi})\RT+C_g\lowconstinv 
    c_2^\alpha(\alpha\nu)^{-1} q^{1-\alpha}
    \sigma(\ltwo{\bx})\log\log\ltwo{\bx}.
    \numberthis\label{eq:hxslarge}
\end{align*}

\vspace{\baselineskip}

\paragraph{\textbf{Bounding $h(\xs)$ when $\ltwo{\xs}$ is 
small:}} 
Suppose $\ltwo{\xs}\leq M_g$.
Since $\by_{qi}\in Q_{\bx/q}(\alpha,\phi_2,C_c,K)$,
\[
    0\leq h(\by_{qi})\leq g_{\bx}(\by_{qi})+C_c 
    q^{-\alpha}\ltwo{\bx}^\alpha\widehat{\phi}(\ltwo{\bx/q}).
\]
Therefore
\[
    g_{\bx}(\by_{qi})\geq - C_c 
    q^{-\alpha}\ltwo{\bx}^\alpha\widehat{\phi}(\ltwo{\bx/q}).
\]
So, letting $I(q):=\LP i\leq 3:g_{\bx}(\by_{qi})<0\RP$, and using definition of $c_0$, and \eqref{A2}, we have 
\begin{multline*}
    g_{\bx}(\xs) 
       = \sum_{i\not\in L(q)}\gamma_{qi}g_{\bx}(\by_{qi})
    \geq - \sum_{i\in I(q)}|g_{\bx}(\by_{qi})|
    \geq - c_0 q^{-\alpha}
        \ltwo{\bx}^\alpha\widehat{\phi}(\ltwo{\bx/q})
    \geq - \lowconstinv c_0     
    q^{-\alpha}\ltwo{\bx}^\alpha\widehat{\phi}(\ltwo{\bx}).
\end{multline*}
Therefore, using $\ltwo{\xs}\leq M_g$ and 
$h(\xs)\leq r\ltwo{\xs}$ 
\begin{equation}\label{eq:hxssmall}
    h(\xs)
    \leq r M_g
    \leq g_{\bx}(\xs) 
    + \lowconstinv c_0 q^{-\alpha} \ltwo{\bx}^\alpha\widehat{\phi}(\ltwo{\bx}) + r M_g.
\end{equation}

\vspace{\baselineskip}

\paragraph{\textbf{Bounding $h(\xl)$ when $\ltwo{\xl}$ is 
large:}}
Suppose $\ltwo{\xl}\geq M$. Using $q\geq c_2\geq 1$, 
$\ltwo{\bx/q}\geq M\geq 1$, \eqref{eq:xstxlxs}, and 
$\ltwo{\by_{qi}}\leq K\ltwo{\bx/q}$, we get
\[
    \ltwo{\xl}
    = \ltwo{\sum_{i\in L(q)} \gamma_{qi} \by_{qi}}
    \leq \sum_{i\in L(q)} \gamma_{qi} \ltwo{\by_{qi}}
    \leq c_2 \ltwo{\bx/q}
    \leq \ltwo{\bx}.
\]
Using this and applying $\text{GAP}(\alpha,\psi_m,M,C)$, which holds by the induction hypothesis for $m$, we get 
\begin{align*}
    h(\xl) 
    & \leq g(\xl) + C \ltwo{\xl}^\alpha \psi_m 
    \Lp\ltwo{\xl}\Rp\\
    ~
    & \leq \sum_{i\in L(q)}\gamma_{qi} g(\by_{qi}) 
    + C\lowconstinv c_2^\alpha 
    q^{-\alpha}\ltwo{\bx}^\alpha\psi_m\Lp\ltwo{\bx}\Rp\\ 
    ~
    & \leq g_{\bx}(\xl)
    +\sum_{i\in L(q)}\gamma_{qi}\LT 
    g(\by_{qi})-g_{\bx}(\by_{qi})\RT 
    + C\lowconstinv c_2^\alpha 
    q^{-\alpha}\ltwo{\bx}^\alpha\psi_m\Lp\ltwo{\bx}\Rp
    \numberthis\label{eq:hxllarge}
\end{align*}

\vspace{\baselineskip}

\paragraph{\textbf{Bounding $h(\xl)$ when $\ltwo{\xl}$ is small:}} 
Suppose $\ltwo{\xl}\leq M$. Then by similar calculations that
lead to \eqref{eq:hxssmall} we get 
\begin{equation}\label{eq:hxlsmall}
         h(\xl)
    \leq r M 
    \leq g_{\bx}(\xl) + 
    \lowconstinv c_0 q^{-\alpha} \ltwo{\bx}^\alpha 
    \widehat{\phi}(\ltwo{\bx}) + rM.
\end{equation}

\vspace{\baselineskip}

\paragraph{\textbf{Overall bound on $h(\bx)$:}} 
Using $\by_{qi}\in Q_{\bx/q}(\alpha,\widehat{\phi},C_c,K)$,
definition of $c_0$, and \eqref{A2},
\begin{multline*}
    \sum_{i=1}^{3}\gamma_{qi}
        \LT g(\by_{qi})-g_{\bx}(\by_{qi})\RT
    \leq \sum_{i=1}^{3}\gamma_{qi}
        \LT h(\by_{qi})-g_{\bx}(\by_{qi})\RT
    \leq c_0 \ltwo{\bx/q}^{\alpha}
        \widehat{\phi}(\ltwo{\bx/q})
    \leq c_0 \lowconstinv q^{-\alpha}\ltwo{\bx}^\alpha
        \widehat{\phi}(\ltwo{\bx}).
\end{multline*}
Combining this with \eqref{eq:hxt2}-\eqref{eq:hxlsmall}, we 
get
\begin{align*}
    h(\bx) 
    & \leq h(\bx^\ast) + h(\xl) + h(\xs)\\
    & \leq g(\bx) + 
    c_3 q^{1-\alpha}\sigma(\ltwo{\bx})\log\log\ltwo{\bx} 
    + C c_4 q^{-\alpha}
    \sigma(\ltwo{\bx})
    (\log\ltwo{\bx})^{(1-\alpha)^m}
    (\log\log\ltwo{\bx})^{1-(1-\alpha)^m} \\ 
    & \quad + rM + rM_g.
    \numberthis\label{eq:hx2}
\end{align*}


\paragraph{\textbf{Optimizing over $q$:}}
The optimal $q$ that 
minimizes the right-hand side \eqref{eq:hx2} is 
\[
    q_0:=
    C\frac{\alpha c_4}{(1-\alpha)c_3}
    \Lp\frac{\log\ltwo{\bx}}
        {\log\log\ltwo{\bx}}\Rp^{(1-\alpha)^m}.
\]
Plugging in $q=q_0$ in \eqref{eq:hx2} we see that if $C$ is 
large enough depending on $r$, $M$, $M_g$, then
\begin{align*}
    h(\bx) & \leq g(\bx) + 
    c_5 C^{1-\alpha}         
        \sigma(\ltwo{\bx})
        (\log\ltwo{\bx})^{(1-\alpha)^{m+1}}
        (\log\log\ltwo{\bx})^{1-(1-\alpha)^{m+1}}
        + rM + rM_g\\
    & \leq g(\bx) + C\sigma(\ltwo{\bx})
    (\log\ltwo{\bx})^{(1-\alpha)^{m+1}} (\log\log\ltwo{\bx})^{1-(1-\alpha)^{m+1}}.
\end{align*}
Thus we get \eqref{eq:toverify} provided we prove $q_0$ is 
feasible.

\vspace{\baselineskip}

\paragraph{\textbf{Feasibility of $q$:}} 
We need to verify that 
$q_0\in [c_2,\ltwo{\bx}/M]$. We get $q_0\geq c_2$ using 
$\ltwo{\bx}\geq M$, $m\leq n$, choosing $C>1$, and assuming 
$M$ is large enough. Suppose $\ltwo{\bx}<q_0 M$ so that we 
have 
\begin{equation}\label{eq:aux1}
    CM\geq c_6\ltwo{\bx}
    \Lp\frac{\log\ltwo{\bx}}
        {\log\log\ltwo{\bx}}\Rp^{-(1-\alpha)^m}.
\end{equation}
This gives an upper bound on $\ltwo{\bx}$. 
So $q_0$ is not feasible when $\ltwo{\bx}$ is too small. 
But we can prove \eqref{eq:toverify} in a different way. 
Consider two cases. \textit{Case I:} Suppose $\bx$ is such 
that
\[
    C\geq 
    C_g\Lp\frac{\log\ltwo{\bx}}
    {\log\log\ltwo{\bx}}\Rp^{1-(1-\alpha)^{m+1}}.
\]
Then from \eqref{eq:gapcgmg} we get
\[
    h(\bx)\leq 
    g(\bx)+
    C\sigma(\ltwo{\bx})
    (\log\ltwo{\bx})^{(1-\alpha)^{m+1}}
    (\log\log\ltwo{\bx})^{1-(1-\alpha)^{m+1}}.
\]
Thus \eqref{eq:toverify} is verified. \textit{Case II:} 
Suppose $\bx$ is such that
\[
    C\leq C_g\Lp\frac{\log\ltwo{\bx}}
        {\log\log\ltwo{\bx}}\Rp^{1-(1-\alpha)^{m+1}}.
\]
Combining with \eqref{eq:aux1} we get
\[
    C_g M \Lp\frac{\log\ltwo{\bx}}
    {\log\log\ltwo{\bx}}\Rp^{1-(1-\alpha)^{m+1}} 
    \geq c_6\ltwo{\bx}\Lp\frac{\log\ltwo{\bx}}
    {\log\log\ltwo{\bx}}\Rp^{-(1-\alpha)^m}.
\]
So
\[
    \frac{C_g M}{c_6} 
    \geq\ltwo{\bx}
    \Lp\frac{\log\ltwo{\bx}}
    {\log\log\ltwo{\bx}}\Rp^{-(1-(1-\alpha)^{m+1}
                                    +(1-\alpha)^m)}
    \geq\ltwo{\bx}\Lp\frac{\log\ltwo{\bx}}
        {\log\log\ltwo{\bx}}\Rp^{-2}.
\]
Therefore, $M\geq F(\ltwo{\bx})$, where 
$F(k):= c_6 C_g^{-1} k(\log k/\log\log k)^{-2}$. Taking 
$C\geq r F^{-1}(M)$, where 
$F^{-1}(k^\prime):=\sup\LP k: F(k)\leq k^\prime\RP$, we get 
\[
    h(\bx)
    \leq r\ltwo{\bx}
    \leq r F^{-1}(M)
    \leq C 
    \leq g(\bx)+C\psi_{m+1}(\ltwo{\bx}).
\] 
Thus \eqref{eq:toverify} is verified. This completes the 
inductive step and proves Theorem~\ref{thm:nrf2}.

\section{Upper bound of the transverse increments}
\label{sec:loglogupinc}

\resetconstant

\begin{figure}[H]
    \centering
    \includegraphics[width=0.65\linewidth]{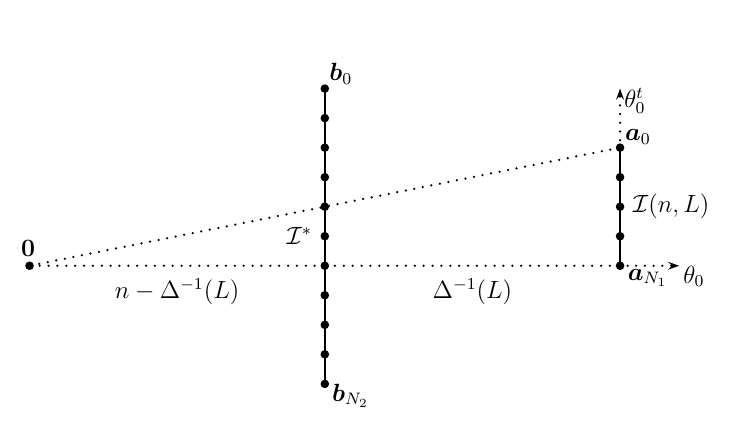}
    \caption{Setup of Theorem~\ref{thm:loglogupinc}: width of
    the portion of $\segment^\ast$ between the dotted lines
    is $L(1-\Delta^{-1}(L)/n)$; width of the portion of
    $\segment^\ast$ above the line joining $\origin$ and
    $\ba_0$ is $t^{1/2}L\Lp\log\Delta^{-1}(L)\Rp^{\eta/2}$;
    width of the portion of $\segment^\ast$ below the line
    joining $\origin$ and $\ba_{N_1}$ is also
    $t^{1/2}L\Lp\log\Delta^{-1}(L)\Rp^{\eta/2}$.}
    \label{Fig:loglogupinc}
\end{figure}

In this section, we prove Theorem~\ref{thm:loglogupinc}, 
which is our main result on upper bound of the transverse 
increments. Fix $L_0>0$, $n_0>0$, $t_0>0$, to be assumed 
large enough whenever required. Consider $n$, $L$, $t$, 
satisfying $n\geq n_0$, $L\geq L_0$, $t\geq t_0$, and 
$L\leq\Delta(n)$. If 
$t\geq 4\mu L(\sigma(\Delta^{-1}(L))
(\log\Delta^{-1}(L))^\eta)^{-1}$, 
where $\mu$ is the expected edge-weight, then we are in a 
large-deviation regime, and the proof is similar to Case I of
Theorem~\ref{thm:logupinc}. Therefore, let us assume 
\begin{equation}\label{eq:loglogupinc0}
    t\leq 
    4\mu L(\sigma(\Delta^{-1}(L))
        (\log\Delta^{-1}(L))^\eta)^{-1}.
\end{equation}
Consider
\[
\begin{aligned}
    &J:=\LT-t^{1/2}L\Lp\log\Delta^{-1}(L)\Rp^{\eta/2},
    \Lp 1-\frac{\Delta^{-1}(L)}{n}\Rp L
    + t^{1/2}L\Lp\log\Delta^{-1}(L)\Rp^{\eta/2}\RT,\\
    &\segment^\ast:=\LP\,\bx\in\RR^2:
    \xprojthnot(\bx)=n-\Delta^{-1}(L),\,
    \yprojthnot(\bx)\in J\RP.
\end{aligned}
\]
Let 
\begin{equation}\label{eq:Mvalue}
M :=(1+\beta)/(2\alpha),\quad
N_1 :=\lfloor(\log L)^M\rfloor,\quad
N_2 :=\lfloor t^{1/2}(\log L)^{M+\eta/2}\rfloor.
\end{equation}
Divide the segment $\segment(n,L)$ in $N_1$ segments of equal
length: $\segment_1,\dots,\segment_{N_1}$, with endpoints
$\ba_0,\dots,\ba_{N_1}$, as shown in 
Figure~\ref{Fig:loglogupinc}. Divide the segment
$\segment^\ast$ in $N_2$ segments of equal length:
$\segment^\ast_1,\dots,\segment^\ast_{N_2}$, with endpoints 
$\bb_0,\dots,\bb_{N_2}$, as shown in 
Figure~\ref{Fig:loglogupinc}. By \eqref{A2}, $\log L$ is of 
the same order as $\log\Delta^{-1}(L)$, i.e., 
\begin{equation}\label{eq:loglogupinc0.5}
    \C\log L\leq\log\Delta^{-1}(L)\leq\C\log L.
\end{equation}
Therefore, length of the segments $\segment_i$ and 
$\segment^\ast_j$ are bounded as 
\begin{gather}
    \C L(\log L)^{-M} \leq
    \yprojthnot(\ba_{i-1}-\ba_i) \leq 
    \C L(\log L)^{-M},\label{eq:loglogupinc01}\\ 
    \C L(\log L)^{-M} \leq
    \yprojthnot(\bb_{j-1}-\bb_j) \leq 
    \C L(\log L)^{-M}.\label{eq:loglogupinc02} 
\end{gather}

Define the event 
\[
    \event:\mbox{$\Gamma(\origin,\ba_i)$ passes through 
    $\segment^\ast$ for all $0\leq i\leq N_1$}. 
\]
Define
\begin{align*}
    D^1 &:=
    \max\LP\,|T(\origin,\ba_{i_1})-T(\origin,\ba_{i_2})|\,:\,
    0\leq i_1 < i_2 \leq N_1\,\RP,\\
    D^2 &:=
    \max\LP\,|T(\bx,\ba_{i_1})-T(\bx,\ba_{i_2})|\,:\,
    0\leq i_1 < i_2 \leq N_1,\,
    \bx\in\segment^\ast\,\RP,\\
    D^3 &:=
    \max\LP\,|T(\bb_j,\ba_{i_1})-T(\bb_j,\ba_{i_2})|\,:\,
    0\leq i_1 < i_2 \leq N_1,\,
    0\leq j\leq N_2\,\RP,\\
    D^4 &:=\max\LP\,|T(\bb,\ba_i)-T(\bb^\prime,\ba_i)|
    \,:\,\bb,\bb^\prime\in\segment^\ast_j,\,
    0\leq j\leq N_2,\,
    0\leq i\leq N_1\,\RP,\\
    D_i &:=\max\LP\,|T(\origin,\bx)-T(\origin,\by)|\,:\,
    \bx,\by\in\segment_i\,\RP,\mbox{ for } 0\leq i\leq N_1.
\end{align*}
Therefore
\begin{align*}
    \diff(n,L) & \leq D^1+\max_{0\leq i\leq N_1} D_i,\\
    D^2 & \leq D^3+2 D^4.
    \numberthis\label{eq:loglogupinc3}
\end{align*}
We claim that if $\mathbb{T}\in\event$ then  
\begin{equation}\label{eq:loglogupinc4}
    D^1\leq D^2.
\end{equation}
To see this, take $0\leq i_1 < i_2 \leq N_1$. If 
$\mathbb{T}\in\event$, then there exist points $\by$ and 
$\bz$ in $\segment^\ast$ such that 
$\Gamma(\origin,\ba_{i_1})$ passes through $\by$, and 
$\Gamma(\origin,\ba_{i_2})$ passes through $\bz$. Therefore 
\[
    T(\origin,\ba_{i_1})-T(\origin,\ba_{i_2}) 
    \leq T(\origin,\bz)+T(\bz,\ba_{i_1})-
    T(\origin,\bz)-T(\bz,\ba_{i_2}) 
    = T(\bz,\ba_{i_1})-T(\bz,\ba_{i_2}).
\]
Similarly,
\[
    T(\origin,\ba_{i_2})-T(\origin,\ba_{i_1})
    \leq T(\by,\ba_{i_2})-T(\by,\ba_{i_1}).
\]
Therefore
\[
    |T(\origin,\ba_{i_1})-T(\origin,\ba_{i_2})|
    \leq
    \max_{\bx\in\segment^\ast}
    |T(\bx,\ba_{i_1})-T(\bx,\ba_{i_2})|.
\]
Thus, \eqref{eq:loglogupinc4} is proved by taking maximum 
over values of $i_1$, $i_2$.  Combining 
\eqref{eq:loglogupinc3} and \eqref{eq:loglogupinc4} we get 
that on the event $\event$,
\[
    \diff(n,L)\leq D^3+2 D^4+\max_{0\leq i\leq N_1} D_i.
\]
Therefore
\begin{align*}
    \Prob\Lp\,
    \diff(n,L)\geq t\sigma\Lp\Delta^{-1}(L)\Rp
    (\log L)^\eta\,\Rp
    \leq & 
    \Prob\Lp\,\event^c\,\Rp + 
    \Prob\Lp\,D^3\geq\frac{t}{4}
    \sigma\Lp\Delta^{-1}(L)\Rp(\log L)^\eta\,\Rp\\ 
    & + 
    \Prob\Lp\,D^4\geq\frac{t}{4}
        \sigma\Lp\Delta^{-1}(L)\Rp(\log L)^\eta\,\Rp\\
    & +
    \Prob\Lp\,\max_{0\leq i\leq N_1}     
    D_i\geq\frac{t}{4}
    \sigma\Lp\Delta^{-1}(L)\Rp(\log L)^\eta\,\Rp.
    \numberthis\label{eq:loglogdecomp}
\end{align*}
First, we show that $\Prob(\event^c)$ is small. If
$\mathbb{T}\not\in\event$, then for some $i$,
$\Gamma(\origin,\ba_i)$ wanders more than
$t^{1/2}L(\log\Delta^{-1}(L))^{\eta/2}$ in $\pm\thnt$
directions when it is at a distance $\Delta^{-1}(L)$ from
$\ba_i$ in $-\thn$ direction. Since $\thn$ is a direction of
both type I and type II, by Theorem~\ref{thm:nrf2} $h$
satisfies GAP with exponent $\alpha$ and correction factor
$\phi_\eta$ in a neighborhood of $\thn$. Thus, applying
Theorem~\ref{thm:endwandlog} with the variables 
\begin{align*}
    &\tilde{\thn} = -\thn,\quad
    \tilde{\eta} = \eta,\quad
    \tilde{n} = n,\\
    &\tilde{k}=\Delta^{-1}(L),\quad
    \tilde{l}=\yprojthnot(\ba_{i}),\quad 
    \tilde{t}=t^{1/2},
\end{align*}
and using \eqref{eq:loglogupinc0.5}, we get 
\begin{equation}\label{eq:1016201}
    \Prob\Lp\,
    \wandering{\ba_{i}}{\origin}{\Delta^{-1}(L)}{-\thn}
    \geq t^{1/2}L(\log\Delta^{-1}(L))^{\eta/2}\,\Rp
    \leq \C e^{-\C t(\log L)^\eta}, 
\end{equation}
provided $\tilde{n}\geq\tilde{n_0}$, 
$\tilde{t}\geq\tilde{t_0}$, $\tilde{k}\geq\tilde{k_0}$,  
$\tilde{t}\Delta(\tilde{k})(\log\tilde{k})^{\tilde{\eta}/2}
\leq\tilde{k}\tilde{\delta_1}$, and 
$\tilde{l}\leq\tilde{n}\tilde{\delta_2}$. We verify these 
conditions now. Taking $n_0$, $L_0$, $t_0$ large enough we 
get $\tilde{n}\geq\tilde{n_0}$, $\tilde{t}\geq\tilde{t_0}$, 
and $\tilde{k}\geq\tilde{k_0}$. Using \eqref{eq:loglogupinc0}
and \eqref{A2} we get
\begin{multline}\label{eq:loglogupinc2}
    \frac{\tilde{t}\Delta(\tilde{k})(\log\tilde{k})^{\eta/2}}
    {\tilde{k}}
    = \frac{t^{1/2}
    L(\log\Delta^{-1}(L))^{\eta/2}}{\Delta^{-1}(L)}
    \leq \C\frac{L^{1/2}}{(\Delta^{-1}(L))^{1/2}}\\
    \leq \C L^{-(1-\beta)/2}
    \leq \C L_0^{-(1-\beta)/2}
    \leq \tilde{\delta_1}.
\end{multline}
Using $L\leq\Delta(n)$ and \eqref{A2}, we get
\[
    \frac{|\tilde{l}|}{\tilde{n}}
    \leq\frac{L}{n}
    \leq\frac{\Delta(n)}{n}
    \leq\C n_0^{-(1-\beta)/2}\leq\tilde{\delta_2}.
\]
Thus all the conditions for \eqref{eq:1016201} to hold are 
true. From \eqref{eq:1016201} taking a union bound over $i$
values we get 
\begin{equation}\label{eq:eventloglogfinal}
    \Prob\Lp\,\event^c\,\Rp
    \leq \C e^{-\C t(\log L)^\eta}.
\end{equation}
Now we show that the second term in the right-hand side of
\eqref{eq:loglogdecomp} is small. Take $\bx$ in
$\segment^\ast$ and $\bu$, $\bv$ in $\segment(n,L)$. Then
\begin{multline*}
         |T(\bx,\bu)-T(\bx,\bv)| 
    \leq |T(\bx,\bu)-h(\bu-\bx)| 
       + |T(\bx,\bv)-h(\bv-\bx)|\\
       + |h(\bu-\bx)-g(\bu-\bx)| 
       + |h(\bv-\bx)-g(\bv-\bx)| 
       + |g(\bu-\bx)-g(\bv-\bx)|.
\end{multline*}
From the definition of $\segment^\ast$ it follows 
\begin{equation}\label{eq:loglogupinc7}
    \xprojthnot(\bu-\bx)
    =\Delta^{-1}(L),\quad 
    |\yprojthnot(\bu-\bx)|\leq \C t^{1/2} L (\log\Delta^{-1}(L))^{\eta/2}.
\end{equation}
By same calculation as in \eqref{eq:loglogupinc2} we get that
the direction of $\bu-\bx$ can be made arbitrarily close to
$\thn$ by choosing $L_0$ large enough. So by 
Theorem~\ref{thm:nrf2} and equation~\eqref{eq:loglogupinc0.5}
we get  
\begin{equation}\label{eq:loglogupinc7.1}
    |h(\bu-\bx)-g(\bu-\bx)|
    \leq\C\sigma(\Delta^{-1}(L))(\log L)^{\eta}.
\end{equation} 
Same holds true for $\bv$ replaced by $\bu$. Using \eqref{eq:loglogupinc7}, Lemma~\ref{lem:auxgeom2}, and \eqref{eq:loglogupinc0.5}, we get 
\[
         |g(\bu-\bx)-g(\bv-\bx)| 
    \leq \C t^{1/2}
    \frac{L^2(\log(\Delta^{-1}(L)))^\eta}{\Delta^{-1}(L)}
    \leq \C t^{1/2}\sigma(\Delta^{-1}(L))(\log L)^\eta.   
\]
Using this with \eqref{eq:loglogupinc7.1} and \eqref{A1} 
implies 
\[
     \Prob\Lp\,|T(\bx,\bu)-T(\bx,\bv)|\geq\frac{t}{4}\sigma(\Delta^{-1}(L))(\log L)^\eta\,\Rp
\leq \C e^{-\C t(\log L)^\eta}.
\]
This is true for fixed $\bx$, $\bu$ and $\bv$. So, for all
$0\leq i_1 < i_2\leq N_1$ and $0\leq j\leq N_2$ we get 
\[
    \Prob\Lp\,|T(\bb_{j},\ba_{i_1})-T(\bb_{j},\ba_{i_2})|
    \geq\frac{t}{4}\sigma(\Delta^{-1}(L))(\log L)^\eta\,\Rp
    \leq \C e^{-\C t(\log L)^\eta}.
\]
By \eqref{eq:Mvalue}, the number of triplets $(i_1,i_2,j)$ is
less than $\C t^{1/2}(\log L)^{3M+\eta/2}$. Therefore, by a
union bound we get
\begin{equation}\label{eq:D3final}
    \Prob\Lp\,
    D^3\geq\frac{t}{4}\sigma(\Delta^{-1}(L))
    (\log L)^\eta\,\Rp
    \leq\C e^{-\C t(\log L)^{\eta}}.
\end{equation}
Now let us consider the third term in the right-hand side of
\eqref{eq:loglogdecomp}. Fix $0\leq i\leq N_1$ and 
$1\leq j\leq N_2$. Applying Corollary~\ref{cor:logupinc} with
the variables 
\[
    \tilde{\thn}:=-\thn,\,
    \tilde{n}:=\Delta^{-1}(L),\,
    \tilde{L}:=|\yprojthnot(\bb_{j-1}-\bb_{j})|,\,
    \tilde{d}:=\yprojthnot(\bb_{j}-\ba_{i}),
\]
we get for all $\tilde{t}\geq\tilde{t_0}$
\begin{equation}\label{eq:D4setup}
    \Prob\Lp\,\max_{\bb,\bb^\prime\in\segment^\ast_j}
    |T(\bb,\ba_{i})-T(\bb^\prime,\ba_{i})|
    \geq\Cl{D4}\tilde{L}\frac{|\tilde{d}|}{\tilde{n}} + 
    \tilde{t}\sigma(\Delta^{-1}(\tilde{L}))
    \log\Delta^{-1}(\tilde{L})\,\Rp
    \leq\C e^{-\C t\log\tilde{L}},
\end{equation}
provided the following conditions are satisfied: 
$|\tilde{d}|\leq\tilde{\delta_1}\tilde{n}$, 
$\tilde{L}\leq\tilde{\delta_2}\Delta(\tilde{n})$, 
$\tilde{n}\geq\tilde{n_0}$, $\tilde{L}\geq\tilde{L_0}$. Let 
us now verify these conditions. From definition of 
$\segment^\ast$ we get 
\[
    |\tilde{d}|
    \leq\C t^{1/2}L(\log\Delta^{-1}(L))^{\eta/2},
\]
so that by calculations similar to \eqref{eq:loglogupinc2} 
we get $|\tilde{d}|\leq\tilde{\delta_1}\tilde{n}$. Combining 
this with \eqref{eq:loglogupinc02} we get 
\begin{equation}\label{eq:loglogupinc11}
     \tilde{L}\frac{|\tilde{d}|}{\tilde{n}} 
\leq \C t^{1/2}\sigma(\Delta^{-1}(L))(\log L)^{-M}
\end{equation}
Let $\tilde{t}$ be such that
\begin{equation}\label{eq:D4tsetup}
    \Cr{D4}\tilde{L}\frac{|\tilde{d}|}{\tilde{n}} + \tilde{t}\sigma(\Delta^{-1}(\tilde{L}))
    \log\Delta^{-1}(\tilde{L})
    =\frac{t}{4}\sigma(\Delta^{-1}(L))(\log L)^{\eta}.
\end{equation}
Using \eqref{eq:loglogupinc11}, \eqref{A2}, lower
bound on $\tilde{L}$ from \eqref{eq:loglogupinc02}, value of
$M$ from \eqref{eq:Mvalue}, and \eqref{eq:loglogupinc0.5}, we
get
\begin{multline*}
    \tilde{t}
    \geq
    \C t \frac{\sigma(\Delta^{-1}(L))(\log L)^\eta}
    {\sigma(\Delta^{-1}(\tilde{L}))
        \log\Delta^{-1}(\tilde{L})}
    \\
    \geq \C t (\log L)^{2\alpha M/(1+\beta)+\eta-1}
    \geq \Cl{D4t} t (\log L)^{\eta}
    \geq \Cr{D4t} t_0 (\log L_0)^{\eta}.
\end{multline*}
So we have $\tilde{t}\geq\tilde{t_0}$, assuming $t_0$ and 
$L_0$ are large enough. Therefore, all the conditions for
\eqref{eq:D4setup} are satisfied. Combining 
\eqref{eq:D4setup} with \eqref{eq:D4tsetup}, we get
\[
    \Prob\Lp\,
    \max_{\bb,\bb^\prime\in\segment^\ast_j}
    |T(\bb,\ba_i)-T(\bb^\prime,\ba_i)|
    \geq\frac{t}{4}\sigma(\Delta^{-1}(L))(\log L)^{\eta}\,\Rp
    \leq \C e^{-\C t(\log L)^{1+\eta}},
\]
The number of choices of $i$ and $j$ is at most 
$\C t^{1/2} (\log L)^{2M+\eta/2}$. Hence using a union bound 
we get 
\begin{equation}\label{eq:D4final}
    \Prob\Lp\,D^4\geq\frac{t}{4}
        \sigma(\Delta^{-1}(L))(\log L)^{\eta}\,\Rp
    \leq \C e^{-\C t(\log L)^{1+\eta}}.
\end{equation}
Now we are going to consider the fourth term in the 
right-hand side of \eqref{eq:loglogdecomp}. Fix an $i$. To
bound $D_i$ we apply Corollary~\ref{cor:logupinc} with the 
following variables:
\[
    \tilde{n}:=n,\,
    \tilde{L}:=|\yprojthnot(\ba_{i-1}-\ba_i)|,\,
    \tilde{d}:=\yprojthnot(\ba_i).
\]
By Corollary~\ref{cor:logupinc} we have for all
$\tilde{t}\geq\tilde{t_0}$
\begin{equation}\label{eq:Disetup}
    \Prob\Lp\, 
    D_i\geq\Cl{Di}\tilde{L}\frac{|\tilde{d}|}{\tilde{n}} + 
    \tilde{t}\sigma(\tilde{k})\log \tilde{k}
    \,\Rp
    \leq\C e^{-\C t\log\tilde{k}},
\end{equation}
provided the following conditions are satisfied:
$|\tilde{d}|\leq\tilde{\delta_1}\tilde{n}$, 
$\tilde{L}\leq\tilde{\delta_2}\Delta(\tilde{n})$, 
$\tilde{n}\geq\tilde{n_0}$, $\tilde{L}\geq\tilde{L_0}$.
From definition of $\segment(n,L)$ we have 
$|\tilde{d}|\leq L$. 
Therefore, 
\[
        \frac{|\tilde{d}|}{\tilde{n}}
    \leq\frac{L}{n}\leq\frac{\Delta(n)}{n}
    \leq\C n_0^{-(1-\beta)/2}
    \leq\tilde{\delta_1}.
\]
From \eqref{eq:loglogupinc01} we get $\tilde{L}\leq L$. 
Further using $L\leq\Delta(n)$ we get 
\[
    \tilde{L}\leq 
    L\leq\Delta(n)\leq\tilde{\delta_2}\tilde{n}.
\]
Also, using the bound on $\tilde{L}$ from 
\eqref{eq:loglogupinc01} we get 
\begin{equation}\label{eq:loglogupinc1015201}
    \tilde{L}\frac{|\tilde{d}|}{\tilde{n}}\leq
    \C L(\log L)^{-M}\frac{L}{n}\leq
    \C \frac{L^2}{\Delta^{-1}(L)}(\log L)^{-M}\leq
    \C \sigma(\Delta^{-1}(L))(\log L)^{-M}.
\end{equation}
Let $\tilde{t}$ be such that
\begin{equation}\label{eq:Ditsetup}
    \Cr{Di}\tilde{L}\frac{|\tilde{d}|}{\tilde{n}} + 
    \tilde{t}\sigma(\Delta^{-1}(\tilde{L}))
    \log\Delta^{-1}(\tilde{L}) =
    \frac{t}{4}\sigma(\Delta^{-1}(L))(\log L)^{\eta}.
\end{equation}
Therefore, using \eqref{eq:loglogupinc1015201}, bound on 
$\tilde{L}$ from \eqref{eq:loglogupinc01}, and value of $M$ 
from \eqref{eq:Mvalue}, we get  
\begin{multline*}
    \tilde{t}\geq 
    \C t \frac{\sigma(\Delta^{-1}(L))
    (\log L)^\eta}{\sigma(\Delta^{-1}(\tilde{L})\log 
    \Delta^{-1}(\tilde{L})}\\\geq
    \C t (\log L)^{2\alpha M/(1+\beta)+\eta-1}\geq
    \Cl{loglog50} t (\log L)^{\eta}\geq
    \Cr{loglog50} t_0 (\log L_0)^{\eta}.
\end{multline*}
So $\tilde{t}\geq\tilde{t_0}$ assuming $t_0$ and $L_0$ are
large enough. Combining this with \eqref{eq:Disetup} and 
\eqref{eq:Ditsetup} we get 
\[
    \Prob\Lp\,D_i\geq\frac{t}{4}\sigma(\Delta^{-1}(L))(\log L)^{\eta}\,\Rp
    \leq \C e^{-\C t(\log L)^{1+\eta}}
\]
Using a union bound over values of $i$ we get 
\begin{equation}\label{eq:DiFinal}
    \Prob\Lp\,\max_i D_i\geq 
    t\sigma(\Delta^{-1}(L))(\log L)^{\eta}\,\Rp
    \leq \C e^{-\C t(\log L)^{1+\eta}}.
\end{equation}
Combining \eqref{eq:eventloglogfinal}, \eqref{eq:D3final}, 
\eqref{eq:D4final}, \eqref{eq:DiFinal}, and 
\eqref{eq:loglogdecomp}, completes the proof.

\section{Lower bound of variance of the fluctuations of the 
transverse increments}

\resetconstant
\resetconstant[epsilon]
\resetconstant[nu]

In this section, we prove Theorem~\ref{thm:lowmain}, which is
our main result on the lower bound of the fluctuations of the
transverse increments. In accordance with the statement of
Theorem~\ref{thm:lowmain}, we consider a direction $\thn$
which is of both type I and type II. Without loss of
generality we assume $\thn\in[0,\pi/4]$. We also fix a
constant $\nu\in(1/2,1)$ as in the statement of
Theorem~\ref{thm:lowmain}. In addition, we fix a constant
$\eta\in(1/2,1)$. Consider $n>0$ and $L>0$ satisfying
$L\leq\Delta(n)$. Define $k$ such that 
\begin{equation}\label{eq:lowoutline0}
    L=\Delta(k)(\log k)^{\eta}.
\end{equation}
Therefore, any lower limit of $k$ yields lower limits of both
$n$ and $L$. Hence, we state results which hold for large
enough $k$, tacitly assuming $n$ and $L$ are large enough so
that the two relations $L\leq\Delta(n)$ and
\eqref{eq:lowoutline0} hold. Let us introduce some more
notations that we use throughout this section.
\begin{notation}\label{notn:shorthand}
Let $\ba:=n\uthn$, $\bb:=n\uthn+L\uthnt$, 
\begin{align*}
    h^\ast & := \max\LP\, h(\bx-\ba) 
    \,:\,\bx\in\ZZ^2,\ltwo{\bx-\ba}\leq k\,\RP,\\
    H & := \LP\,\bx\in\ZZ^2\,:\,h(\bx-\ba)\leq h^\ast\,\RP,\\
    \tau & := \min\LP\,T(\origin,\bx)\,:\,\bx\in H\,\RP,\\
    F & := \LP\,\bx\in\ZZ^2
    \,:\,T(\origin,\bx)\leq\tau\,\RP,\\
    \partial H & := \LP\,\bx\in H
    \,:\,\bx\pm\Unit{i}\in H^c\mbox{ for some } i=1,2\,\RP.  
\end{align*}
So $H$ is the smallest $h$-ball around $\ba$ which contains 
the Euclidean ball of radius $k$ around $\ba$; $\partial H$ 
is the vertex-boundary of $H$ i.e., the set of vertices in 
$H$ which are also adjacent to some vertex outside $H$; 
$\tau$ is the time required to reach $H$ from $\origin$; $F$ 
is the set of vertices that can be reached by time $\tau$ 
from $\origin$ i.e., the wet region $\Bb(\tau)$. Furthermore,
let $\edge\subset\edge(\ZZ^2)$ be the set of nearest-neighbor
edges which have at least one endpoint in $F$, and let $\FF$
be the sigma-field generated by $\tau$, $F$, and 
$\LP \tau_e:e\in\edge\RP$.
\end{notation}

\begin{figure}[H]
    \centering
    \includegraphics[width=0.5\linewidth]{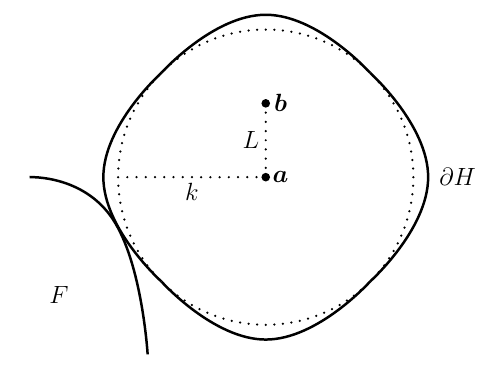}
    \caption{Illustration for $\ba$, $\bb$, $H$, $F$,
    $\partial H$: direction of $\ba$ is $\thn$, direction of 
    $\bb-\ba$ is $\thnt$. We want to prove a lower bound of 
    the variance of the transverse increment 
    $T(\origin,\ba)-T(\origin,\bb)$.}
\end{figure}

\begin{remark}\label{remark:boxh}
    Since $h$ is sublinear and $g$ is a norm, 
    we have for $\ltwo{\bu}\geq\C$, 
    \begin{equation}\label{eq:boxh1}
        \C\ltwo{\bu}\leq h(\bu)\leq\C\ltwo{\bu}.
    \end{equation}
    Therefore 
    \[
        \C k\leq h^\ast\leq\C k.
    \]
    Every $\by\in\partial H$ has an adjacent vertex that 
    does not belong to $H$. Therefore, for all 
    $\by\in\partial H$ we have
    \begin{equation}\label{eq:boxh3}
        h^\ast\geq h(\by-\ba)\geq h^\ast-\C. 
    \end{equation}
    Combining \eqref{eq:boxh1}-\eqref{eq:boxh3} we get that 
    for all $\by\in\partial H$,
    \begin{equation}\label{eq:boxh4}
        \C k\leq\ltwo{\by-\ba}\leq\C k.
    \end{equation}
    Therefore, $H$ can be inscribed in a square whose sides 
    are of the order of $k$ around $\ba$.
\end{remark}

\begin{remark} 
    Using \eqref{A2}, \eqref{eq:lowoutline0}, and 
    $L\leq\Delta(n)$, we get that for any $\delta>0$, 
    $k\leq\delta n$ for large enough $k$  depending on 
    $\delta$. Therefore, using Remark~\ref{remark:boxh}, we
    get that $\origin$ lies outside $H$. In this case, the
    region $F$ touches the region $H$ i.e., 
    $F\cap H\subset\partial H$. 
    Since the edge-weights are continuous, $F$ touches $H$ at
    only one point almost surely. So we assume that $k$ is 
    large enough such that $\origin$ lies outside $H$.
\end{remark} 

Now we state three propositions, and then we complete the 
proof of Theorem~\ref{thm:lowmain} using these propositions. 
We prove these propositions separately in later subsections. 

\begin{proposition}\label{prop:lowfast}
    Under the assumptions of Theorem~\ref{thm:lowmain}, 
    there exist constants ${\Cl{nuone}}>0$, 
    $\Cl[epsilon]{eps1}>0$, such that for large enough $k$
    \[
        \Prob\Lp\,\Prob\Lp\,T(\origin,\ba)\leq 
        h^\ast+\tau-\Cr{eps1}\sigma(k)\mid\FF\,\Rp
        \leq\Cr{eps1}\,\Rp
        \leq e^{-k^{\Cr{nuone}}}.
    \]
\end{proposition} 

By definition of $\tau$, $T(\origin,\ba)-\tau$ is an upper
bound of the time it takes for $\Gamma(\origin,\ba)$ to exit
$H$ after starting from $\ba$. Recall from \eqref{eq:boxh3}
that $h^\ast$ is approximately the expected passage time from
$\ba$ to any point on $\partial H$. Also recall from
\eqref{eq:boxh4} that points of $\partial H$ are at a
distance of the order of $k$ from $\ba$. Therefore,
Proposition~\ref{prop:lowfast} implies that, roughly
speaking, there is a nonnegligible probability that the time
taken for $\Gamma(\origin,\ba)$ to exit $H$ starting from
$\ba$ is less than $h^\ast$ by a fraction of $\sigma(k)$
i.e., the exit-time is faster than usual with nonnegligible
probability. 

\begin{proposition}\label{prop:lowslow}
    Under the assumptions of Theorem~\ref{thm:lowmain}, there
    exist constants $\Cl[nu]{nutwo}\in(1/2,\nu)$,
    $\Cl[epsilon]{eps2}>0$, such that for large enough $k$
    \[
        \Prob\Lp\,T(\origin,\ba)\geq 
        h^\ast+\tau+\Cr{eps2}\sigma(k)\,\Rp
        \geq e^{-(\log k)^{\Cr{nutwo}}}.
    \]
\end{proposition} 

Proposition~\ref{prop:lowfast} implies that, roughly
speaking, there is a nonnegligible probability that the time
taken for $\Gamma(\origin,\ba)$ to exit $H$ starting from
$\ba$ is greater than $h^\ast$ by a fraction of $\sigma(k)$
i.e., the exit-time is slower than usual with nonnegligible
probability.

\begin{proposition}\label{prop:cov}
    Under the assumptions of Theorem~\ref{thm:lowmain}, 
    there exist a constant $\Cl{pcov}>0$, such that for 
    large enough $k$   
    \[
        0\leq
        \Exp\Cov\Lp T(\origin,\ba),T(\origin,\bb)\mid\FF\Rp
        \leq\Cr{pcov}.
    \]
\end{proposition} 

Now we prove Theorem~\ref{thm:lowmain} using these
propositions. Let 
$\Cl[epsilon]{min12}:=\min\LP\Cr{eps1},\Cr{eps2}\RP$. 
Expanding the expectation of the conditional variance given 
$\FF$, we get 
\begin{align*}
           & \Var\Lp T(\origin,\ba)-T(\origin,\bb)\Rp\\
    \geq\; & \Exp\Var\Lp T(\origin,\ba)\mid\FF\Rp
    +\Exp\Var\Lp T(\origin,\bb)\mid\FF\Rp
    -2\Exp\Cov\Lp T(\origin,\ba),T(\origin,\bb)\mid\FF\Rp.
\numberthis\label{eq:lowoutline1}
\end{align*}
For any random variable $X$, 
$\Var(X)=\Exp\Lp X-X^\prime\Rp^2/2$, where $X^\prime$ is
another random variable with the same distribution as $X$ and
is independent of $X$. Therefore, for any random variable $X$
and for any $a>b$, 
\[
    \Var(X)\geq
    \frac{1}{2}(a-b)^2
    \Prob\Lp\,X\geq a\,\Rp
    \Prob\Lp\,X\leq b\,\Rp.
\]
Thus
\begin{align*}
           &  \Var\Lp T(\origin,\ba)\mid\FF\Rp\\
    \geq\; & \C \sigma^2(k)
    \Prob\Lp T(\origin,\ba)\geq h^\ast+\tau+\Cr{min12}\sigma(k)\mid\FF\Rp
    \Prob\Lp T(\origin,\ba)\leq h^\ast+\tau-\Cr{min12}\sigma(k)\mid\FF\Rp.
\numberthis\label{eq:lowoutline3}
\end{align*}
As a shorthand notation let us use 
\begin{align*}
    X & :=\Prob\Lp T(\origin,\ba)\geq
    h^\ast+\tau+\Cr{min12}\sigma(k)\mid\FF\Rp
            \Prob\Lp T(\origin,\ba)\leq 
            h^\ast+\tau-\Cr{min12}\sigma(k)\mid\FF\Rp,\\
    Y & :=\Prob\Lp T(\origin,\ba)\geq 
    h^\ast+\tau+\Cr{min12}\sigma(k)\mid\FF\Rp\Cr{min12}.
\end{align*}
Then Proposition~\ref{prop:lowslow} implies that
\begin{equation}\label{eq:lowoutline4}
    \Exp Y \geq \Cr{min12} e^{-(\log k)^{{\Cr{nutwo}}}}.
\end{equation}
Furthermore, using Proposition~\ref{prop:lowfast} and 
$0\leq X,Y\leq 1$, we get
\[
    \Exp(Y-X)^+
    \leq\Prob\Lp Y\geq X\Rp
    \leq e^{-k^{{\Cr{nuone}}}}.
\]
Therefore, using \eqref{eq:lowoutline4} and the inequality 
$\Exp X \geq \Exp Y - \Exp(Y-X)^+$, we get
\[
    \Exp X
    \geq 
    \Cr{min12} e^{-(\log k)^{\Cr{nutwo}}}
    - e^{-k^{\Cr{nuone}}}
    \geq\C e^{-(\log k)^{\Cr{nutwo}}}.
\]
Combining this with \eqref{eq:lowoutline3}, we get 
\[
    \Exp\Var\Lp T(\origin,\ba)\mid\FF\Rp
    \geq\C\sigma^2(k)e^{-(\log k)^{\Cr{nutwo}}}. 
\]
Therefore, by \eqref{eq:lowoutline1}, 
Proposition~\ref{prop:cov}, and \eqref{A2}, we get
\[
    \Var\Lp T(\origin,\ba)-T(\origin,\bb)\Rp\geq\C\sigma^2(k)e^{-(\log k)^{\Cr{nutwo}}}.
\]
Therefore, using \eqref{eq:lowoutline0}, $\nu>{\Cr{nutwo}}$, 
and \eqref{A2}, we get  
\[
    \Var\Lp T(\origin,\ba)-T(\origin,\bb)\Rp
    \geq\sigma^2\Lp\Delta^{-1}(L)\Rp e^{-(\log L)^\nu}.
\]
This completes the proof of Theorem~\ref{thm:lowmain} using 
Propositions~\ref{prop:lowfast}-\ref{prop:cov}. Now we prove 
these propositions.

\subsection{Proof of Proposition~\ref{prop:lowfast}}

\begin{figure}[H]
    \centering
    \includegraphics[width=0.5\linewidth]{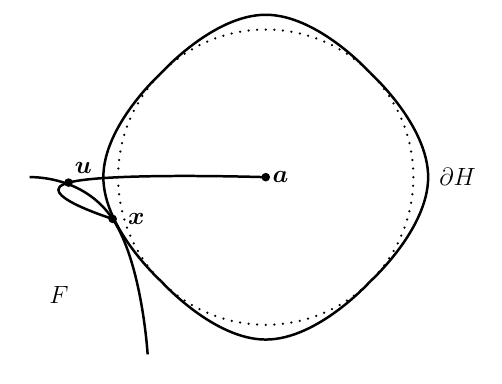}
    \caption{Setup for Proposition~\ref{prop:lowfast}: given
    a realization of the edge-weights on the whole lattice,
    we take another configuration on $\edge$, which is the
    set of edges having at least one endpoint in $F$. The
    geodesic $\Gamma^\prime(\ba,\bx)$ is then constructed in
    the environment where we have the new edge-weight
    configuration on $\edge$ and the original realization of
    edge-weights in $\edge^c$.}
    \label{fig:9}
\end{figure}

Recall the definitions of $\FF$, $\tau$, and $\edge$ from 
Notation~\ref{notn:shorthand}. Conditioned on $\FF$, consider
i.i.d.\ random variables $\LP\tau^\prime_e:e\in\edge\RP$ each
having distribution of the original edge-weights. For a path 
$\gamma$ let 
\[
    T^\prime(\gamma)=
    \sum_{\mbox{$\gamma$ contains $e$, and $e\in\edge$}}
    \tau^\prime_e +
    \sum_{\mbox{$\gamma$ contains $e$, and $e\in\edge^c$}}
    \tau_e.
\]
For any two points $\by,\bz\in\RR^2$, let 
\[
T^\prime(\by,\bz):=\inf\LP\,T^\prime(\gamma):\gamma\mbox{ is a path from $\by$ to $\bz$}\,\RP. 
\]
Therefore, the conditional distribution of all the passage
times $\{\,T^\prime(\by,\bz):\by,\bz\in\RR^2\,\}$ given $\FF$
is same as the unconditional distribution of all the passage 
times $\{\,T(\by,\bz):\by,\bz\in\RR^2\,\}$. For all 
$\by,\bz\in\RR^2$ let $\Gamma^\prime(\by,\bz)$ be the
geodesic corresponding to $T^\prime(\by,\bz)$. Let $\bx$ be 
the point where $F$ touches $H$. Let $\bu$ be the first point
belonging to $F$ when the geodesic $\Gamma^\prime(\ba,\bx)$
is traced starting from $\ba$, see Figure~\ref{fig:9}.  
Since $\bu\in F$, we have $T(\origin,\bu)\leq\tau$. Let $\bw$
be the lattice point on $\Gamma^\prime(\ba,\bx)$ preceding
$\bu$ when $\Gamma^\prime(\ba,\bx)$ is traced from $\ba$. 
Then $\Gamma^\prime(\ba,\bw)$ consists of edges only in 
$\edge^c$, and hence $T(\ba,\bw)\leq T^\prime(\ba,\bw)$. 
Therefore, 
\[
    T(\ba,\bu)
    \leq T(\ba,\bw)+T(\bw,\bu)
    \leq T^\prime(\ba,\bw)+T(\bw,\bu).
\]
For any $\bv\in\ZZ^2$ let 
\[
    d(\bv):=\max\LP\,T(\bv,\bv\pm\Unit{i}):i=1,2\,\RP.
\]
So $T(\bw,\bu)\leq d(\bu)$, and hence 
$T(\bu,\ba)\leq d(\bu) + T^\prime(\bu,\ba)$. 
Therefore,
\begin{equation}\label{eq:lowfast1}
    T(\origin,\ba) 
    \leq T(\origin,\bu) + T(\bu,\ba) 
    \leq \tau + d(\bu) + T^\prime(\bu,\ba) 
    \leq \tau + d(\bu) + T^\prime(\bx,\ba).
\end{equation}
By \eqref{A3} we get there exists $\Cl[epsilon]{fast5}>0$ 
such that
\begin{equation}\label{eq:lowfast22}
    \Prob\Lp\,
    T^\prime(\bx,\ba)\leq h(\bx-\ba) 
    -\Cr{fast5}\sigma\Lp\ltwo{\bx-\ba}\Rp
    \mid\FF\,\Rp
    \geq\Cr{fast5}.
\end{equation}
Since $F$ touches $H$ at $\bx$, we have $\bx\in\partial H$. 
Hence, $h(\bx-\ba)\leq h^\ast$ and $\ltwo{\bx-\ba}\geq k$.
Therefore, by \eqref{A2}, \eqref{eq:lowfast1}, and 
\eqref{eq:lowfast22}, we get for some $\Cl[epsilon]{fast6}>0$
\begin{align*}
    \Cr{fast6}
    \leq & \Prob\Lp\,T^\prime(\bx,\ba)\leq h^\ast-\Cr{fast6}\sigma(k)\mid\FF\,\Rp\\
    \leq & \Prob\Lp\,T(\origin,\ba) \leq h^\ast+\tau-\frac{\Cr{fast6}}{2}\sigma(k)\mid\FF\,\Rp 
        + \Prob\Lp\, d(\bu)\geq\frac{\Cr{fast6}}{2}\sigma(k)\mid\FF\,\Rp.
\end{align*}
Therefore, taking $\Cl[epsilon]{fast7}:=\Cr{fast6}/2$, we get
\begin{align*}
    & \Prob\Lp\,\Prob\Lp\,T(\origin,\ba)\leq h^\ast+\tau-\Cr{fast7}\sigma(k)\mid\FF\,\Rp
    \leq\Cr{fast7}\,\Rp\\
    \leq & \Prob\Lp\,\Prob\Lp\, d(\bu)\geq\Cr{fast7}\sigma(k)\mid\FF\,\Rp\geq\Cr{fast7}\,\Rp\\
    \leq & \Cr{fast7}^{-1}\Prob\Lp\, d(\bu)\geq\Cr{fast7}\sigma(k)\,\Rp.
    \numberthis\label{eq:lowfast2}
\end{align*}
Since $\bx\in\partial H$, by Remark~\ref{remark:boxh} we get
$\ltwo{\bx-\ba}\leq \C k$. Since $\bu$ lies on the geodesic
$\Gamma^\prime(\ba,\bx)$, by Lemma~\ref{lem:boxwand} we get
\[
    \Prob\Lp\,\ltwo{\bu-\ba}\geq\Cl{614}k\,\mid\,\FF\,\Rp
    \leq e^{-\Cl{615}k}.
\]
Therefore,
\begin{equation}\label{eq:lowfast3}
    \Prob\Lp\,\ltwo{\bu-\ba}\geq\Cr{614}k\,\Rp
    \leq e^{-\Cr{615}k}.
\end{equation}
Since edge-weights have exponential moments, for any $\bv$ 
and all $t>0$ we have
\begin{equation}\label{eq:lowfast4}
    \Prob\Lp\,d(\bv)\geq t\,\Rp\leq\C e^{-\C t}
\end{equation}
Therefore, using \eqref{eq:lowfast3} and \eqref{eq:lowfast4},
we get
\[
    \Prob\Lp\,d(\bu)\geq\Cr{fast7}\sigma(k)\,\Rp
    \leq\Cl{616} k^2 e^{-\Cl{617}\sigma(k)}. 
\]
Combining this with \eqref{eq:lowfast2} and using \eqref{A2} 
we get
\[
    \Prob\Lp\,\Prob\Lp\,T(\origin,\ba)\leq h^\ast+\tau-\Cr{fast7}\sigma(k)\,\mid\,\FF\,\Rp
    \leq \Cr{fast7}\,\Rp
    \leq \Cr{fast7}^{-1}\Cr{616} k^2 e^{-\Cr{617}\sigma(k)}
    \leq e^{-k^{\Cr{nuone}}}.
\]
This proves Proposition~\ref{prop:lowfast}.

\subsection{Proof of Proposition~\ref{prop:lowslow}}

Here is an outline of the proof. Recall from 
Notation~\ref{notn:shorthand} that $\tau$ is the passage time
from the origin to $H$, and $h^\ast$ is the maximum passage
time from $\ba$ to any point on $\partial H$. By definition,
$H$ contains an Euclidean ball of size $k$ around $\ba$, and
by Remark~\ref{remark:boxh} $H$ is contained in a ball of
radius of the order of $k$ around $\ba$. Thus, $\sigma(k)$ is
the order of fluctuation of passage times from $\ba$ to any
point on $\partial H$. Therefore, to prove
Proposition~\ref{prop:lowslow}, it suffices to show that the 
time it takes for $\Gamma(\origin,\ba)$ to reach $\partial H$
after starting from $\ba$ is slower than $h^\ast$ by a
fraction of $\sigma(k)$ with a non-negligible probability.
Here, by a non-negligible probability we mean probability at
least $e^{-(\log k)^\epsilon}$ for some $\epsilon\in(0,1)$.
To prove this, we define $H^\ast$, a subset of $\partial H$,
such that with high probability $\Gamma(\origin,\ba)$ does
not intersect $H\backslash H^\ast$. Then, we show that
passage times from $\ba$ to the points of $H^\ast$ can be
uniformly slow with non-negligible probability. To achieve
this, we further define $G^\ast$, a polygonal line which is
roughly a discrete approximation of a sector of a $g$-ball
around $\ba$. We show that passage times from $\ba$ to points
on $G^\ast$ are uniformly slow with non-negligible
probability, and passage time from $G^\ast$ to $H^\ast$ are
sufficiently small, so that passage times from $\ba$ to
points of $H^\ast$ are also slow. Moreover, $G^\ast$ is 
constructed in a way so that $\Gamma(\origin,\ba)$, when
traced from starting from $\ba$, intersect $G^\ast$ before 
$H^\ast$ with high probability. 

Let us now begin the proof formally. By 
Remark~\ref{remark:boxh}, the maximum distance in the 
direction $-\thn$ from $\ba$ of a point in $\partial H$ is at
most of the order of $k$, i.e.,
\begin{equation}\label{eq:maxH}
    \max_{\bx\in\partial H}
    \xprojthnot(\ba-\bx)\leq\Cl{slow1}k.
\end{equation}
We do not need to use any absolute value in the above 
equation because we have assumed $\thn\in[0,\pi/4]$.
Using Corollary~\ref{cor:endwandlogspl} and \eqref{A2}, we 
have for all $\bx$ on $\Gamma(\origin,\ba)$ with 
$\xprojthnot(\ba-\bx)\leq\Cr{slow1}k$ 
\[
    |\yprojthnot(\ba-\bx)|
    \leq\Cl{slow2}\Delta(k)(\log k)^{1/2},
\]
with probability at least $1-e^{-\Cl{slow3}\log k}$. 
Thus, defining the event
\begin{align*}
    \event: & \mbox{ $\Gamma(\origin,\ba)$ does not wander 
    more than $\Cr{slow2}\Delta(k)(\log k)^{1/2}$ in 
    $\pm\thnt$ directions}\\
    & \mbox{ before exiting $H$ when traced  starting from 
    $\ba$},
\end{align*}
we get 
\begin{equation}\label{eq:event0}
    \Prob\Lp\,\event\,\Rp\geq 1-e^{-\Cr{slow3}\log k}.
\end{equation}
This motivates us to define $H^\ast$ as the portion of
$\partial H$, facing towards the origin i.e., towards
direction $-\thn$ from $\ba$, having width
$2\Cr{slow2}\Delta(k)(\log k)^{1/2}$ in $\thnt$ direction,
i.e.,
\begin{equation}\label{eq:defHstar}
    H^\ast:=\LP\,\bx\in\partial H\,:\,
    |\yprojthnot(\ba-\bx)|\leq\Cr{slow2}\Delta(k)(\log k)^{1/2},\;
    \xprojthnot(\ba-\bx)\geq 0\,\RP.
\end{equation}
Thus, \eqref{eq:event0} implies that the geodesic 
$\Gamma(\origin,\ba)$ do not pass through any point in the 
set $H\backslash H^\ast$ with probability at least 
$1-e^{-\Cr{slow3}\log k}$. Now we establish a bound on the 
width of $H^\ast$.

\begin{lemma} 
    Under the assumptions of Theorem~\ref{thm:lowmain} we have for large enough $k$
    \begin{equation}\label{eq:minH}
        \min_{\bx\in\partial H}\xprojthnot(\ba-\bx)\geq\C k,
    \end{equation}
    and
    \begin{equation}\label{eq:rangeh}
        \max_{\bx\in H^\ast}\xprojthnot(\ba-\bx)
        - \min_{\bx\in H^\ast}\xprojthnot(\ba-\bx)
        \leq \Cl{21}\Delta(k)(\log k)^{1/2}.
    \end{equation}
\end{lemma}

\begin{proof} 
    Recall the definition of $g_{\ba}$ from 
    Notation~\ref{Notation:CHAP}. Consider $\bx\in H^\ast$.
    Using the fact that $g$ is a norm, and
    \eqref{eq:defHstar}, we get 
    \begin{equation}\label{eq:width3}
        g_{\ba}(\ba-\bx)
        \geq g(\ba-\bx)-\C\yprojthnot(\ba-\bx)
        \geq g(\ba-\bx)-\C\Delta(k)(\log k)^{1/2}.
    \end{equation}
    Using Proposition~\ref{prop:nrflog}, \eqref{eq:boxh3}, 
    \eqref{eq:boxh4}, and \eqref{A2}, we get
    \begin{equation}\label{eq:width4}
        g(\ba-\bx)\geq h^\ast-\C\sigma(k)\log k.
    \end{equation}
    Combining \eqref{eq:width3} and \eqref{eq:width4} we get 
    \begin{equation}\label{eq:width7}
        g_{\ba}(\ba-\bx)\geq h^\ast-\C\Delta(k)
        (\log k)^{1/2}.
    \end{equation}
    Using \eqref{eq:boxh3} and 
    Remark~\ref{remark:nonlattice1}, we get
    \begin{equation}\label{eq:width8}
             g_{\ba}(\ba-\bx)
        \leq g(\ba-\bx)
        \leq h(\ba-\bx)+\Cl{633}
        \leq h^\ast+\Cr{633}.
    \end{equation}
    Combining \eqref{eq:width7} and \eqref{eq:width8} we get 
    \[
        \max_{\bx\in H^\ast}g_{\ba}(\ba-\bx) - 
        \min_{\bx\in H^\ast}g_{\ba}(\ba-\bx)
        \leq \C\Delta(k)(\log k)^{1/2}.
    \]
    This establishes \eqref{eq:rangeh} because 
    $g_{\ba}(\ba-\bx)$
    is proportional to $\xprojthnot(\ba-\bx)$. Combining
    \eqref{eq:maxH} and \eqref{eq:rangeh} we get 
    \eqref{eq:minH}.
\end{proof}

\paragraph{\textbf{Construction and properties of $G^\ast$:}}
We denote the vertices of $G^\ast$ by 
$\Lp\bb_{i}\Rp_{i=-N_2}^{N_1}$. In addition, for each
$i$, we denote the direction of $\ba-\bb_{i}$ by $\theta_i$.
We construct $G^\ast$ in a manner such that the vertices of
$G^\ast$ have the following properties: (1) $\bb_{0}$ is
situated on the line joining $\origin$ and $\ba$; (2) the
points $\LP\bb_{i}:0<i\leq N_1\RP$ are above the line joining
$\origin$ and $\ba$; (3) the points 
$\LP\bb_{i}:-N_2\leq i<0\RP$ are below the line joining 
$\origin$ and $\ba$; (4) $\theta_i$s are in a small
neighborhood of $\thn$, say $[\thn-\delta,\thn+\delta]$, so 
that $\theta_i^t$ exists for all $i$; (5) for $0\leq i<N_1$, 
direction of $\bb_{i+1}-\bb_{i}$ is $\theta_i^t$; (6) for 
$-N_2\leq i< 0$ direction of $\bb_{i-1}-\bb_i$ is 
$-\theta_i^t$. 

Let us now proceed with the construction. Let $\Cr{nutwo}$, 
$\Cl[nu]{nu3}$ be constants such that 
\begin{equation}\label{eq:defeta}
    \frac{1}{2}<\Cr{nu3}<\Cr{nutwo}<\nu. 
\end{equation}
By Remark~\ref{remark:linear}, we choose a $\delta>0$ such 
that the limit shape boundary is differentiable in the sector
$[\thn-\delta,\thn+\delta]$, and for $\theta$ belonging to 
this sector, we have
\begin{equation}\label{eq:linear}
    |\theta^t-\thnt|\leq\Cl{tangent}|\theta-\thn|.
\end{equation}
Let 
\begin{equation}\label{eq:Gwidth}
    \ell:=
    \frac{\Delta(k)(\log k)^{1/2}}{(\log k)^{\Cr{nu3}}}. 
\end{equation} 
The point $\bb_{0}$ is defined by the conditions 
$\yprojthnot(\bb_{0})=0$ and 
\begin{equation}\label{eq:adjacent}
    \xprojthnot(\bb_{0})=\max_{y\in H^\ast}\xprojthnot(\by).
\end{equation}
We construct $\LP\bb_{i}:0<i\leq N_1\RP$ inductively as
follows. Suppose for some $j\geq 0$ we have defined 
$(\bb_{0},\dots,\bb_{j})$. Further assume that
$(\bb_{0},\dots,\bb_{j})$ satisfies the conditions:
(i) $|\theta_i-\thn|\leq\delta$ for all $0\leq i\leq j$; 
(ii) $\yprojthnot(\bb_{j})<
\Cr{slow2}\Delta(k)(\log k)^{1/2}$;
(iii) direction of $\bb_{i+1}-\bb_{i}$ is $\theta_i^t$ for 
all $0\leq i<j$. Due to convexity of $\partial\Bb$,
$\xprojthnot(\bb_{i+1}-\bb_{i})\geq 0$ and 
$\yprojthnot(\bb_{i+1}-\bb_{i})\geq 0$ for all $0\leq i<j$.
We construct $\bb_{j+1}$ as follows. Let $\bb_{j}^\prime$ be 
the point such that direction of $\bb_{j}^\prime-\bb_{j}$ is 
$\theta_j^t$, $\ltwo{\bb_{j}^\prime-\bb_{j}}=\ell$. We define
$\bb_{j+1}$ to be $\bb_{j}^\prime$, if
$\yprojthnot(\bb_{j}^\prime)
\leq\Cr{slow2}\Delta(k)(\log k)^{1/2}$. Otherwise, we take 
$\bb_{j+1}$ to be the point $\bb_{j}^{\prime\prime}$ on the 
line joining $\bb_{j}$ and $\bb_{j}^\prime$ which satisfies 
$\yprojthnot(\bb_{j}^{\prime\prime})
=\Cr{slow2}\Delta(k)(\log k)^{1/2}$, and end the 
construction. To establish that this construction is 
well-defined, it suffices to show that 
$|\theta_{j+1}-\thn|\leq\delta$, where $\theta_{j+1}$ is the 
direction of $\ba-\bb_{j+1}$. Assuming $\delta$ is small 
enough and using \eqref{eq:linear} we get for all 
$0\leq i\leq j$
\begin{equation}\label{eq:wd1}
    \xprojthnot(\bb_{i+1}-\bb_{i})
    =\ltwo{\bb_{i+1}-\bb_{i}}
    \frac{|\sin(\theta_i^t-\thnt)|}{|\sin(\thnt-\thn)|}
    \leq\C\delta\ltwo{\bb_{i+1}-\bb_{i}},
\end{equation}
and
\begin{equation}\label{eq:wd2}
     \yprojthnot(\bb_{i+1}-\bb_{i})
    =\ltwo{\bb_{i+1}-\bb_{i}}
     \frac{|\sin(\theta_i^t-\thn)|}{|\sin(\thnt-\thn)|}
    \geq 
     \C\ltwo{\bb_{i+1}-\bb_{i}}.
\end{equation}
By construction we have
\begin{equation}\label{eq:wd3}
     \yprojthnot(\bb_{j+1}-\ba)
    =\yprojthnot(\bb_{j+1}-\bb_{0})
    =\yprojthnot(\bb_{j+1})
    \leq\Cr{slow2}\Delta(k)(\log k)^{1/2}.
\end{equation}
Taking sum over $0\leq i\leq j$ in \eqref{eq:wd1} and 
\eqref{eq:wd2}, and using \eqref{eq:wd3} we get 
\begin{equation}\label{eq:wd5}
        \xprojthnot(\bb_{j+1}-\bb_{0})
    \leq\C\delta\yprojthnot(\bb_{j+1}-\bb_{0})
    \leq\C\delta\Delta(k)(\log k)^{1/2}.
\end{equation}
Therefore, using \eqref{eq:minH} and \eqref{eq:adjacent}, we get
\begin{equation}\label{eq:wd4}
\xprojthnot(\ba-\bb_{j+1})\geq\C k.
\end{equation}
Combining this with \eqref{eq:wd3} and using \eqref{A2} we 
get $|\theta_{j+1}-\thn|\leq\delta$ for large enough $k$. 
This shows that the construction is well-defined. In a 
similar way we construct $\LP\bb_{i}:-N_2\leq i<0\RP$. 
We require $\xprojthnot(\bb_{i})$ decreases from
$0$ to $-\Cr{slow2}\Delta(k)(\log k)^{1/2}$ as $i$ runs from
$0$ to $-N_2$. Equations~\eqref{eq:wd3} and \eqref{eq:wd4}
also yield for all $\bx$ in the part of $G^\ast$ joining
$\bb_{0}$ and $\bb_{N_1}$
\begin{equation}\label{eq:Glength}
    \C k \leq \ltwo{\ba-\bx} \leq \C k, 
\end{equation}
but the same holds for all $\bx$ in the part joining 
$\bb_{0}$ and $\bb_{-N_2}$. By \eqref{eq:wd2}, 
\eqref{eq:wd3}, and \eqref{eq:Gwidth} we get 
$N_1\leq\C(\log k)^{\Cr{nu3}}$, and the same is true for 
$N_2$. Hence the total number of sides of $G^\ast$ is bounded
as
\begin{equation}\label{eq:Gtotal}
    N := N_1 + N_2 \leq \C (\log k)^{\Cr{nu3}}. 
\end{equation}
By \eqref{eq:wd5}, width of the part of $G^\ast$ joining 
$\bb_{0}$ and $\bb_{N_1}$ is at most 
$\C\Delta(k)(\log k)^{1/2}$ in $\thn$ direction, and the same
holds for the part joining $\bb_{0}$ and $\bb_{-N_2}$. 
Therefore, for all $\bx,\by\in G^\ast$ 
\begin{equation}\label{eq:Gwidth2}
    |\xprojthnot(\bx-\by)|\leq\C\Delta(k)(\log k)^{1/2}
\end{equation}
This ends the discussion on construction and properties of 
$G^\ast$.

We now state two propositions and complete the proof of     
Proposition~\ref{prop:lowslow} using them. Then we proceed to
prove these propositions.

\begin{proposition}\label{prop:lowslowline}
Under the assumptions of Theorem~\ref{thm:lowmain}, there 
exists $\Cl[epsilon]{line1}>0$ such that for all large enough
$k$ we have
\[
     \Prob\Lp\,T(\bx,\ba)\geq h(\bx-\ba)+\Cr{line1}\sigma(k)\mbox{ for all } \bx\in G^\ast\,\Rp
\geq e^{-\Cl{slow4}(\log k)^{\Cr{nu3}}}.
\]
\end{proposition}

\begin{proposition}\label{prop:varlowslow3}
Under the assumptions of Theorem~\ref{thm:lowmain}, we have 
\[
     \Prob\Lp\,
     |T(\bx,\by)-h(\bx-\by)|\geq\frac{\Cr{line1}}{2}
     \sigma(k) \mbox{ for some } \bx\in H^\ast, 
     \by\in G^\ast\,\Rp
\leq e^{-k^{\Cl{slow5}}}.
\]
where $\Cr{line1}$ is the constant from 
\eqref{prop:lowslowline}.
\end{proposition}

Let us now complete proof of Proposition~\ref{prop:lowslow} using these propositions. Define the event 
\[
    \event_1: \mbox{The geodesic $\Gamma(\origin,\ba)$, when 
    traced from $\ba$ to $\origin$, intersects first $G^\ast$
    and then $H^\ast$}.
\]
Since $\event_1\subset\event$, we have using 
\eqref{eq:event0}
\begin{equation}\label{eq:slowevent1}
    \Prob\Lp\,\event_1^c\,\Rp\leq e^{-\Cr{slow3}\log k}.
\end{equation}
Define the events 
\[
    \event_2: \mbox{$T(\bx,\ba)\geq 
    h(\bx-\ba)+\epsilon_{8}\sigma(k)$ for all $\bx\in 
    G^\ast$},
\]
and 
\[
    \event_3: \;|T(\bx,\by)-h(\bx-\by)|\leq
    \frac{\epsilon_{8}}{2}\sigma(k)\mbox{ for all }
    \bx\in H^\ast,\; \by\in G^\ast.
\]
Using \eqref{eq:slowevent1}, 
Proposition~\ref{prop:lowslowline}, and 
Proposition~\ref{prop:varlowslow3}, and \eqref{eq:defeta}, we
get 
\[
    \Prob\Lp\,\event_1\cap\event_2\cap\event_3\,\Rp
    \geq e^{-(\log k)^{\Cr{nutwo}}}.
\]
So let us suppose 
$\mathbb{T}\in\event_1\cap\event_2\cap\event_3$.

\begin{figure}[H]
    \centering
    \includegraphics[width=0.5\linewidth]{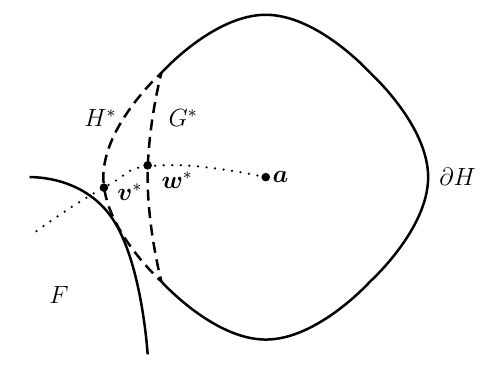}
    \caption{Setup of Proposition~\ref{prop:lowslow}: we show
    that (i) the geodesic $\Gamma(\origin,\ba)$ when traced
    starting from $\ba$ intersects first $G^\ast$ and then
    $H^\ast$ with high probability; (ii) passage time from
    $\ba$ to points in $G^\ast$ can be large with
    non-negligible probability; (iii) passage time between
    points in $G^\ast$ and points in $H^\ast$ are not too
    small with high probability. So passage time from $\ba$
    to $G^\ast$ can be large with non-negligible
    probability.}
    \label{Fig:10}
\end{figure}

Since $\mathbb{T}\in\event_1$, there exist points 
$\bw^\ast\in G^\ast$ and $\bv^\ast\in H^\ast$, both situated 
on the geodesic $\Gamma(\origin,\ba)$, such that
\begin{equation}\label{eq:slowevent11}
    T(\ba,\bv^\ast)=T(\ba,\bw^\ast)+T(\bw^\ast,\bv^\ast),
\end{equation}
see Figure~\ref{Fig:10}. Since $\bv^\ast\in 
H^\ast\subset\partial H$, we have
\begin{equation}\label{eq:slowevent12}
    T(\origin,\bv^\ast)\geq\tau.
\end{equation}
Combining \eqref{eq:slowevent11} and \eqref{eq:slowevent12} 
we get
\begin{equation}\label{eq:slowevent13}
    T(\origin,\ba) = T(\origin,\bv^\ast) 
    + T(\bv^\ast,\ba) \geq \tau + T(\ba,\bw^\ast) 
    + T(\bw^\ast,\bv^\ast).
\end{equation}
From $\mathbb{T}\in\event_2\cap\event_3$ we get
\begin{multline}
    T(\ba,\bw^\ast)+T(\bw^\ast,\bv^\ast)
    \geq h(\bw^\ast-\ba) + \epsilon_{8}\sigma(k) + 
    h(\bv^\ast-\bw^\ast)
    - \frac{\epsilon_{8}}{2}\sigma(k)
    \\ 
    \geq h(\bv^\ast-\ba)
    + \frac{\epsilon_{8}}{2}\sigma(k)
    - \Cl{623},
\label{eq:slowevent14}
\end{multline}
where we get the extra constant at the end by 
Remark~\ref{remark:nonlattice1}. Since $\bv^\ast$ is on 
$\partial H$, from \eqref{eq:boxh3} we have 
\begin{equation}\label{eq:slowevent14.5}
    h(\bv^\ast-\ba)\geq h^\ast-\Cl{624}.
\end{equation} 
For large enough $k$ we have 
\begin{equation}\label{eq:slowevent15}
    \frac{\epsilon_{8}}{2}\sigma(k)-\Cr{623}-\Cr{624}
    \geq\frac{\epsilon_{8}}{4}\sigma(k). 
\end{equation} 
Therefore, combining 
\eqref{eq:slowevent13}-\eqref{eq:slowevent15} we get 
\[ 
    T(\origin,\ba)\geq
    \tau+h^\ast+\frac{\epsilon_{8}}{4}\sigma(k).
\]
Letting $\epsilon_2:=\epsilon_{8}/4$ completes proof of
Proposition~\ref{prop:lowslow}. Let us now prove
Propositions~\ref{prop:lowslowline} and 
\ref{prop:varlowslow3}.  

\subsubsection{Proof of Proposition~\ref{prop:lowslowline}}

Let $G_i$ be the segment of $G^\ast$ joining $\bb_{i}$ and
$\bb_{i+1}$ for all $-N_2\leq i<N_1$. For all $i$, let
$\ba_{i}$ be the point on $G_i$ which has the maximum 
expected passage time from $\ba$. By \eqref{A3}, there exists
$\Cl[epsilon]{line2}>0$ such that for each $i$, 
\[
    \Prob\Lp\,T(\ba,\ba_{i})\geq h(\ba_{i}-\ba) + 
    \Cr{line2}\sigma(\ltwo{\ba_{i}-\ba})\,\Rp\geq\Cr{line2}.
\]
Therefore, using \eqref{eq:Glength} and \eqref{A2} we get, 
for some $\Cl[epsilon]{line3}>0$ and for all $i$,
\begin{equation}\label{eq:pointtoline1}
    \Prob\Lp\,T(\ba,\ba_{i})\geq h(\ba_{i}-\ba) + 
    \Cr{line3}\sigma(k)\,\Rp\geq\Cr{line3}.
\end{equation}
Define for each $i$
\[
    D_i:=\max\LP\,|T(\ba,\bx)-T(\ba,\by)|\,:\,
    \bx,\by\in G_i\,\RP.
\]
Recall $\Cr{nu3}>1/2$ from \eqref{eq:defeta}. Choose 
$\Cl[nu]{nu4}>0$ such that 
\[
    \Cr{nu4}<\frac{2\alpha}{(1+\beta)}(\Cr{nu3}-1/2).
\]
To bound $D_i$ for $i<0$ we can use 
Theorem~\ref{thm:loglogupinc} with the variables
\begin{align*}
    \tilde{\eta}:=\Cr{nu4},\quad
    \tilde{\thn}:=-\theta_i,\quad
    \tilde{L}:=\ell,\quad
    \tilde{n}:=\ltwo{\bb_{i}-\ba}. 
\end{align*}
For $i\geq 0$ we have a minor technical issue because of the 
direction of the tangents, but same bounds hold. Using 
\eqref{eq:Gwidth}, \eqref{eq:Glength}, and \eqref{A2}, we get
$\tilde{L}\geq\tilde{L_0}$, $\tilde{n}\geq\tilde{n_0}$,
$\tilde{L}\leq\Delta(\tilde{n})$, as required. By 
Theorem~\ref{thm:loglogupinc} and \eqref{A2}, we get for
large enough $t$
\begin{equation}\label{eq:pointtoline2}
    \Prob\Lp\,D_i\geq t(\log k)^{\Cr{nu4}}\sigma(\Delta^{-1}(\ell))\,\Rp
    \leq\C \exp\Lp-\C t (\log k)^{\Cr{nu4}}\Rp.
\end{equation}
Using \eqref{eq:Gwidth} and \eqref{A2}, we get
\begin{align*}
    \frac{\sigma(k)}{
    \sigma(\Delta^{-1}(\ell))(\log k)^{\Cr{nu4}}}
    \geq 
    \C(\log k)^{-\Cr{nu4}+(\Cr{nu3}-1/2)(2\alpha)/(1+\beta)}.
\end{align*}
This can be made arbitrarily large by choosing $k$ large. 
Therefore, in \eqref{eq:pointtoline2} we can choose 
\[
    t=\frac{\Cr{line3}}{2}\frac{\sigma(k)}
        {\sigma(\Delta^{-1}(\ell))(\log k)^{\Cr{nu4}}},
\]
and we get for large enough $k$
\[
    \Prob\Lp\,D_i\geq\frac{\Cr{line3}}{2} 
    \sigma(k)\,\Rp\leq\frac{\Cr{line3}}{2}.
\]
Combining this with \eqref{eq:pointtoline1} we get for all 
$i$
\[
    \Prob\Lp\,
    T(\ba,\bx)\geq h(\bx-\ba)+\frac{\Cr{line3}}{2}\sigma(k)
    \mbox{ for all }\bx\in G_i\,\Rp
    \geq\frac{\Cr{line3}}{2}.
\]
Since for each $i$, 
$\inf_{\bx\in G_i}\LP\,T(\ba,\bx)-h(\bx-\ba)\,\RP$ is an 
increasing function of the edge-weight configuration, by the 
FKG inequality we get
\[
    \Prob\Lp\,
    T(\ba,\bx)\geq h(\bx-\ba)+\frac{\Cr{line3}}{2}\sigma(k)
    \mbox{ for all }\bx\in G_i
    \,\Rp
    \geq\Lp\frac{\Cr{line3}}{2}\Rp^N,
\]
where recall from \eqref{eq:Gtotal} that $N$ is the total 
number of segments in $G^\ast$. Using \eqref{eq:Gtotal} we 
get 
\[
    \Prob\Lp\,T(\ba,\bx)\geq h(\bx-\ba)+\frac{\Cr{line3}}{2}\sigma(k)
    \mbox{ for all }\bx\in G_i\,\Rp
    \geq\exp\Lp-\C (\log k)^{\Cr{nu3}}\Rp.
\]
This completes the proof of 
Proposition~\ref{prop:lowslowline}.

\subsubsection{Proof of Proposition~\ref{prop:varlowslow3}}

By \eqref{eq:rangeh}, width of $H^\ast$ in $\thn$ direction
is $2\Cr{slow2}\Delta(k)(\log k)^{1/2}$. By 
\eqref{eq:Gwidth2}, width of $G^\ast$ in $\thn$ direction is 
$\C \Delta(k)(\log k)^{1/2}$. By construction of $G^\ast$ we 
have 
\[
    \min_{\bx\in G^\ast}\xprojthnot(\bx)
    =\max_{\by\in H^\ast}\xprojthnot(\by).
\]
Both $G^\ast$ and $H^\ast$ are centered around the line
joining $\origin$ and $\ba$, and both go up to distance
$\Cr{slow2}\Delta(k)(\log k)^{1/2}$ in $\pm\thnt$ directions.
Therefore, for all $\bx\in G^\ast$ and $\by\in H^\ast$
\[
    \ltwo{\bx-\by}\leq\C\Delta(k)(\log k)^{1/2}.
\]
Using \eqref{A2}, the number of pairs of such points is at
most $\C k^2$. Hence, using \eqref{A1}, \eqref{A2}, and a 
union bound, we get 
\begin{align*}
     & \Prob\Lp\;|T(\bx,\by)-h(\bx-\by)|\geq\frac{\Cr{line1}}{2}\sigma(k)\mbox{ for some }
     \bx\in G^\ast\mbox{ and }\by\in H^\ast\;\Rp\\
\leq & \C k^2 \exp\Lp-\C
        \frac{\sigma(k)\log k}
        {\sigma(\Delta(k)(\log k)^{1/2})}\Rp\\
\leq & \exp\Lp-k^{\C}\Rp.
\end{align*}
This completes the proof of 
Proposition~\ref{prop:varlowslow3}.

\subsection{Proof of Proposition~\ref{prop:cov}}

The passage times $T(\origin,\ba)$ and $T(\origin,\bb)$ 
are increasing functions of the edge-weight configuration 
$\mathbb{T}$. Therefore, using the FKG inequality and taking 
expectation we get 
\[
    \Exp\Cov\Lp T(\origin,\ba),T(\origin,\bb)\mid\FF\Rp
    \geq 0.
\]
Define two collections of paths 
\begin{align*}
    \Pi(\origin,\ba):=&
    \left\{\vphantom{|\yprojthnot(\by-\ba)|<\frac{1}{2}
    \Delta(k)(\log k)^\eta}
    \;\gamma\,:\,\mbox{$\gamma$ is a path from $\origin$ to
    $\ba$, and for all lattice points}\right.\\
    &\left.\;\mbox{ $\by$ in $\gamma$ outside $F$ we have 
    $|\yprojthnot(\by-\ba)|<
    \frac{1}{2}\Delta(k)(\log k)^\eta$}\;
    \vphantom{|\yprojthnot(\by-\ba)|<\frac{1}{2}
    \Delta(k)(\log k)^\eta}\right\},\\
    \Pi(\origin,\bb):=&
    \left\{\vphantom{|\yprojthnot(\by-\bb)|<\frac{1}{2}
    \Delta(k)(\log k)^\eta}
    \;\gamma\,:\,\mbox{$\gamma$ is a path from $\origin$ to
    $\bb$, and for all lattice points}\right.\\
    &\left.\;\mbox{ $\by$ in $\gamma$ outside $F$ we have 
    $|\yprojthnot(\by-\bb)|<
    \frac{1}{2}\Delta(k)(\log k)^\eta$}\;
    \vphantom{|\yprojthnot(\by-\bb)|<\frac{1}{2}
    \Delta(k)(\log k)^\eta}\right\}.
\end{align*}
Since $|\yprojthnot(\ba-\bb)|=L=\Delta(k)(\log k)^\eta$, 
a path in $\Pi(\origin,\ba)$ do not touch a path in 
$\Pi(\origin,\bb)$ outside $F$. Let 
\[
    T^\prime(\origin,\ba):=
    \min_{\gamma\in\Pi(\origin,\ba)}T(\gamma),\quad
    T^\prime(\origin,\bb):=
    \min_{\gamma\in\Pi(\origin,\bb)}T(\gamma).
\]
Since paths in $\Pi(\origin,\ba)$ do not intersect paths in 
$\Pi(\origin,\bb)$ outside $F$, $T^\prime(\origin,\ba)$ and 
$T^\prime(\origin,\bb)$ are independent conditioned on $\FF$.

Let $\Gamma(\ba,F)$ be the geodesic from $\ba$ to $F$ i.e.,
the path with minimum passage time from $\ba$ to a point in 
$F$. Let $T(\ba,F):=T(\Gamma(\ba,F))$. Since $F$ touches $H$,
we have 
\[
    T(\ba,F)\leq\max_{\by\in H}T(\ba,\by).
\]
By Remark~\ref{remark:boxh}, $H$ is contained in a square of 
side length $\Cl{644} k$. Therefore, by 
Lemma~\ref{lem:boxwand} we get that there exist positive 
constants $\Cl{645}$ and $\Cl{645a}$ such that
\begin{equation}\label{eq:cov1}
    \Prob\Lp\,T(\ba,F)\geq 
    \Cr{645} k\,\Rp\leq e^{-\Cr{645a} k}.
\end{equation}
Let 
\[
    T^t(\origin,\ba):=
    \min\LP T^\prime(\origin,\ba),\tau+\Cr{645}k\RP,\quad 
    T^t(\origin,\bb):=
    \min\LP T^\prime(\origin,\bb),\tau+\Cr{645}k\RP.
\]
Because $T^\prime(\origin,\ba)$ and $T^\prime(\origin,\bb)$ 
are independent conditioned on $\FF$, $T^t(\origin,\ba)$ and 
$T^t(\origin,\bb)$ are also independent conditioned on $\FF$.
Therefore
\begin{align*}
    & \Exp\Cov\Lp T(\origin,\ba),T(\origin,\bb)|\FF\Rp\\
    = & \Exp\Cov\Lp T(\origin,\ba)-T^t(\origin,\ba),
    T(\origin,\bb)-\tau\mid\FF\Rp\\
    & +\Exp\Cov\Lp T^t(\origin,\ba)-\tau,
    T(\origin,\bb)-T^t(\origin,\bb)\mid\FF\Rp\\
    \leq & \Exp\LT\Lp\Exp\Lp\Lp 
    T(\origin,\ba)-T^t(\origin,\ba)\Rp^2\mid\FF\Rp\Rp^{1/2}
    \Lp\Exp\Lp\Lp 
    T(\origin,\bb)-\tau\Rp^2\mid\FF\Rp\Rp^{1/2}\RT\\
    & +\Exp\LT\Lp\Exp\Lp\Lp 
    T(\origin,\bb)-T^t(\origin,\bb)\Rp^2\mid\FF\Rp\Rp^{1/2}
      	\Lp\Exp\Lp\Lp 
      	 T^t(\origin,\ba)-\tau\Rp^2\mid\FF\Rp\Rp^{1/2}\RT\\
    \leq & \Lp\Exp\Lp 
    T(\origin,\ba)-T^t(\origin,\ba)\Rp^2\Rp^{1/2}
    	   \Lp\Exp\Lp T(\origin,\bb)-\tau\Rp^2\Rp^{1/2}\\
         & + \Lp\Exp\Lp 
         T(\origin,\bb)-T^t(\origin,\bb)\Rp^2\Rp^{1/2}
      	   \Lp\Exp\Lp T^t(\origin,\ba)-\tau\Rp^2\Rp^{1/2}.  
    \numberthis\label{eq:cov2}
\end{align*}
Let $\bx$ be the point where $F$ touches $H$. Therefore
$T(\origin,\bx)=\tau\leq T(\origin,\ba)$. Hence $0\leq
T(\origin,\ba)-\tau\leq T(\ba,\bx)$. Using that $H$ is
contained in a box of size $\Cr{644}k$ we get
\begin{equation}\label{eq:cov21}
    \Exp\Lp T(\origin,\ba)-\tau\Rp^4
    \leq\Exp T(\ba,\bx)^4\leq\C k^4.
\end{equation}
Since $T^\prime(\origin,\ba)$ is the minimum passage time 
restricted to some paths from $\origin$ to $\ba$, we have
$T^\prime(\origin,\ba)\geq T(\origin,\ba)$. Therefore 
$T^t(\origin,\ba)\geq \tau$. Therefore
\begin{equation}\label{eq:cov22}
     |T^t(\origin,\ba)-\tau|\leq \C k.
\end{equation}
Combining \eqref{eq:cov21} and \eqref{eq:cov22} we get
\begin{align*}
         & \Lp\Exp\Lp 
         T(\origin,\ba)-T^t(\origin,\ba)\Rp^2\Rp^{1/2}\\
    \leq & \Lp\Exp\Lp 
    T(\origin,\ba)-T^t(\origin,\ba)\Rp^4\Rp^{1/4}
    \Prob\Lp 
    T(\origin,\ba)\neq T^t(\origin,\ba)\Rp^{1/4}\\
    \leq & \LT\Lp\Exp\Lp T(\origin,\ba)-\tau\Rp^4\Rp^{1/4}+
        \Lp\Exp\Lp 
        T^t(\origin,\ba)-\tau\Rp^4\Rp^{1/4}\RT\Prob\Lp 
        T(\origin,\ba)\neq T^t(\origin,\ba)\Rp^{1/4}\\
    \leq & 
    \C k \Prob\Lp T(\origin,\ba)\neq 
    T^t(\origin,\ba)\Rp^{1/4}.
\numberthis\label{eq:cov3}
\end{align*}
Similarly we have
\begin{equation}\label{eq:cov3.1}
    \Exp\Lp T(\origin,\bb)-\tau\Rp^2
\leq\C k^2.
\end{equation}
and
\begin{equation}\label{eq:cov3.2}
    \Lp\Exp\Lp T(\origin,\bb)-T^t(\origin,\bb)\Rp^2\Rp^{1/2}
    \leq \C k \Prob\Lp T(\origin,\bb)\neq 
    T^t(\origin,\bb)\Rp^{1/4}.
\end{equation}
Therefore, combining \eqref{eq:cov2}, \eqref{eq:cov3}-\eqref{eq:cov3.2} we 
get
\begin{multline}\label{eq:covprefinal}
     \Exp\Cov\Lp T(\origin,\ba),T(\origin,\bb)|\FF\Rp
    \\ \leq\C k^2\Lp\Prob\Lp T(\origin,\ba)\neq T^t(\origin,\ba)\Rp^{1/4}+\Prob\Lp T(\origin,\bb)\neq T^t(\origin,\bb)\Rp^{1/4}\Rp. 
\end{multline}
Let us consider the term 
$\Prob\Lp T(\origin,\ba)\neq T^t(\origin,\ba)\Rp$, the other 
one can be dealt similarly. Depending on $\Cr{645}$ there 
exist positive constants $\Cl{646}$ and $\Cl{647}$ such that 
for large enough $k$ we have
\begin{equation}\label{eq:cov4}
    \Prob\Lp\,\min\LP\, T(\ba,\by):\ltwo{\by-\ba}=\Cr{646}k\,\RP\leq\Cr{645} k\,\Rp\leq e^{-\Cr{647} k}.
\end{equation}
Define the events
\begin{align*}
    \event_1 & : T(\ba,F)\leq\Cr{645}k,\\
    \event_2 & : \min\LP T(\ba,\by):\ltwo{\by-\ba}
    =\Cr{646} k\RP>\Cr{645}k\\
    \event_3 & : \text{Diam}(\Gamma(\ba,F))\leq 2\Cr{646}k.
\end{align*}
Therefore, $\event_1$ and $\event_2$ implies $\event_3$. 
Recall Notation~\ref{notn:wandering}. Let
\[
    \event_4 : \max_{k^\prime\leq 2\Cr{646}k} 
    \wandering{\ba}{\origin}{k^\prime}{-\thn}
    <\frac{1}{2}\Delta(k)(\log k)^\eta.
\]
Using Corollary~\ref{cor:endwandlogspl} and $\eta>1/2$, we 
get
\begin{equation}\label{eq:cov5}
    \Prob\Lp\,\event_4^c\,\Rp
    \leq e^{-\C (\log k)^{2\eta}}.
\end{equation}
If $\mathbb{T}\in\event_1\cap\event_3\cap\event_4$ then 
$T(\origin,\ba)=T^t(\origin,\ba)$. Therefore, combining 
\eqref{eq:cov1}, \eqref{eq:cov4}, and \eqref{eq:cov5}, we get
\[
    \Prob\Lp\,T(\origin,\ba)\neq T^t(\origin,a)\,\Rp
    \leq e^{-\C(\log k)^{2\eta}}.
\]
Similar bound holds 
$\Prob\Lp\,T(\origin,\bb)\neq T^t(\origin,\bb)\,\Rp$. 
Therefore, using \eqref{eq:covprefinal} and $\eta>1/2$, we 
get 
\[
    \Exp\Cov(T(\origin,\ba),T(\origin,\bb)|\FF)
    \leq\Cl{pcov2}
\]
for any $\Cr{pcov2}>0$, provided $k$ is large enough. This 
completes the proof of Proposition~\ref{prop:cov}.

\section{Upper bound of the long-range correlations}

\resetconstant

In this section, our objective is to prove 
Theorem~\ref{thm:longrange}.
Fix $J_0>0$, $n_0>0$ to be assumed large enough whenever 
required. 
Consider $n\geq n_0$, $J\in[\upconst^{1/2}J_0,n^\delta]$.  
Let $m:=f^{-1}(Jf(n)/J_0)$, so that 
\begin{equation}\label{eq:long1}
     J_0\frac{\Delta(m)(\log m)^{1/2}}{m}
    =J\frac{\Delta(n)(\log n)^{1/2}}{n}.
\end{equation}
Using $J\geq\upconst^{1/2}J_0$ and \eqref{A2} we get 
$m\leq n$. Using $\delta<(1-\beta)/2$, $J\leq n^\delta$, 
\eqref{A2}, and assuming $n_0$ is large enough, we get 
$\log m\geq\Cl{long1}\log n$. As a shorthand notation let us 
use 
\[
\begin{aligned}
    \ba:=n\uthn+J\Delta(n)(\log n)^{1/2}\uthnt,\quad
    \bb:=n\uthn-J\Delta(n)(\log n)^{1/2}\uthnt.
\end{aligned}
\]
Using $\delta<(1-\beta)/2$ and $J\leq n^\delta$, we get 
$\ltwo{\ba}$ and $\ltwo{\bb}$ are at most $\Cl{long2}n$. Let 
\[
    H:=\LP\,\bx\in\RR^2:\xprojthnot(\bx)\geq m\,\RP. 
\]
Let $\FF$ be the sigma-field generated by all the 
edge-weights $\tau_e$ such that both endpoints of the edge 
$e$ are in $H$. Our objective is to establish upper bound of 
$\Cov(T(\origin,\ba),T(\origin,\bb))$. We split the 
covariance in expectation of conditional covariances given 
$\FF$ and covariance of conditional expectations given $\FF$,
and establish upper bound of them separately.

\begin{proposition}\label{prop:longrange1}
    Assuming $J_0$ and $n_0$ are large enough, we have
    \begin{equation}\label{eq:longprop1.1}
        |\Cov\Lp\Exp\Lp T(\origin,\ba)\mid\FF\Rp,\Exp\Lp 
        T(\origin,\bb)\mid\FF\Rp\Rp|\leq\C.
    \end{equation}
\end{proposition}

\begin{figure}[H]
    \centering
    \includegraphics[width=0.65\linewidth]{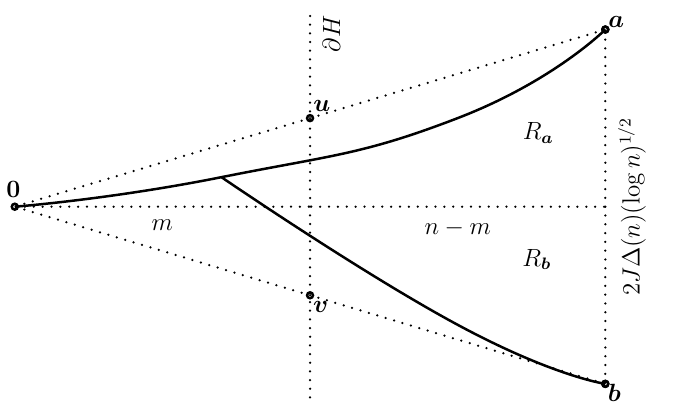}
    \caption{Setup of Proposition~\ref{prop:longrange1}: 
    $\ltwo{\bu-\bv}=2J\Delta(m)(\log m)^{1/2}$; 
    $H$ is the region to the right of $\partial 
    H$; $R_{\ba}$ is the subset of $H$ above the line in 
    direction $\thn$, $R_{\bb}$ is the region below the line;
    with high probability $\Gamma(\origin,\ba)$ stays in 
    $R_{\ba}$ while it is in $H$, $\Gamma(\origin,\bb)$ stays
    in $R_{\bb}$ while it is in $H$.}
    \label{Fig:11}
\end{figure}%
 
\begin{proof}
    Define two regions
    \[
    \begin{aligned}
        R_{\ba}:=\LP\,\bx\in\RR^2:\xprojthnot(\bx)\geq 
        m,\;\yprojthnot(\bx)>0\,\RP,\\
        R_{\bb}:=\LP\,\bx\in\RR^2:\xprojthnot(\bx)\geq 
        m,\;\yprojthnot(\bx)<0\,\RP.
        \end{aligned}
    \]
    Define the event
    \[
        \event_1: \mbox{ $\Gamma(\origin,\ba)$ stays inside 
        $R_{\ba}$ while it is in the region $H$,}
    \]
    i.e., for all $\bu\in\Gamma(\origin,\ba)$ with 
    $\xprojthnot(\bu)\geq m$ we have $\yprojthnot(\bu)>0$.
    Similarly, define the event 
    \[
        \event_2: \mbox{ $\Gamma(\origin,\bb)$ stays in the
        region $R_{\bb}$ while it is in the region $H$.} 
    \]
    Thus, if $\mathbb{T}\not\in\event_1$, then 
    either
    $\wandering{\origin}{\ba}{k}{\thn}\geq 
    J_0(k/m)\Delta(m)(\log m)^{1/2}$ 
    for some $k\in[m,n]$,  
    or 
    $\wandering{\origin}{\ba}{k}{\thn}\geq
    J\Delta(n)(\log n)^{1/2}$
    for some $k\geq n$.
    By \eqref{A2}, 
    \[
    J_0\frac{k}{m}\Delta(m)(\log m)^{1/2}\geq 
    J_0\upconst^{-1/2}\Delta(k)(\log k)^{1/2}
    \]
    for $k\geq m$.
    Therefore, for each $k\in[m,n]$, using 
    Proposition~\ref{prop:nrflog} and 
    Theorem~\ref{thm:endwandlog}, we get 
    \begin{equation}\label{eq:longprop1.2}
        \Prob\Lp\,\wandering{\origin}{\ba}{k}{\thn}\geq 
        J_0\frac{k}{m}\Delta(m)(\log m)^{1/2}\,\Rp
        \leq\C e^{-\C J_0^2\log k}.
    \end{equation}
    Using Corollary~\ref{cor:endwandlogspl} we get
    \begin{equation}\label{eq:longprop1.3}
        \Prob\Lp\,\wandering{\origin}{\ba}{k}{\thn}\geq 
        J_0\frac{k}{m}\Delta(m)(\log m)^{1/2}
        \mbox{ for some $k\geq n$}\,\Rp
        \leq \C e^{- \C J_0^2 \log n}.
    \end{equation}
    Using union bound over $k\in[m,n]$ in \eqref{eq:longprop1.2}, 
    $\log m\geq\Cr{long1}\log n$, and \eqref{eq:longprop1.3}, we 
    get
    \begin{equation}\label{eq:longprop1.4}
        \Prob\Lp\,\event_1^c\,\Rp\leq \C e^{-\C J_0^2\log n}.
    \end{equation} 
    The same holds for $\event_2$.
    
    Let $\hat{T}(\origin,\ba)$ be the minimum passage time 
    among of all paths from $\origin$ to $\ba$ which stays in
    $R_{\ba}$ when in $H$. Similarly define 
    $\hat{T}(\origin,\bb)$. 
    Then $\Exp(\hat{T}(\origin,\ba)|\FF)$ and 
    $\Exp(\hat{T}(\origin,\bb)|\FF)$ are independent because 
    $R_{\ba}$ and $R_{\bb}$ are disjoint. 
    If $\mathbb{T}\in\event_1$, then 
    $T(\origin,\ba)=\hat{T}(\origin,\ba)$, and if 
    $\mathbb{T}\in\event_2$, then 
    $T(\origin,\bb)=\hat{T}(\origin,\bb)$.
    Using $\ltwo{\ba}\leq\Cr{long2}n$, 
    $\ltwo{\bb}\leq\Cr{long2}n$, \eqref{eq:longprop1.4}, and the 
    same bound for $\event_2$, we get 
    \begin{align*}
             & \Cov\Lp\Exp\Lp T(\origin,\ba)\mid\FF\Rp,
        \Exp\Lp T(\origin,\bb)\mid\FF\Rp\Rp\\ 
        \leq & \Lp\Exp
                \Lp T(\origin,\ba)-\hat{T}(\origin,\ba)\Rp^2
               \Rp^{1/2}
               \Lp\Exp T(\origin,\bb)^2\Rp^{1/2}\\
            & +\Lp\Exp
                \Lp T(\origin,\bb)-\hat{T}(\origin,\bb)\Rp^2
               \Rp^{1/2}
               \Lp\Exp T(\origin,\ba)^2\Rp^{1/2}\\
        \leq & \Prob(\,\event_1^c\,)^{1/2}
               \Lp\Exp\Lp 
               T(\origin,\ba)-\hat{T}(\origin,\ba)
               \Rp^4\Rp^{1/4}
               \Lp\Exp T(\origin,\bb)^2\Rp^{1/2}\\
             &+\Prob(\,\event_2^c\,)^{1/2}
               \Lp\Exp\Lp
               T(\origin,\bb)-\hat{T}(\origin,\bb)
               \Rp^4\Rp^{1/4}
               \Lp\Exp T(\origin,\bb)^2\Rp^{1/2}\\
        \leq & \C n^2 e^{-\C J_0^2 \log n}\\
        \leq & \C. 
    \end{align*}
    This completes the proof of 
    Proposition~\ref{prop:longrange1}.
\end{proof}

Now we consider the expected conditional covariance. We have
\begin{equation}\label{eq:long2}
    \Exp\Cov\Lp T(\origin,\ba),T(\origin,\bb)|\FF\Rp
    \leq\Lp\Exp\Var\Lp T(\origin,\ba)|\FF\Rp\Rp^{1/2}
    \Lp\Exp\Var\Lp T(\origin,\bb)|\FF\Rp\Rp^{1/2}.
\end{equation}
We establish an upper bound of 
$\Exp\Var(T(\origin,\ba)|\FF)$. Similar bound also 
holds for $\Exp\Var(T(\origin,\ba)|\FF)$. Consider an 
independent edge-weight configuration on the edges which have
at least one endpoint not in the half-space $H$. Let
$T^\prime(\origin,\ba)$ be the passage time from $\origin$ to
$\ba$ on the new configuration. Let 
$\Gamma^\prime(\origin,\ba)$ be the corresponding geodesic. 
Then $T(\origin,\ba)$ and $T^\prime(\origin,\ba)$ are 
independent given $\FF$. Therefore
\begin{equation}\label{eq:long3}
     \Exp\Var(T(\origin,\ba)|\FF)
    =\Exp\Lp
        \Exp\Lp
            \frac{1}{2}
            \Lp T(\origin,\ba)-T^\prime(\origin,\ba)\Rp^2|\FF
        \Rp
    \Rp
    =\frac{1}{2}
     \Exp\Lp T(\origin,\ba)-T^\prime(\origin,\ba)\Rp^2.
\end{equation}

\begin{proposition}\label{prop:longrange2}
    Assuming $J_0$ and $n_0$ are large enough, we have
    \begin{equation}\label{eq:longprop2.1}
        \Exp(T(\origin,\ba)-T^\prime(\origin,\ba))^2
        \leq\Cl{longrangeclaim2}\sigma^2(m)\log n.
    \end{equation}
\end{proposition}

\begin{figure}[H]
    \centering
    \includegraphics[width=0.65\linewidth]{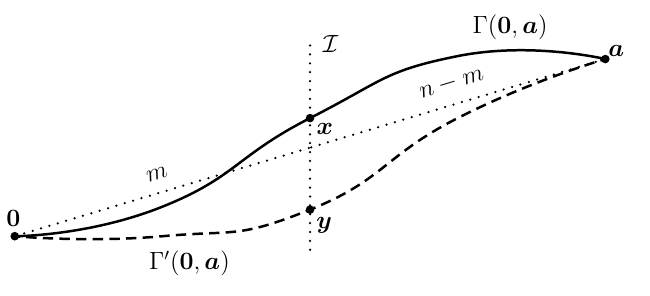}
    \caption{Setup of Proposition~\ref{prop:longrange2}: the segment $\segment$ is a part of $\partial H$, see
    Figure~\ref{Fig:11} for the location of $\partial H$; the
    geodesic $\Gamma^\prime(\origin,\ba)$ is constructed by
    taking a new configuration in the left-side of 
    $\partial H$. With high probability, both geodesics 
    $\Gamma(\origin,\ba)$ and $\Gamma^\prime(\origin,\ba)$
    intersect $\segment$ when they intersect $\partial H$.}
\end{figure}

\begin{proof} 
    Consider the line segment
    \[
    \segment:=\LP\,\bx:\xprojthnot(\bx)=m,\, 0\leq\yprojthnot(\bx)\leq 2J_0\Delta(m)(\log m)^{1/2}\,\RP. 
    \]
    Define the event $\event_3$
    \[
    \event_3: \mbox{ $\Gamma(\origin,\ba)$ and $\Gamma^\prime(\origin,\ba)$ passes through $\segment$.}
    \]
    If $\mathbb{T}\not\in\event_3$ then both $\Gamma(\origin,\ba)$ and $\Gamma^\prime(\origin,\ba)$ wanders more than $J_0\Delta(m)(\log m)^{1/2}$ in $\pm\thnt$ directions when they are at distance $m$ in $\thn$ direction from $\origin$. 
    So, using Theorem~\ref{thm:midptwand} and $\log m\geq\Cr{long1}\log n$, we get 
    \[
    \Prob\Lp\,\event_3^c\,\Rp\leq \C e^{-\C J_0^2\log n}.
    \]
    Therefore, using $\ltwo{\ba}\leq\Cr{long2}n$ and assuming $J_0$ is large enough, we get  
    \begin{equation}\label{eq:longprop2.2}
    \Exp\Lp\Lp T(\origin,\ba)-T^\prime(\origin,\ba)\Rp^2
    \mathbb{1}(\event_3^c)\Rp\leq\C.
    \end{equation}
    If we have $\mathbb{T}\in\event_3$, $\Gamma(\origin,\ba)$ passes through $\bx\in\segment$, and $\Gamma^\prime(\origin,\ba)$ passes through $\by\in\segment$, then
    \[
    T(\origin,\bx)-T^\prime(0,\bx)
    \leq T(\origin,\ba)-T^\prime(\origin,\ba)
    \leq T(\origin,\by)-T^\prime(\origin,\by).
    \]
    Therefore
    \begin{equation}\label{eq:longprop2.3}
    \Exp\Lp\Lp T(\origin,\ba)-T^\prime(\origin,\ba)\Rp^2\mathbb{1}(\event_3)\Rp
    \leq\Exp\max_{\bz\in\segment}\Lp T(\origin,\bz)-T^\prime(\origin,\bz)\Rp^2.
    \end{equation}
    For every $\bz\in\segment$, $T(\origin,\bz)$ and $T^\prime(\origin,\bz)$ have the same mean. Using \eqref{eq:long1}, $J\leq n^\delta$, $\delta\leq(1-\beta)/2$, and \eqref{A2}, we get $\ltwo{\bz}\leq\C m$ for all $\bz\in\segment$. 
    Therefore, by \eqref{A1} and \eqref{A2}, for all $t>0$
    \[ 
    \Prob\Lp\,\max_{\bz\in\segment} |T(\origin,\bz)-T^\prime(\origin,\bz)|\geq t\sigma(m)\,\Rp\leq \C m e^{-\C t}.
    \]
    Therefore, 
    \[
    \Exp\max_{\bz\in\segment}\Lp T(\origin,\bz)-T^\prime(\origin,\bz)\Rp^2
    \leq\C\sigma^2(m)\log m.
    \]
    Combining this with \eqref{eq:longprop2.2}, \eqref{eq:longprop2.3}, and using $m\leq n$ proves Proposition~\ref{prop:longrange2}. 
\end{proof}

Combining \eqref{eq:long3} and \eqref{eq:longprop2.1} we get
\[
    \Exp\Var\Lp T(\origin,\ba)\,\mid\,\FF\Rp
    \leq\C\sigma^2(m)\log n.
\]
By symmetry, the same statement holds if we replace $\ba$ by $\bb$. Therefore by \eqref{eq:long2} we get
\begin{equation}\label{eq:long4}
    \Exp\Cov\Lp T(\origin,\ba),T(\origin,\bb)\,\mid\,\FF\Rp
    \leq\C\sigma^2(m)\log n. 
\end{equation}
Therefore, by \eqref{A2}, the bound on the covariance of the 
conditional expectations in \eqref{eq:longprop1.1} is 
negligible compared to the bound on the expectation of the 
conditional covariance in \eqref{eq:long4}. Thus, combining 
\eqref{eq:longprop1.1} and \eqref{eq:long4}, proves 
Theorem~\ref{thm:longrange}.

\end{document}